\documentclass{amsart}%
\usepackage{hyperref}
\usepackage{amssymb}
\usepackage{amsmath}
\usepackage{amsthm}
\usepackage{amsfonts}
\usepackage{graphicx}%
\allowdisplaybreaks
\makeatletter
\@ifundefined{pdfoutput}
{
\DeclareGraphicsRule{.wmf}{bmp}{}{}
\DeclareGraphicsRule{.jpg}{bmp}{}{}
\DeclareGraphicsRule{.png}{bmp}{}{}
\DeclareGraphicsRule{.cdr}{bmp}{}{}
\DeclareGraphicsRule{.gif}{bmp}{}{}
}
{
\ifnum\pdfoutput=0\relax
\DeclareGraphicsRule{.wmf}{bmp}{}{}
\DeclareGraphicsRule{.jpg}{bmp}{}{}
\DeclareGraphicsRule{.png}{bmp}{}{}
\DeclareGraphicsRule{.cdr}{bmp}{}{}
\DeclareGraphicsRule{.gif}{bmp}{}{}
\fi
\ifnum\pdfoutput=1\relax
\DeclareGraphicsRule{.eps}{pdf}{}{}
\DeclareGraphicsRule{.wmf}{jpg}{}{}
\fi
}
\makeatother
\newtheorem{thm}{Theorem}[section]
\newtheorem{cor}[thm]{Corollary}
\newtheorem{ex}[thm]{Example}

\newtheorem{lem}[thm]{Lemma}
\newtheorem{nota}[thm]{Notation}
\newtheorem{prop}[thm]{Proposition}

\theoremstyle{definition}
\newtheorem*{acknowledgement}{Acknowledgement}
\newtheorem{df}[thm]{Definition}

\newtheorem{rem}[thm]{Remark}

\numberwithin{equation}{section}
\begin{document}
\def\H{\mathcal{H}}
\def\P{\mathcal{P}}
\def\WC{W_{\mathbb C}}
\def\HC{H_{\mathbb C}}

\def\enddoc{\end{document}}

\title[Heat Kernels on Infinite-dimensional Heisneberg Groups]{Heat Kernel Analysis on Infinite-Dimensional Heisenberg Groups}
\author[Driver]{Bruce K. Driver{$^{\dagger }$}}
\thanks{\footnotemark {$^\dagger$}This research was supported in part by NSF Grant DMS-0504608 and the Miller Institute at the University of California, at Berkeley.}
\address{Department of Mathematics, 0112\\
University of California, San Diego \\
La Jolla, CA 92093-0112 } \email{driver{@}euclid.ucsd.edu}
\author[Gordina]{Maria Gordina{$^{*}$}}
\thanks{\footnotemark {$*$}Research was supported in part by NSF Grant DMS-0706784.}
\address{Department of Mathematics\\
University of Connecticut\\
Storrs, CT 06269, U.S.A. } \email{gordina@math.uconn.edu}


% 35K05 Heat equation (Primary)
% 58G32 Diffusion processes and stochastic analysis on manifolds (Primary)
% 43A15 Lp spaces and other function spaces on groups. (Secondary)

\date{\today \ \emph{File:\jobname{.tex}}}

\begin{abstract}
We introduce a class of non-commutative Heisenberg like infinite dimensional
Lie groups based on an abstract Wiener space. The Ricci curvature tensor for these groups is computed and shown to be bounded. Brownian motion and the corresponding heat kernel measures, $\{\nu_t\}_{t>0},$ are also studied. We show that these heat kernel measures admit: 1) Gaussian like upper bounds, 2) Cameron-Martin type quasi-invariance results, 3) good $L^p$ -- bounds on the corresponding Radon-Nykodim derivatives, 4) integration by parts formulas, and 5) logarithmic Sobolev inequalities. The last three results heavily rely on the boundedness of the Ricci tensor.

\end{abstract}

\maketitle

\tableofcontents

\renewcommand{\contentsname}{Table of Contents}

\section{Introduction\label{s.1}}

Both authors have been greatly influenced by Professor Malliavin and his work
over the years. In particular this paper is partially an attempt to better
understand Malliavin's paper, \cite{Malliavin1990}. It is with great pleasure
to us that this article appears (assuming it is accepted) in this special
edition of JFA dedicated to Professor Paul Malliavin.

The aim of this paper is to construct and study properties of heat kernel
measures on certain infinite-dimensional Heisenberg groups. In this paper the
Heisenberg groups will be constructed from a skew symmetric form on an
abstract Wiener space. A typical example of such a group is the Heisenberg
group of a symplectic vector space. Before describing our results let us
recall some typical heat kernel results for finite-dimensional Riemannian manifolds.

\subsection{A finite-dimensional paradigm\label{s.1.1}}

Let $\left(  M,g\right)  $ be a complete connected $n$ -- dimensional
Riemannian manifold $\left(  n<\infty\right)  ,$ $\Delta=\Delta_{g}$ be the
Laplace Beltrami operator acting on $C^{2}\left(  M\right)  ,$ and
$\operatorname{Ric}$ denote the associated Ricci tensor. Recall (see for
example Strichartz \cite{Str1983}, Dodziuk \cite{Dodziuk83} and Davies
\cite{Davies89}) that the closure, $\bar{\Delta},$ of $\Delta|_{C_{c}^{\infty
}(M)}$ is self-adjoint on $L^{2}(M,dV),$ where $dV=\sqrt{g}dx^{1}\dots dx^{n}$
is the Riemann volume measure on $M.$ Moreover, the semi-group $P_{t}%
:=e^{t\bar{\Delta}/2}$ has a symmetric positive integral (heat) kernel,
$p_{t}(x,y),$ such that $\int_{M}p_{t}(x,y)dV(y)\leqslant1$ for all $x\in M$
and%
\begin{equation}
P_{t}f(x):=\left(  e^{t\bar{\Delta}/2}f\right)  (x)=\int_{M}p_{t}%
(x,y)f(y)dV(y)\text{ for all }f\in L^{2}(M). \label{e.1.1}%
\end{equation}
Theorem \ref{t.1.2} summarizes some of the results that we would like to
extend to our infinite-dimensional Heisenberg group setting.

\begin{nota}
\label{n.1.1}If $\mu$ is a probability measure on a measure space $\left(
\Omega,\mathcal{F}\right)  $ and $f\in L^{1}\left(  \mu\right)  =L^{1}\left(
\Omega,\mathcal{F},\mu\right)  ,$ we will often write $\mu(f)$ for the
integral, $\int_{\Omega}fd\mu.$
\end{nota}

\begin{thm}
\label{t.1.2}Beyond the assumptions above, let us further assume that
$\operatorname{Ric}\geqslant kI$ for some $k\in\mathbb{R}.$ Then

\begin{enumerate}
\item $p_{t}\left(  x,y\right)  $ is a smooth function. (The Ricci curvature
assumption is not needed here.)

\item $\int_{M}p_{t}\left(  x,y\right)  dV\left(  y\right)  =1,$ (see for
example Davies \cite[Theorem 5.2.6 ]{Davies89}).

\item Given a point $o\in M,$ let $d\nu_{t}\left(  x\right)  :=p_{t}\left(
o,x\right)  dV\left(  x\right)  $ for all $t>0.$ Then $\left\{  \nu
_{t}\right\}  _{t>0}$ may be characterized as the unique family of probability
measures such that the function $t\rightarrow\nu_{t}(f):=\int_{M}fd\nu_{t}$ is
continuously differentiable,
\begin{equation}
\frac{d}{dt}\nu_{t}(f)=\frac{1}{2}\nu_{t}(\Delta f),\text{ and }%
\lim_{t\downarrow0}\nu_{t}(f)=f(o) \label{e.1.2}%
\end{equation}
for all $f\in BC^{2}(M)$, the bounded $C^{2}$--functions on $M.$

\item There exist constants, $c=c\left(  K, n , T\right)  $ and $C=C\left(  K,
n, T\right)  ,$ such that,%

\begin{equation}
p(t,x,y)\leqslant\frac{C}{V\left(  x,\sqrt{t/2}\right)  }\exp\left(
-c\frac{d^{2}(x,y)}{t}\right)  , \label{e.1.3}%
\end{equation}
for all $x,y\in M$ and $t\in(0,T],$ where $d\left(  x,y\right)  $ is the
Riemannian distance from $x$ to $y$ and $V\left(  x,r\right)  $ is the volume
of the $r$ -- ball centered at $x.$

\item The \textbf{heat kernel measure, }$\nu_{T},$ for any $T>0$ satisfies the
following logarithmic Sobolev inequality;
\begin{equation}
\nu_{T}(f^{2}\log f^{2})\leqslant2k^{-1}\left(  1-e^{-kT}\right)  v_{T}\left(
|\nabla f|^{2}\right)  +\nu_{T}(f^{2})\log\nu_{T}\left(  f^{2}\right)  ,
\label{e.1.4}%
\end{equation}
for $f\in C_{c}^{\infty}\left(  M\right)  .$
\end{enumerate}
\end{thm}

These results are fairly standard. For item 3. see \cite[Theorem 2.6
]{Driver2003b}, for Eq. (\ref{e.1.3}) see for example Theorems 5.6.4, 5.6.6,
and 5.4.12 in Sallof-Coste \cite{Saloff-Coste2002} and for more detailed
bounds see \cite{Li-Yau86,Davies89,Saloff-Coste92,Davies93,Grigoryan1995}).
The logarithmic Sobolev inequality in Eq. (\ref{e.1.4}) generalizes Gross'
\cite{Gross75} original Logarithmic Sobolev inequality valid for
$M=\mathbb{R}^{n}$ and is due in this generality to D. Bakry and M. Ledoux,
see \cite{BaLe1,BaLe2,Ledoux2000}. Also see
\cite{Hsu95,Hsu97b,Wang96a,Wang97b,Driver-Hu96} and Driver and Lohrenz
\cite[Theorem 2.9]{Driver1996b} for the case of interest here, namely when $M$
is a uni-modular Lie group with a left invariant Riemannian metric.

When passing to infinite-dimensional Riemannian manifolds we will no longer
have available the Riemannian volume measure. Because of this problem, we will
take item 3. of Theorem \ref{t.1.2} as our definition of the heat kernel
measure. The heat kernel upper bound in Eq. (\ref{e.1.3}) also does not make
sense in infinite dimensions. However, the following consequence almost does:
there exists $c\left(  T\right)  >0$ such that
\begin{equation}
\int_{M}\exp\left(  c\frac{d^{2}(o,x)}{t}\right)  d\nu_{t}\left(  x\right)
<\infty\text{ for all }0<t\leq T. \label{e.1.5}%
\end{equation}
In fact Eq. (\ref{e.1.5}) will not hold in infinite dimensions either. It will
be necessary to replace the distance function, $d,$ by a weaker distance
function as happens in Fernique's theorem for Gaussian measure spaces. With
these results as background we are now ready to summarize the results of this paper.

\subsection{Summary of results\label{s.1.2}}

Let us describe the setting informally, for precise definitions see Sections
\ref{s.2} and \ref{s.3}. Let $\left(  W,H,\mu\right)  $ be an abstract Wiener
space, $\mathbf{C}$ be a finite-dimensional inner product space, and
$\omega:W\times W\rightarrow\mathbf{C}$ be a continuous skew symmetric
bilinear quadratic form on $W.$ The set $\mathfrak{g}=W\times\mathbf{C}$ can
be equipped with a Lie bracket by setting%

\[
\lbrack\left(  A,a\right)  ,\left(  B,b\right)  ]=\left(  0,\omega\left(
A,B\right)  \right)  .
\]
As in the case for the Heisenberg group of a symplectic vector space, the Lie
algebra $\mathfrak{g}=W\times\mathbf{C}$ can be given the group structure by
defining
\[
\left(  w_{1},c_{1}\right)  \cdot\left(  w_{2},c_{2}\right)  =\left(
w_{1}+w_{2},c_{1}+c_{2}+\frac{1}{2}\omega\left(  \omega_{1},w_{2}\right)
\right)  .
\]
The set $W\times\mathbf{C}$ with the group structure will be denoted by $G$ or
$G\left(  \omega\right)  $. The Lie subalgebra $\mathfrak{g}_{CM}%
=H\times\mathbf{C}$ is called the Cameron-Martin subalgebra, and
$\mathfrak{g}_{CM}$ equipped with the same group multiplication denoted by
$G_{CM}$ and called the Cameron-Martin subgroup. We equip $G_{CM}$ with the
left invariant Riemannian metric which agrees with the natural Hilbert inner
product,
\[
\left\langle \left(  A,a\right)  ,\left(  B,b\right)  \right\rangle
_{\mathfrak{g}_{CM}}:=\left\langle A,B\right\rangle _{H}+\left\langle
a,b\right\rangle _{\mathbf{C}},
\]
on $\mathfrak{g}_{CM}\cong T_{\mathbf{e}}G_{CM}.$ In Section \ref{s.3} we give
several examples of this abstract setting including the standard
finite-dimensional Heisenberg group.

The main objects of our study are a Brownian motion in $G$ and the
corresponding heat kernel measure defined in Section \ref{s.4}. Namely, let
$\left\{  \left(  B\left(  t\right)  ,B_{0}\left(  t\right)  \right)
\right\}  _{t\geqslant0}$ be a Brownian motion on $\mathfrak{g}$ with variance
determined by
\begin{multline*}
\mathbb{E}\left[  \left\langle \left(  B\left(  s\right)  ,B_{0}\left(
s\right)  \right)  , \left(  A,a\right)  \right\rangle _{\mathfrak{g}_{CM}%
}\cdot\left\langle \left(  B\left(  t\right)  ,B_{0}\left(  t\right)  \right)
,\left(  C,c\right)  \right\rangle _{\mathfrak{g}_{CM}}\right] \\
=\operatorname{Re}\left\langle \left(  A, a\right)  ,\left(  C,c\right)
\right\rangle _{\mathfrak{g}_{CM}}\min\left(  s,t\right)
\end{multline*}
for all $s,t\in\lbrack0,\infty),$ $A,C\in H_{\ast}$ and $a,c\in\mathbf{C}.$
Then the Brownian motion on $G$ is the continuous $G$--valued process defined
by%
\[
g\left(  t\right)  =\left(  B\left(  t\right)  ,B_{0}\left(  t\right)
+\frac{1}{2}\int_{0}^{t}\omega\left(  B\left(  \tau\right)  ,dB\left(
\tau\right)  \right)  \right)  .
\]
For $T>0$ the heat kernel measure on $G$ is $\nu_{T}=\operatorname{Law}\left(
g\left(  T\right)  \right)  .$ It is shown in Corollary \ref{c.4.5} that
$\left\{  \nu_{t}\right\}  _{t>0}$ satisfies item 3. of Theorem \ref{t.1.2}
with $o=\left(  0,0\right)  \in G\left(  \omega\right)  .$

Theorem \ref{t.4.16} gives heat kernel measure bounds that may be viewed as a
non-commutative version of Fernique's theorem for $G(\omega)$. In light of
Theorem \ref{t.3.12} this result is also analogous to the integrated
integrated Gaussian upper bound in Eq. (\ref{e.1.5}).

In Theorem \ref{t.5.2} we prove quasi-invariance for the path space measure
associated to the Brownian motion, $g,$ on $G$ with respect to multiplication
on the left by finite energy paths in the Cameron-Martin subgroup $G_{CM}$.
(In light of the results in Malliavin \cite{Malliavin1990} it is surprising
that Theorem \ref{t.5.2} holds.) Theorem \ref{t.5.2} is then used to prove
quasi-invariance of the heat kernel measures with respect to both right and
left multiplication (Theorem \ref{t.6.1} and Corollary \ref{c.6.2}), as well
as integration by parts formulae on the path space and for the heat kernel
measures, see Corollaries \ref{c.5.6} -- \ref{c.6.5}. These results can be
interpreted as the first steps towards proving $\nu_{t}$ is a
\textquotedblleft strictly positive\textquotedblright\ smooth measure. In this
infinite-dimensional it is natural to interpret quasi-invariance and
integration by parts formulae as properties as smoothness of the heat kernel
measure, see \cite[Theorem 3.3]{Driver2003b} for example.

In Section \ref{s.7} we compute the Ricci curvature and check that not only it
is bounded from below (see Proposition \ref{p.7.2}), but also that the Ricci
curvature of certain finite-dimensional \textquotedblleft
approximations\textquotedblright\ are bounded from below with constants
independent of the approximation. Based on results in \cite{DG07a}, these
bounds allow us to give another proof of the quasi-invariance result for
$\nu_{t}$ and at the same time to get $L^{p}$-- estimates on the corresponding
Radon-Nikodym derivatives, see Theorem \ref{t.8.1}. These estimates are
crucial for the heat kernel analysis on the spaces of holomorphic functions
which is the subject of our paper \cite{DG07c}. In Theorem \ref{t.8.3} we show
that analogue of the logarithmic Sobolev inequality in Eq. (\ref{e.1.4}) holds
in our setting as well.

In Section \ref{s.9} we give a list of open questions and further possible
developments of the results of this paper. We expect our methods to be
applicable to a much larger class of infinite-dimensional nilpotent groups.

Finally, we refer to papers of H.~Airault, P.~Malliavin, D.~Bell, Y.~Inahama
concerning quasi-invariance, integration by parts formulae and the logarithmic
Sobolev inequality on certain infinite-dimensional curved spaces
(\cite{AirMall2006, Bell05, Bell06a, Bell06b, Inahama04}).

\section{Abstract Wiener Space Preliminaries\label{s.2}}

Suppose that $X$ is a real separable Banach space and $\mathcal{B}_{X}$ is the
Borel $\sigma$--algebra on $X$.

\begin{df}
\label{d.2.1} A measure $\mu$ on $(X,\mathcal{B}_{X})$ is called a (mean zero,
non-degenerate) \textbf{Gaussian measure} provided that its characteristic
functional is given by
\begin{equation}
\hat{\mu}(u):=\int_{X}e^{iu\left(  x\right)  }d\mu\left(  x\right)
=e^{-\frac{1}{2}q(u,u)}\text{ for all }u\in X^{\ast}, \label{e.2.1}%
\end{equation}
where $q=q_{\mu}:X^{\ast}\times X^{\ast}\rightarrow\mathbb{R}$ is a quadratic
form such that $q(u,v)=q(v,u)$ and $q(u)=q(u,u)\geqslant0$ with equality iff
$\ u=0,$ i.e. $q$ is a real inner product on $X^{\ast}.$
\end{df}

In what follows we frequently make use of the fact that%
\begin{equation}
C_{p}:=\int_{X}\left\Vert x\right\Vert _{X}^{p}d\mu\left(  x\right)
<\infty\text{ for all }1\leqslant p<\infty. \label{e.2.2}%
\end{equation}
This is a consequence of Skorohod's inequality (see for example \cite[Theorem
3.2]{Kuo75})
\begin{equation}
\int_{X}e^{\lambda\left\Vert x\right\Vert _{X}}d\mu\left(  x\right)
<\infty\text{ for all }\lambda<\infty; \label{e.2.3}%
\end{equation}
or the even stronger Fernique's inequality (see for example \cite[Theorem
2.8.5]{Bog98} or \cite[Theorem 3.1]{Kuo75})
\begin{equation}
\int_{X}e^{\delta\left\Vert x\right\Vert _{X}^{2}}d\mu\left(  x\right)
<\infty\text{ for some }\delta>0. \label{e.2.4}%
\end{equation}

\begin{lem}
\label{l.2.2}If $u,v\in X^{\ast},$ then
\begin{equation}
\int_{X}u\left(  x\right)  v\left(  x\right)  d\mu\left(  x\right)  =q\left(
u,v\right)  \label{e.2.5}%
\end{equation}
and%
\begin{equation}
\left\vert q(u,v)\right\vert \leqslant C_{2}\Vert u\Vert_{X^{\ast}}\Vert
v\Vert_{X^{\ast}}. \label{e.2.6}%
\end{equation}

\end{lem}

\begin{proof}
Let $u_{\ast}\mu:=\mu\circ u^{-1}$ denote the law of $u$ under $\mu.$ Then by
Equation \eqref{e.2.1},%
\[
\left(  u_{\ast}\mu\right)  \left(  dx\right)  =\frac{1}{\sqrt{2\pi q\left(
u,u\right)  }}e^{-\frac{1}{2q\left(  u,u\right)  }x^{2}}dx
\]
and hence,%
\begin{equation}
\int_{X}u^{2}\left(  x\right)  d\mu\left(  x\right)  =q_{\mu}\left(
u,u\right)  =q\left(  u,u\right)  . \label{e.2.7}%
\end{equation}
Polarizing this identity gives Equation \eqref{e.2.5} which along with
Equation \eqref{e.2.2} implies Equation \eqref{e.2.6}.
\end{proof}

\iffalse

\textbf{Alternate Proof of the Existence of }$C_{2}<\infty$ such that Equation
\eqref{e.2.6} holds. Note that $q(u,u)=-2\ln[\hat{\mu}(u)]$, and so if
$u_{n}\in X^{\ast}$ and $u_{n}\rightarrow0$ weakly, that is, $u_{n}%
(x)\rightarrow0$ as $n\rightarrow\infty$ for all $x\in X$, then by the
dominated convergence theorem, $q(u_{n})\rightarrow-2\ln[\hat{\mu}(0)]=0.$
More generally, if $u_{n},v_{n}\rightarrow0$ weakly, then $|q(u_{n}%
,v_{n})|\leqslant\sqrt{q(u_{n})}\sqrt{q(v_{n})}\rightarrow0$ as $n\rightarrow
\infty$. In particular, there exists $C<\infty$ such that $\left\vert
q(u,v)\right\vert \leqslant C\Vert u\Vert_{X^{\ast}}\Vert v\Vert_{X^{\ast}}.$
Indeed, if no such $C<\infty$ existed, then there would exist $u_{n},v_{n}$
such that $\Vert u_{n}\Vert_{X^{\ast}}=\Vert v_{n}\Vert_{X^{\ast}}=1$ while
$\left\vert q(u_{n},v_{n})\right\vert \geqslant n^{2}.$ We would then have
$\left\vert q(u_{n}/n,v_{n}/n)\right\vert \geqslant1$ while $u_{n}%
/n\rightarrow0$ and, $v_{n}/n\rightarrow0.$

\fi

The next theorem summarizes some well known properties of Gaussian measures
that we will use freely below.

\begin{thm}
\label{t.2.3}Let $\mu$ be a Gaussian measure on a real separable Banach space,
$X.$ For $x\in X$ let
\begin{equation}
\left\Vert x\right\Vert _{H}:=\sup\limits_{u\in X^{\ast}\setminus\left\{
0\right\}  }\frac{|u(x)|}{\sqrt{q(u,u)}} \label{e.2.8}%
\end{equation}
and define the \textbf{Cameron-Martin subspace, }$H\subset X,$ by
\begin{equation}
H=\left\{  h\in X:\left\Vert h\right\Vert _{H}<\infty\right\}  . \label{e.2.9}%
\end{equation}
Then

\begin{enumerate}
\item $H$ is a dense subspace of $X$;

\item there exists a unique inner product, $\left\langle \cdot,\cdot
\right\rangle _{H}$ on $H$ such that $\left\Vert h\right\Vert _{H}%
^{2}=\left\langle h,h\right\rangle $ for all $h\in H.$ Moreover, with this
inner product $H$ is a separable Hilbert space.

\item For any $h\in H$
\begin{equation}
\left\Vert h\right\Vert _{X}\leqslant\sqrt{C_{2}}\left\Vert h\right\Vert _{H},
\label{e.2.10}%
\end{equation}
where $C_{2}$ is as in \eqref{e.2.2}.

\item If $\left\{  e_{j}\right\}  _{j=1}^{\infty}$ is an orthonormal basis for
$H,$ then for any $u, v \in H^{\ast}$
\begin{equation}
q\left(  u,v\right)  =\left\langle u, v\right\rangle _{H^{\ast}}=\sum
_{j=1}^{\infty}u\left(  e_{j}\right)  v\left(  e_{j}\right)  . \label{e.2.11}%
\end{equation}

\end{enumerate}
\end{thm}

The proof of this standard theorem is relegated to Appendix \ref{s.10} -- see
Theorem \ref{t.10.1}.

\begin{rem}
\label{r.2.4} It follows from Equation \eqref{e.2.10} that any $u\in X^{\ast}$
restricted to $H$ is in $H^{\ast}.$ Therefore we have
\begin{equation}
\int_{X}u^{2}d\mu=q\left(  u,u\right)  =\left\Vert u\right\Vert _{H^{\ast}%
}^{2}=\sum_{j=1}^{\infty}\left\vert u\left(  e_{j}\right)  \right\vert ^{2},
\label{e.2.12}%
\end{equation}
where $\left\{  e_{j}\right\}  _{j=1}^{\infty}$ is an orthonormal bases for
$H.$ More generally, if $\varphi$ is a linear bounded map from $W$ to
$\mathbf{C}$, where $\mathbf{C}$ is a real Hilbert space, then%
\begin{equation}
\left\Vert \varphi\right\Vert _{H^{\ast}\otimes\mathbf{C}}^{2}=:\sum
_{j=1}^{\infty}\left\Vert \varphi\left(  e_{j}\right)  \right\Vert
_{\mathbf{C}}^{2}=\int_{X}\left\Vert \varphi\left(  x\right)  \right\Vert
_{\mathbf{C}}^{2}d\mu\left(  x\right)  <\infty. \label{e.2.13}%
\end{equation}

\end{rem}

To prove Equation \eqref{e.2.13}, let $\left\{  f_{j}\right\}  _{j=1}^{\infty
}$ be an orthonormal basis for $\mathbf{C}$. Then
\begin{align*}
\int_{X}\left\Vert \varphi\left(  x\right)  \right\Vert _{\mathbf{C}}^{2}%
d\mu\left(  x\right)   &  =\int_{X}\sum_{j=1}^{\infty}\left\vert \left\langle
\varphi\left(  x\right)  ,f_{j}\right\rangle _{\mathbf{C}}\right\vert ^{2}%
d\mu\left(  x\right)  =\sum_{j=1}^{\infty}\int_{X}\left\vert \left\langle
\varphi\left(  x\right)  ,f_{j}\right\rangle _{\mathbf{C}}\right\vert ^{2}%
d\mu\left(  x\right) \\
&  =\sum_{j=1}^{\infty}\left\Vert \left\langle \varphi\left(  \cdot\right)
,f_{j}\right\rangle _{\mathbf{C}}\right\Vert _{H^{\ast}}^{2}=\sum
_{j=1}^{\infty}\sum_{k=1}^{\infty}\left\vert \left\langle \varphi\left(
e_{k}\right)  ,f_{j}\right\rangle _{\mathbf{C}}\right\vert ^{2}\\
&  =\sum_{k=1}^{\infty}\sum_{j=1}^{\infty}\left\vert \left\langle
\varphi\left(  e_{k}\right)  ,f_{j}\right\rangle _{\mathbf{C}}\right\vert
^{2}=\sum_{k=1}^{\infty}\left\Vert \varphi\left(  e_{k}\right)  \right\Vert
_{\mathbf{C}}^{2}=\left\Vert \varphi\right\Vert _{H^{\ast}\otimes\mathbf{C}%
}^{2}.
\end{align*}
A simple consequence of Eq. (\ref{e.2.14}) is that%
\begin{equation}
\left\Vert \varphi\right\Vert _{H^{\ast}\otimes\mathbf{C}}^{2}\leqslant
\left\Vert \varphi\right\Vert _{W^{\ast}\otimes\mathbf{C}}^{2}\int
_{X}\left\Vert x\right\Vert _{W}^{2}d\mu\left(  x\right)  =C_{2}\left\Vert
\varphi\right\Vert _{W^{\ast}\otimes\mathbf{C}}^{2}. \label{e.2.14}%
\end{equation}

\section{Infinite-Dimensional Heisenberg Type Groups\label{s.3}}

Throughout the rest of this paper $\left(  X,H,\mu\right)  $ will denote a
\textbf{real abstract Wiener space}, i.e. $X$ is a real separable Banach
space, $H$ is a real separable Hilbert space densely embedded into $X$, and
$\mu$ is a Gaussian measure on $\left(  X,\mathcal{B}_{X}\right)  $ such that
Equation \eqref{e.2.1} holds with $q\left(  u,u\right)  :=\left\langle
u|_{H},u|_{H}\right\rangle _{H^{\ast}}.$

Following the discussion in \cite{Kaplan1980} and \cite{Eberlein2004} we will
say that a (possibly infinite--dimensional) Lie algebra, $\mathfrak{g},$ is of
\textbf{Heisenberg type} if $\mathbf{C}:=\left[  \mathfrak{g,g}\right]  $ is
contained in the center of $\mathfrak{g}.$ If $\mathfrak{g}$ is of Heisenberg
type and $W$ is a complementary subspace to $\mathbf{C}$ in $\mathfrak{g},$ we
may define a bilinear map, $\omega:W\times W\rightarrow\mathbf{C},$ by
$\omega\left(  w,w^{\prime}\right)  =\left[  w,w^{\prime}\right]  $ for all
$w,w^{\prime}\in\mathbf{C}.$ Then for $\xi_{i}:=w_{i}+c_{i}\in W\oplus
\mathbf{C}=\mathfrak{g}$, $i=1,2$ we have
\[
\left[  \xi_{1},\xi_{2}\right]  =\left[  w_{1}+c_{1},w_{2}+c_{2}\right]
=0+\omega\left(  w_{1},w_{2}\right)  .
\]
If we now suppose $G$ is a finite-dimensional Lie group with Lie algebra
$\mathfrak{g},$ then by the Baker-Campbell-Dynkin-Hausdorff formula
\[
e^{\xi_{1}}e^{\xi_{2}}=e^{\xi_{1}+\xi_{2}+\frac{1}{2}\left[  \xi_{1},\xi
_{2}\right]  }=e^{w_{1}+w_{2}+c_{1}+c_{2}+\frac{1}{2}\omega\left(  w_{1}%
,w_{2}\right)  }.
\]
In particular, we may introduce a group structure on $\mathfrak{g}$ by
defining%
\[
\left(  w_{1}+c_{1}\right)  \cdot\left(  w_{2}+c_{2}\right)  =w_{1}%
+w_{2}+c_{1}+c_{2}+\frac{1}{2}\omega\left(  \omega_{1},w_{2}\right)  .
\]
With this as motivation, we are now going to introduce a class of Heisenberg
type Lie groups based on the following data.

\begin{nota}
\label{n.3.1}Let $\left(  W,H,\mu\right)  $ be an abstract Wiener space,
$\mathbf{C}$ be a finite-dimensional inner product space, and $\omega:W\times
W\rightarrow\mathbf{C}$ be a continuous skew symmetric bilinear quadratic form
on $W.$ Further let
\begin{equation}
\left\Vert \omega\right\Vert _{0}:=\sup\left\{  \left\Vert \omega\left(
w_{1},w_{2}\right)  \right\Vert _{\mathbf{C}}:w_{1},w_{2}\in W\text{ with
}\left\Vert w_{1}\right\Vert _{W}=\left\Vert w_{2}\right\Vert _{W}=1\right\}
. \label{e.3.1}%
\end{equation}
be the uniform norm on $\omega$ which is finite by the assumed continuity of
$\omega.$
\end{nota}

We now define $\mathfrak{g}:=W\times\mathbf{C}$ which is a Banach space in the
norm
\begin{equation}
\left\Vert \left(  w,c\right)  \right\Vert _{\mathfrak{g}}:=\left\Vert
w\right\Vert _{W}+\left\Vert c\right\Vert _{\mathbf{C}}. \label{e.3.2}%
\end{equation}
We further define $\mathfrak{g}_{CM}:=H\times\mathbf{C}$ which is a Hilbert
space relative to the product inner product
\begin{equation}
\left\langle \left(  A,a\right)  ,\left(  B,b\right)  \right\rangle
_{\mathfrak{g}_{CM}}:=\left\langle A,B\right\rangle _{H}+\left\langle
a,b\right\rangle _{\mathbf{C}}. \label{e.3.3}%
\end{equation}
The associated Hilbertian norm on $\mathfrak{g}_{CM}$ is given by
\begin{equation}
\left\Vert \left(  A,a\right)  \right\Vert _{\mathfrak{g}_{CM}}:=\sqrt
{\left\Vert A\right\Vert _{H}^{2}+\left\Vert a\right\Vert _{\mathbf{C}}^{2}}.
\label{e.3.4}%
\end{equation}
It is easily checked that defining%
\begin{equation}
\left[  \left(  w_{1},c_{1}\right)  ,\left(  w_{2},c_{2}\right)  \right]
:=\left(  0,\omega\left(  w_{1},w_{2}\right)  \right)  \label{e.3.5}%
\end{equation}
for all $\left(  w_{1},c_{1}\right)  ,\left(  w_{2},c_{2}\right)
\in\mathfrak{g}$ makes $\mathfrak{g}$ into a Lie algebra such that
$\mathfrak{g}_{CM}$ is Lie subalgebra of $\mathfrak{g}.$ Note that this
definition implies that $\mathbf{C}=[\mathfrak{g},\mathfrak{g}]$ is contained
in the center of $\mathfrak{g}$. It is also easy to verify that we may make
$\mathfrak{g}$ into a group using the multiplication rule
\begin{equation}
\left(  w_{1},c_{1}\right)  \cdot\left(  w_{2},c_{2}\right)  =\left(
w_{1}+w_{2},c_{1}+c_{2}+\frac{1}{2}\omega\left(  w_{1},w_{2}\right)  \right)
\label{e.3.6}%
\end{equation}
The latter equation may be more simply expressed as%
\begin{equation}
g_{1}g_{2}=g_{1}+g_{2}+\frac{1}{2}\left[  g_{1},g_{2}\right]  , \label{e.3.7}%
\end{equation}
where $g_{i}=\left(  w_{i},c_{i}\right)  $, $i=1,2$. As sets $G$ and
$\mathfrak{g}$ are the same.

The identity in $G$ is $\mathbf{e}=\left(  0,0\right)  $ and the inverse is
given by $g^{-1}=-g$ for all $g=\left(  w,c\right)  \in G.$ Let us observe
that $\left\{  0\right\}  \times\mathbf{C}$ is in the center of both $G$ and
$\mathfrak{g}$ and for $h$ in the center of $G,~$ $g\cdot h=g+h$. In
particular, since $\left[  g,h\right]  \in\left\{  0\right\}  \times
\mathbf{C}$ it follows that $k\cdot\left[  g,h\right]  =k+\left[  g,h\right]
$ for all $k,g,h\in G.$

\begin{df}
\label{d.3.2} When we want to emphasize the group structure on $\mathfrak{g}$
we denote $\mathfrak{g}$ by $G$ or $G\left(  \omega\right)  $. Similarly, when
we view $\mathfrak{g}_{CM}$ as a subgroup of $G$ it will be denoted by
$G_{CM}$ and will be called the \textbf{Cameron--Martin }subgroup.
\end{df}

\begin{lem}
\label{l.3.3}The Banach space topologies on $\mathfrak{g}$ and $\mathfrak{g}%
_{CM}$ make $G$ and $G_{CM}$ into topological groups.
\end{lem}

\begin{proof}
Since $g^{-1}=-g,$ the map $g\mapsto g^{-1}$ is continuous in the
$\mathfrak{g}$ and $\mathfrak{g}_{CM}$ topologies. Since $(g_{1},g_{2})\mapsto
g_{1}+g_{2}$ and $(g_{1},g_{2})\mapsto\left[  g_{1},g_{2}\right]  $ are
continuous in both the $\mathfrak{g}$ and $\mathfrak{g}_{CM}$ topologies, it
follows from Equation \eqref{e.3.7} that $\left(  g_{1},g_{2}\right)  \mapsto
g_{1}\cdot g_{2}$ is continuous as well.
\end{proof}

For later purposes it is useful to observe, by Equations \eqref{e.3.5} and
\eqref{e.3.7}, that
\begin{equation}
\left\Vert g_{1}g_{2}\right\Vert _{\mathfrak{g}}\leqslant\left\Vert
g_{1}\right\Vert _{\mathfrak{g}}+\left\Vert g_{2}\right\Vert _{\mathfrak{g}%
}+\frac{1}{2}\left\Vert \omega\right\Vert _{0}\left\Vert g_{1}\right\Vert
_{\mathfrak{g}}\left\Vert g_{2}\right\Vert _{\mathfrak{g}}\text{ for any
}g_{1},g_{2}\in G. \label{e.3.8}%
\end{equation}

\begin{nota}
\label{n.3.4}To each $g\in G,$ let $l_{g}:G\rightarrow G$ and $r_{g}%
:G\rightarrow G$ denote left and right multiplication by $g$ respectively.
\end{nota}

\begin{nota}
[Linear differentials]\label{n.3.5} Suppose $f:G\rightarrow\mathbb{C}$ is a
Frech\'{e}t smooth function. For $g\in G$ and $h,k\in\mathfrak{g}$ let
\[
f^{\prime}\left(  g\right)  h:=\partial_{h}f\left(  g\right)  =\frac{d}%
{dt}\Big|_{0}f\left(  g+th\right)
\]
and%
\[
f^{\prime\prime}\left(  g\right)  \left(  h\otimes k\right)  :=\partial
_{h}\partial_{k}f\left(  g\right) .
\]
Here and in the sequel a prime on a symbol will be used denote its derivative
or differential.
\end{nota}

As $G$ is a vector space, to each $g\in G$ we can associate the tangent space
(as in the following notation) to $G$ at $g,$ $T_{g}G,$ which is naturally
isomorphic to $G.$

\begin{nota}
\label{n.3.6}For $v,g\in G,$ let $v_{g}\in T_{g}G$ denote the tangent vector
satisfying, $v_{g}f=f^{\prime}\left(  g\right)  v$ for all Frech\'{e}t smooth
functions, $f:G\rightarrow\mathbb{C}.$
\end{nota}

We will write $\mathfrak{g}$ and $\mathfrak{g}_{CM}$ for $T_{\mathbf{e}}G$ and
$T_{\mathbf{e}}G_{CM}$ respectively. Of course as sets we may view
$\mathfrak{g}$ and $\mathfrak{g}_{CM}$ as $G$ and $G_{CM}$ respectively. For
$h\in\mathfrak{g},$ let $\tilde{h}$ be the \textbf{left invariant vector
field} on $G$ such that $\tilde{h}\left(  g\right)  =h$ when $g=\mathbf{e}.$
More precisely if $\sigma\left(  t\right)  \in G$ is any smooth curve such
that $\sigma\left(  0\right)  =\mathbf{e}$ and $\dot{\sigma}\left(  0\right)
=h$ (e.g. $\sigma\left(  t\right)  =th),$ then
\begin{equation}
\tilde{h}\left(  g\right)  =l_{g\ast}h:=\frac{d}{dt}\Big|_{0}g\cdot
\sigma\left(  t\right)  . \label{e.3.9}%
\end{equation}
As usual we view $\tilde{h}$ as a first order differential operator acting on
smooth functions, $f:G\rightarrow\mathbb{C},$ by
\begin{equation}
\left(  \tilde{h}f\right)  \left(  g\right)  =\frac{d}{dt}\Big|_{0}f\left(
g\cdot\sigma\left(  t\right)  \right)  . \label{e.3.10}%
\end{equation}

\begin{prop}
\label{p.3.7} Let $f:G\rightarrow\mathbb{C}$ be a smooth function,
$h=(A,a)\in\mathfrak{g}$ and $g=\left(  w,c\right)  \in G.$ Then
\begin{equation}
\widetilde{h}\left(  g\right)  :=l_{g\ast}h=\left(  A,a+\frac{1}{2}%
\omega\left(  w,A\right)  \right)  _{g}\text{ for all }g=\left(  w,c\right)
\in G \label{e.3.11}%
\end{equation}
and in particular using Notation \ref{n.3.6}
\begin{equation}
\widetilde{\left(  A,a\right)  }f\left(  g\right)  =f^{\prime}\left(
g\right)  \left(  A,a+\frac{1}{2}\omega\left(  w,A\right)  \right)  .
\label{e.3.12}%
\end{equation}
Furthermore, if $h=\left(  A,a\right)  ,$ $k=\left(  B,b\right)  ,$ and then%
\begin{equation}
\left(  \tilde{h}\tilde{k}f-\tilde{k}\tilde{h}f\right)  =\widetilde{\left[
h,k\right]  }f. \label{e.3.13}%
\end{equation}
In other words, the Lie algebra structure on $\mathfrak{g}$ induced by the Lie
algebra structure on the left invariant vector fields on $G$ is the same as
the Lie algebra structure defined in Eq. \eqref{e.3.5}.
\end{prop}

\begin{proof}
Since $th=t\left(  A,a\right)  $ is a curve in $G$ passing through the
identity at $t=0,$ we have%
\begin{align*}
\widetilde{h}\left(  g\right)   &  =\frac{d}{dt}\Big|_{0}\left[  g\cdot\left(
th\right)  \right]  =\frac{d}{dt}\Big|_{0}\left[  \left(  w,c\right)  \cdot
t(A,a)\right] \\
&  =\frac{d}{dt}\Big|_{0}\left[  \left(  w+tA,c+t\delta+\frac{t}{2}%
\omega\left(  w,A\right)  \right)  \right] \\
&  =\left(  A,a+\frac{1}{2}\omega\left(  w,A\right)  \right)  .
\end{align*}
So by the chain rule, $\left(  \tilde{h}f\right)  \left(  g\right)
=f^{\prime}\left(  g\right)  \tilde{h}\left(  g\right)  $ and hence%
\begin{align}
\left(  \tilde{h}\tilde{k}f\right)  \left(  g\right)   &  =\frac{d}%
{dt}\Big|_{0}\left[  f^{\prime}\left(  g\cdot th\right)  \tilde{k}\left(
g\cdot th\right)  \right] \nonumber\\
&  =f^{\prime\prime}\left(  g\right)  \left(  \tilde{h}\left(  g\right)
\otimes\tilde{k}\left(  g\right)  \right)  +f^{\prime}\left(  g\right)
\frac{d}{dt}\Big|_{0}\tilde{k}\left(  g\cdot th\right)  , \label{e.3.14}%
\end{align}
where%
\[
\frac{d}{dt}\Big|_{0}\tilde{k}\left(  g\cdot th\right)  =\frac{d}{dt}%
\Big|_{0}\left(  B,a+\frac{1}{2}\omega\left(  w+tA,B\right)  \right)  =\left(
0,\frac{1}{2}\omega\left(  A,B\right)  \right)  .
\]
Since $f^{\prime\prime}\left(  g\right)  $ is symmetric, it now follows by
subtracting Equation \eqref{e.3.14} with $h$ and $k$ interchanged from itself
that%
\[
\left(  \tilde{h}\tilde{k}f-\tilde{k}\tilde{h}f\right)  \left(  g\right)
=f^{\prime}\left(  g\right)  \left(  0,\omega\left(  A,B\right)  \right)
=f^{\prime}\left(  g\right)  \left[  h,k\right]  =\left(  \widetilde{\left[
h,k\right]  }f\right)  \left(  g\right)
\]
as desired.
\end{proof}

\begin{lem}
\label{l.3.8}The one parameter group in $G,$ $e^{th},$ determined by
$h=\left(  A,a\right)  \in\mathfrak{g},$ is given by%
\begin{equation}
e^{th}=th=t\left(  A,\delta\right)  . \label{e.3.15}%
\end{equation}

\end{lem}

\begin{proof}
Letting $\left(  w\left(  t\right)  , c\left(  t\right)  \right)  :=e^{th},$
according to Equation \eqref{e.3.11} we have that
\[
\frac{d}{dt}\left(  w\left(  t\right)  , c\left(  t\right)  \right)  =\left(
A,a+\frac{1}{2}\omega\left(  w\left(  t\right)  , A\right)  \right)  \text{
with }w\left(  0\right)  =0\text{ and }c\left(  0\right)  =0.
\]
The solution to this differential equation is easily seen to be given by
Equation \eqref{e.3.15}.
\end{proof}

\subsection{Length and distance estimates\label{s.3.1}}

\begin{nota}
\label{n.3.9} Let $T>0$ and $C_{CM}^{1}$ denote the collection of $C^{1}%
$-paths, $g:\left[  0,T\right]  \rightarrow G_{CM}.$ The length of $g$ is
defined as
\begin{equation}
\ell_{G_{CM}}\left(  g\right)  =\int_{0}^{T}\left\Vert l_{g^{-1}\left(
s\right)  \ast}g^{\prime}\left(  s\right)  \right\Vert _{\mathfrak{g}_{CM}}ds.
\label{e.3.16}%
\end{equation}
As usual, the Riemannian distance between $x,y\in G_{CM}$ is defined as%
\[
d_{G_{CM}}(x,y)=\inf\left\{  \ell_{G_{CM}}\left(  g\right)  :g\in C_{CM}%
^{1}\text{ such that }g\left(  0\right)  =x\text{ and }g\left(  T\right)
=y\right\}  .
\]
It will also be convenient to define $\left\vert y\right\vert :=d_{G_{CM}%
}\left(  \mathbf{e},y\right)  $ for all $y\in G_{CM}.$ (The value of $T>0$
used in defining $d_{C_{CM}}$ is irrelevant since the length functional is
reparametrization invariant.)
\end{nota}

Let%
\begin{equation}
C:=\sup\left\{  \left\Vert \omega\left(  h,k\right)  \right\Vert _{\mathbf{C}%
}:\left\Vert h\right\Vert _{H}=\left\Vert k\right\Vert _{H}=1\right\}
\leqslant C_{2}\left\Vert \omega\right\Vert _{0}<\infty. \label{e.3.17}%
\end{equation}
The inequality in Eq. (\ref{e.3.17}) is a consequence of Eq. (\ref{e.2.10})
and the definition of $\left\Vert \omega\right\Vert _{0}$ in Eq. (\ref{e.3.1}).

\begin{prop}
\label{p.3.10}Let $\varepsilon:=1/C$ where $C$ is as in Eq. (\ref{e.3.17}).
Then for all $x,y\in G_{CM},$
\begin{equation}
d_{G_{CM}}\left(  x,y\right)  \leqslant\left(  1+\frac{C}{2}\left\Vert
x\right\Vert _{\mathfrak{g}_{CM}}\wedge\left\Vert y\right\Vert _{\mathfrak{g}%
_{CM}}\right)  \left\Vert y-x\right\Vert _{\mathfrak{g}_{CM}} \label{e.3.18}%
\end{equation}
and in particular, $\left\vert x\right\vert =d_{G_{CM}}\left(  \mathbf{e}%
,x\right)  \leqslant\left\Vert x\right\Vert _{\mathfrak{g}_{CM}}.$ Moreover,
there exists $K<\infty$ such that if $x,y\in G_{CM}$ with $d_{G_{CM}}\left(
x,y\right)  <\varepsilon/2=1/2C,$ then
\begin{equation}
\left\Vert y-x\right\Vert _{\mathfrak{g}_{CM}}\leqslant K\left(  1+\left\Vert
x\right\Vert _{\mathfrak{g}_{CM}}\wedge\left\Vert y\right\Vert _{\mathfrak{g}%
_{CM}}\right)  d_{G_{CM}}\left(  x,y\right)  . \label{e.3.19}%
\end{equation}
As a consequence of Eqs. (\ref{e.3.18}) and (\ref{e.3.19}) we see that the
topology on $G_{CM}$ induced by $d_{G_{CM}}$ is the same as the Hilbert
topology induced by $\left\Vert \cdot\right\Vert _{\mathfrak{g}_{CM}}.$
\end{prop}

\begin{rem}
\label{r.3.11} The equivalence of these two topologies in an
infinite-dimensional setting has been addressed in \cite{Gordina2000b} in the
case of Hilbert-Schmidt groups of operators.
\end{rem}

\begin{proof}
For notational simplicity, let $T=1.$ If $g\left(  s\right)  =\left(  w\left(
s\right) , a\left(  s\right)  \right)  $ is a path in $C_{CM}^{1}$ for
$0\leqslant s\leqslant1,$ then by Equation \eqref{e.3.11}%
\begin{align}
l_{g^{-1}\left(  s\right)  \ast}g^{\prime}\left(  s\right)   &  =\left(
w^{\prime}(s), a^{\prime}(s)-\frac{1}{2}\omega(w(s),w^{\prime}(s))\right)
\label{e.3.20}\\
&  =g^{\prime}\left(  s\right)  -\frac{1}{2}\left[  g\left(  s\right)  ,
g^{\prime}\left(  s\right)  \right] \nonumber
\end{align}
and we may write Equation \eqref{e.3.16} more explicitly as%
\begin{equation}
\ell_{G_{CM}}\left(  g\right)  =\int_{0}^{1}\left\Vert g^{\prime}\left(
s\right)  -\frac{1}{2}\left[  g\left(  s\right)  ,g^{\prime}\left(  s\right)
\right]  \right\Vert _{\mathfrak{g}_{CM}}ds. \label{e.3.21}%
\end{equation}
If we now apply Equation \eqref{e.3.21} to $g\left(  s\right)  =x+s\left(
y-x\right)  $ for $0\leqslant s\leqslant1,$ we see that
\begin{align*}
d_{G_{CM}}\left(  x,y\right)   &  \leqslant\ell_{G_{CM}}\left(  g\right)
=\int_{0}^{1}\left\Vert \left(  y-x\right)  -\frac{1}{2}\left[  x+s\left(
y-x\right)  ,\left(  y-x\right)  \right]  \right\Vert _{\mathfrak{g}_{CM}}ds\\
&  =\left\Vert \left(  y-x\right)  -\frac{1}{2}\left[  x,\left(  y-x\right)
\right]  \right\Vert _{\mathfrak{g}_{CM}}\leqslant\left(  1+\frac{C}%
{2}\left\Vert x\right\Vert _{\mathfrak{g}_{CM}}\right)  \left\Vert
y-x\right\Vert _{\mathfrak{g}_{CM}}.
\end{align*}
As we may interchange the roles of $x$ and $y$ in this inequality, the proof
of Equation \eqref{e.3.18} is complete.

Let
\[
B_{\varepsilon}:=\left\{  x\in\mathfrak{g}_{CM}:\left\Vert x\right\Vert
_{\mathfrak{g}_{CM}}\leqslant\varepsilon\right\}  ,
\]
$y\in B_{\varepsilon},$ and $g:\left[  0,1\right]  \rightarrow G_{CM}$ be a
$C^{1}$--path such that $g\left(  0\right)  =\left(  0,0\right)  =\mathbf{e}$
and $g\left(  1\right)  =y.$ Further let $T\in\left[  0,1\right]  $ be the
first time that $g$ exits $B_{\varepsilon}$ with the convention that $T=1$ if
$g\left(  \left[  0,1\right]  \right)  \subset B_{\varepsilon}.$ Then from
Equation \eqref{e.3.21}
\begin{align}
\ell_{G_{CM}}\left(  g\right)   &  \geqslant\ell_{G_{CM}}\left(  g|_{\left[
0,T\right]  }\right) \nonumber\\
&  \geqslant\int_{0}^{T}\left[  \left\Vert g^{\prime}\left(  s\right)
\right\Vert _{\mathfrak{g}_{CM}}-\frac{1}{2}\left\Vert \left[  g\left(
s\right)  ,g^{\prime}\left(  s\right)  \right]  \right\Vert _{\mathfrak{g}%
_{CM}}\right]  ds\nonumber\\
&  \geqslant\left(  1-\frac{C}{2}\varepsilon\right)  \cdot\int_{0}%
^{T}\left\Vert g^{\prime}\left(  s\right)  \right\Vert _{\mathfrak{g}_{CM}%
}ds\geqslant\left(  1-\frac{C}{2}\varepsilon\right)  \cdot\left\Vert g\left(
T\right)  \right\Vert _{\mathfrak{g}_{CM}}\nonumber\\
&  \geqslant\frac{1}{2}\left\Vert g\left(  T\right)  \right\Vert
_{\mathfrak{g}_{CM}}\geqslant\frac{1}{2}\left\Vert y\right\Vert _{\mathfrak{g}%
_{CM}}. \label{e.3.22}%
\end{align}
Optimizing Equation \eqref{e.3.22} over $g$ implies
\[
\left\vert y\right\vert =d_{G_{CM}}\left(  \mathbf{e},y\right)  \geqslant
\frac{1}{2}\left\Vert y\right\Vert _{\mathfrak{g}_{CM}}\text{ for all }y\in
B_{\varepsilon}.
\]
If in the above argument $y$ was not in $B_{\varepsilon},$ then the path $g$
would have had to exit $B_{\varepsilon}$ and we could conclude that
$\ell_{G_{CM}}\left(  g\right)  \geqslant\left\Vert g\left(  T\right)
\right\Vert _{\mathfrak{g}_{CM}}/2=\varepsilon/2$ and therefore that
$d_{G_{CM}}\left(  \mathbf{e},y\right)  \geqslant\varepsilon/2.$ Hence we have
shown that
\[
\left\vert y\right\vert =d_{G_{CM}}\left(  \mathbf{e},y\right)  \geqslant
\frac{1}{2}\min\left(  \varepsilon,\left\Vert y\right\Vert _{\mathfrak{g}%
_{CM}}\right)  \text{ for all }y\in G_{CM}.
\]

Now suppose that $x,y\in G_{CM}$ and (without loss of generality) that
$\left\Vert x\right\Vert _{\mathfrak{g}_{CM}}\leqslant\left\Vert y\right\Vert
_{\mathfrak{g}_{CM}}.$ Using the left invariance of $d_{G_{CM}},$ it follows
that
\begin{equation}
d_{G_{CM}}\left(  x,y\right)  =d_{G_{CM}}\left(  \mathbf{e},x^{-1}y\right)
\geqslant\frac{1}{2}\min\left(  \varepsilon,\left\Vert x^{-1}y\right\Vert
_{\mathfrak{g}_{CM}}\right)  . \label{e.3.23}%
\end{equation}
If we further suppose that $d_{G_{CM}}\left(  x,y\right)  <\frac{\varepsilon
}{2},$ we may conclude from Equation \eqref{e.3.23} that
\[
\left\Vert y-x-\frac{1}{2}\left[  x,y\right]  \right\Vert _{\mathfrak{g}_{CM}%
}=\left\Vert x^{-1}y\right\Vert _{\mathfrak{g}_{CM}}\leqslant2d_{G_{CM}%
}\left(  x,y\right)  .
\]
If we write $x=\left(  A,a\right)  $ and $y=\left(  B,b\right)  ,$ it follows
that
\[
\left\Vert B-A\right\Vert _{H}^{2}+\left\Vert b-a-\frac{1}{2}\omega\left(
A,B\right)  \right\Vert _{\mathbf{C}}^{2}\leqslant4d_{G_{CM}}^{2}\left(
x,y\right)
\]
and therefore $\left\Vert B-A\right\Vert _{H}\leqslant2d_{G_{CM}}\left(
x,y\right)  $ and%
\begin{align*}
\left\Vert b-a\right\Vert _{\mathbf{C}}  &  \leqslant\left\Vert b-a-\frac
{1}{2}\omega\left(  A,B\right)  \right\Vert _{\mathbf{C}}+\left\Vert \frac
{1}{2}\omega\left(  A,B\right)  \right\Vert _{\mathbf{C}}\\
&  \leqslant2d_{G_{CM}}\left(  x,y\right)  +\frac{1}{2}\left\Vert
\omega\left(  A,B-A\right)  \right\Vert _{\mathbf{C}}\\
&  \leqslant2d_{G_{CM}}\left(  x,y\right)  +\frac{C}{2}\left\Vert A\right\Vert
_{H}\left\Vert B-A\right\Vert _{H}\\
&  \leqslant2d_{G_{CM}}\left(  x,y\right)  \left(  1+\frac{C}{2}\left\Vert
A\right\Vert _{H}\right)  \leqslant2d_{G_{CM}}\left(  x,y\right)  \left(
1+\frac{C}{2}\left\Vert x\right\Vert _{\mathfrak{g}_{CM}}\right)  .
\end{align*}
Combining these results shows that if $d_{G_{CM}}\left(  x,y\right)
<\frac{\varepsilon}{2}$ then%
\[
\left\Vert y-x\right\Vert _{\mathfrak{g}_{CM}}^{2}\leqslant4d_{G_{CM}}%
^{2}\left(  x,y\right)  \left(  1+\left(  1+\frac{C}{2}\left\Vert x\right\Vert
_{\mathfrak{g}_{CM}}\right)  ^{2}\right)
\]
from which Equation \eqref{e.3.19} easily follows.
\end{proof}

We are most interested in the case where $\left\{  \omega\left(  A,B\right)
:A,B\in H\right\}  $ is a total subset of $\mathbf{C,}$ i.e.
$\operatorname*{span}\left\{  \omega\left(  A,B\right)  :A,B\in H\right\}
=\mathbf{C}.$ In this case it turns out that straight line paths are bad
approximations to the geodesics joining $\mathbf{e}\in G_{CM}$ to points $x\in
G_{CM}$ far away from $\mathbf{e}.$ For points $x\in G_{CM}$ distant from
$\mathbf{e}$ it is better to use \textquotedblleft
horizontal\textquotedblright\ paths instead which leads to the following
distance estimates.

\begin{thm}
\label{t.3.12}Suppose that $\left\{  \omega\left(  A,B\right)  :A,B\in
H\right\}  $ is a total subset of $\mathbf{C}.$ Then there exists $C\left(
\omega\right)  <\infty$ such that%
\begin{equation}
d_{CM}\left(  \mathbf{e},\left(  A,a\right)  \right)  \leqslant C\left(
\omega\right)  \left(  \left\Vert A\right\Vert _{H}+\sqrt{\left\Vert
a\right\Vert _{\mathbf{C}}}\right)  \text{ for all }\left(  A,a\right)
\in\mathfrak{g}_{CM}. \label{e.3.24}%
\end{equation}
Moreover, for any $\varepsilon_{0}>0$ there exists $\gamma\left(
\varepsilon_{0}\right)  >0$ such that and
\begin{equation}
\gamma\left(  \varepsilon_{0}\right)  \left(  \left\Vert A\right\Vert
_{H}+\sqrt{\left\Vert a\right\Vert _{\mathbf{C}}}\right)  \leqslant
d_{CM}\left(  \mathbf{e},\left(  A,a\right)  \right)  \text{ if }d_{CM}\left(
\mathbf{e},\left(  A,a\right)  \right)  \geq\varepsilon_{0}. \label{e.3.25}%
\end{equation}
Thus away from any neighborhood of the identity, $d_{CM}\left(  \mathbf{e}%
,\left(  A,a\right)  \right)  $ is comparable to $\left\Vert A\right\Vert
_{H}+\sqrt{\left\Vert a\right\Vert _{\mathbf{C}}}.$
\end{thm}

Since this theorem is not central to the rest of the paper we will relegate
its proof to Appendix \ref{s.12}. The main point of Theorem \ref{t.3.12} is to
explain why Theorem \ref{t.4.16} is an infinite dimensional analogue of the
integrated Gaussian heat kernel bound in Eq. (\ref{e.1.5}).

\subsection{Norm estimates\label{s.3.2}}

\begin{nota}
\label{n.3.13} Suppose $H$ and $\mathbf{C}$ are real (complex) Hilbert spaces,
$L:H\rightarrow\mathbf{C}$ is a bounded operator, $\omega:H\times
H\rightarrow\mathbf{C}$ is a continuous (complex) bilinear form, and $\left\{
e_{j}\right\}  _{j=1}^{\infty}$ is an orthonormal basis for $H.$ The
Hilbert--Schmidt norms of $L$ and $\omega$ are defined by
\begin{equation}
\left\Vert L\right\Vert _{H^{\ast}\otimes\mathbf{C}}^{2}:=\sum_{j=1}^{\infty
}\left\Vert Le_{j}\right\Vert _{\mathbf{C}}^{2}, \label{e.3.26}%
\end{equation}
and%
\begin{equation}
\left\Vert \omega\right\Vert _{2}^{2}=\left\Vert \omega\right\Vert _{H^{\ast
}\otimes H^{\ast}\otimes\mathbf{C}}:=\sum_{i,j=1}^{\infty}\left\Vert
\omega\left(  e_{i},e_{j}\right)  \right\Vert _{\mathbf{C}}^{2}.
\label{e.3.27}%
\end{equation}

\end{nota}

It is easy to verify directly that these definitions are basis independent.
Also see Equation \eqref{e.3.29} below.

\begin{prop}
\label{p.3.14} Suppose that $\left(  W,H,\mu\right)  $ is a real abstract
Wiener space, $\omega:W\times W\rightarrow\mathbf{C}$ is as in Notation
\ref{n.3.1}, and $\left\{  e_{j}\right\}  _{j=1}^{\infty}$ is an orthonormal
basis for $H.$ Then%
\begin{equation}
\left\Vert \omega\left(  w,\cdot\right)  \right\Vert _{H^{\ast}\otimes
\mathbf{C}}^{2}\leqslant C_{2}\left\Vert \omega\right\Vert _{0}^{2}\left\Vert
w\right\Vert _{W}^{2}~\text{ for all }w\in W \label{e.3.28}%
\end{equation}
and
\begin{equation}
\left\Vert \omega\right\Vert _{2}^{2}=\int_{W\times W}\left\Vert \omega\left(
w,w^{\prime}\right)  \right\Vert _{\mathbf{C}}^{2}d\mu\left(  w\right)
d\mu\left(  w^{\prime}\right)  \leqslant\left\Vert \omega\right\Vert _{0}%
^{2}C_{2}^{2}<\infty, \label{e.3.29}%
\end{equation}
where $C_{2}$ is as in Equation \eqref{e.2.2}.
\end{prop}

\begin{proof}
From Equation \eqref{e.2.13},%
\begin{align*}
\left\Vert \omega\left(  w,\cdot\right)  \right\Vert _{H^{\ast}\otimes
\mathbf{C}}^{2}  &  =\int_{W}\left\Vert \omega\left(  w,w^{\prime}\right)
\right\Vert _{\mathbf{C}}^{2}d\mu\left(  w^{\prime}\right) \\
&  \leqslant\left\Vert \omega\right\Vert _{0}^{2}\left\Vert w\right\Vert
_{W}^{2}\int_{W}\left\Vert w^{\prime}\right\Vert _{W}^{2}d\mu\left(
w^{\prime}\right)  =C_{2}\left\Vert \omega\right\Vert _{0}^{2}\left\Vert
w\right\Vert _{W}^{2}.
\end{align*}
Similarly, viewing $w\rightarrow\omega\left(  w,\cdot\right)  $ as a
continuous linear map from $W$ to $H^{\ast}\otimes\mathbf{C}$ it follows from
Eqs. \eqref{e.2.13} and \eqref{e.2.14}, that%
\begin{align*}
\left\Vert \omega\right\Vert _{2}^{2}  &  =\left\Vert h\mapsto\omega\left(
h,\cdot\right)  \right\Vert _{H^{\ast}\otimes\left(  H^{\ast}\otimes
\mathbf{C}\right)  }^{2}=\int_{W}\left\Vert \omega\left(  w,\cdot\right)
\right\Vert _{H^{\ast}\otimes\mathbf{C}}^{2}d\mu\left(  w\right) \\
&  \leqslant\int_{W}C_{2}\left\Vert \omega\right\Vert _{0}^{2}\left\Vert
w\right\Vert _{W}^{2}d\mu\left(  w\right)  =C_{2}^{2}\left\Vert \omega
\right\Vert _{0}^{2}.
\end{align*}

\end{proof}

\begin{rem}
\label{r.3.15}
%\label{l.3.13}
The Lie bracket on $\mathfrak{g}_{CM}$ has the following continuity property,
\[
\left\Vert \left[  \left(  A,a\right)  ,\left(  B,b\right)  \right]
\right\Vert _{\mathfrak{g}_{CM}}\leqslant C\left\Vert \left(  A,a\right)
\right\Vert _{\mathfrak{g}_{CM}}\left\Vert \left(  B,b\right)  \right\Vert
_{\mathfrak{g}_{CM}}%
\]
where $C\leq\left\Vert \omega\right\Vert _{2}$ as in Eq. (\ref{e.3.17}). This
is a consequence of the following simple estimates
\begin{align*}
\left\Vert \left[  \left(  A,a\right)  ,\left(  B,b\right)  \right]
\right\Vert _{\mathfrak{g}_{CM}}  &  =\left\Vert \left(  0,\omega\left(
A,B\right)  \right)  \right\Vert _{\mathfrak{g}_{CM}}=\left\Vert \omega\left(
A,B\right)  \right\Vert _{\mathbf{C}}\\
&  \leqslant C\left\Vert A\right\Vert _{H}\left\Vert B\right\Vert
_{H}\leqslant C\left\Vert \left(  A,a\right)  \right\Vert _{\mathfrak{g}_{CM}%
}\left\Vert \left(  B,b\right)  \right\Vert _{\mathfrak{g}_{CM}}.
\end{align*}
This continuity property of the Lie bracket is often used to prove that the
exponential map is a local diffeomorphism (e.g. see \cite{Gordina2000b} in the
case of infinite-dimensional matrix groups). In the Heisenberg group setting
the exponential map is a global diffeomorphism as follows from Lemma
\ref{l.3.8}, where we have not used continuity of the Lie bracket.
\end{rem}

\begin{lem}
\label{l.3.16}Suppose that $H$ is a real Hilbert space, $\mathbf{C}$ is a real
finite-dimensional inner product space, and $\ell:H\rightarrow\mathbf{C}$ is a
continuous linear map. Then for any orthonormal basis $\left\{  e_{j}\right\}
_{j=1}^{\infty}$ of $H$ the series
\begin{equation}
\sum_{j=1}^{\infty}\ell\left(  e_{j}\right)  \otimes\ell\left(  e_{j}\right)
\in\mathbf{C}\otimes\mathbf{C} \label{e.3.30}%
\end{equation}
and
\begin{equation}
\sum_{j=1}^{\infty}\ell\left(  e_{j}\right)  \otimes e_{j}\in\mathbf{C}\otimes
H \label{e.3.31}%
\end{equation}
are convergent and independent of the basis.\footnote{If we were to allow
$\mathbf{C}$ to be an infinite-dimensional Hilbert space here, we would have
to assume that $\ell$ is Hilbert--Schmidt. When $\dim\mathbf{C<}\infty,$
$\ell:H\rightarrow\mathbf{C}$ is Hilbert--Schmidt iff it is bounded.}
\end{lem}

\begin{proof}
If $\left\{  f_{i}\right\}  _{i=1}^{\dim\mathbf{C}}$ is an orthonormal basis
for $\mathbf{C},$ then%
\begin{align*}
\sum_{j=1}^{\infty}\left\Vert \ell\left(  e_{j}\right)  \otimes\ell\left(
e_{j}\right)  \right\Vert _{\mathbf{C}\otimes\mathbf{C}}  &  =\sum
_{j=1}^{\infty}\left\Vert \ell\left(  e_{j}\right)  \right\Vert _{\mathbf{C}%
}^{2}\\
&  =\sum_{i=1}^{\dim\mathbf{C}}\sum_{j=1}^{\infty}\left(  f_{i},\ell\left(
e_{j}\right)  \right)  ^{2}=\sum_{i=1}^{\dim\mathbf{C}}\left\Vert \left(
f_{i},\ell\left(  \cdot\right)  \right)  \right\Vert _{H^{\ast}}^{2}<\infty
\end{align*}
which shows that the sum in Equation \eqref{e.3.30} is absolutely convergent
and that $\ell$ is Hilbert--Schmidt. Similarly, since $\left\{  \ell\left(
e_{j}\right)  \otimes e_{j}\right\}  _{j=1}^{\infty}$ is an orthogonal set in
$\mathbf{C}\otimes H$ and
\[
\sum_{j=1}^{\infty}\left\Vert \ell\left(  e_{j}\right)  \otimes e_{j}%
\right\Vert _{\mathbf{C}\otimes H}^{2}=\sum_{j=1}^{\infty}\left\Vert
\ell\left(  e_{j}\right)  \right\Vert _{\mathbf{C}}^{2}<\infty,
\]
the sum in Equation \eqref{e.3.31} is convergent as well.

Recall that if $H$ and $K$ are two real Hilbert spaces then the Hilbert space
tensor product, $H\otimes K,$ is unitarily equivalent to the space of
Hilbert--Schmidt operators, $HS\left(  H,K\right)  ,$ from $H$ to $K.$ Under
this identification, $h\otimes k\in H\otimes K$ corresponds to the operator
(still denoted by $h\otimes k)$ in $HS\left(  H,K\right)  $ defined by; $H\ni
h^{\prime}\rightarrow\left(  h,h^{\prime}\right)  _{H}\,k\in K.$ Using this
identification we have that for all $c\in\mathbf{C;}$%
\begin{align*}
\left(  \sum_{j=1}^{\infty}\ell\left(  e_{j}\right)  \otimes\ell\left(
e_{j}\right)  \right)  c  &  =\sum_{j=1}^{\infty}\ell\left(  e_{j}\right)
\left\langle \ell\left(  e_{j}\right)  ,c\right\rangle _{\mathbf{C}}%
=\sum_{j=1}^{\infty}\ell\left(  e_{j}\right)  \left\langle e_{j},\ell^{\ast
}c\right\rangle _{\mathbf{C}}\\
&  =\ell\left(  \sum_{j=1}^{\infty}\left\langle e_{j},\ell^{\ast
}c\right\rangle _{\mathbf{C}}e_{j}\right)  =\ell\ell^{\ast}c
\end{align*}
and%
\[
\left(  \sum_{j=1}^{\infty}\ell\left(  e_{j}\right)  \otimes e_{j}\right)
c=\sum_{j=1}^{\infty}e_{j}\left\langle \ell\left(  e_{j}\right)
,c\right\rangle _{\mathbf{C}}=\sum_{j=1}^{\infty}e_{j}\left\langle e_{j}%
,\ell^{\ast}c\right\rangle _{\mathbf{C}}=\ell^{\ast}c,
\]
which clearly shows that Equations \eqref{e.3.30} and\eqref{e.3.31} are basis--independent.
\end{proof}

\subsection{Examples\label{s.3.3}}

Here we describe several examples including finite-dimensional Heisenberg
groups. As we mentioned earlier a typical example of such a group is the
Heisenberg group of a symplectic vector space. For each of the examples
presented we will explicitly compute the norm $\left\Vert \omega\right\Vert
_{2}^{2}$ of the form $\omega$ as defined in Equation \eqref{e.3.27}. In
Section \ref{s.7} we will also explicitly compute the Ricci tensor for each of
the examples introduced here.

To describe some of the examples below, it is convenient to use complex Banach
and Hilbert spaces. However, for the purposes of this paper the complex
structure on these spaces should be forgotten. In doing so we will use the
following notation. If $X$ is a complex vector space, let
$X_{\operatorname{Re}}$ denote $X$ thought of as a real vector space. If
$\left(  H,\left\langle \cdot,\cdot\right\rangle _{H}\right)  $ is a complex
Hilbert space, we define $\left\langle \cdot,\cdot\right\rangle
_{H_{\operatorname{Re}}}=\operatorname{Re}\left\langle \cdot,\cdot
\right\rangle _{H}$ in which case $\left(  H_{\operatorname{Re}},\left\langle
\cdot,\cdot\right\rangle _{H_{\operatorname{Re}}}\right)  $ becomes a real
Hilbert space. Before going to the examples, let us record the relationship
between the complex and real Hilbert Schmidt norms of Notation \ref{n.3.13}.

\begin{lem}
\label{l.3.17}Suppose $H$ and $\mathbf{C}$ are complex Hilbert spaces,
$L:H\rightarrow\mathbf{C}$ and $\omega:H\times H\rightarrow\mathbf{C}$ are as
in Notation \ref{n.3.13}, and $c\in\mathbf{C}.$ Then
\begin{align}
\left\Vert L\right\Vert _{H_{\operatorname{Re}}^{\ast}\otimes\mathbf{C}%
_{\operatorname{Re}}}^{2}  &  =2\left\Vert L\right\Vert _{H^{\ast}%
\otimes\mathbf{C}}^{2},\label{e.3.32}\\
\left\Vert \left\langle \omega\left(  \cdot,\cdot\right)  , c\right\rangle
_{\mathbf{C}_{\operatorname{Re}}}\right\Vert _{H_{\operatorname{Re}}^{\ast
}\otimes H_{\operatorname{Re}}^{\ast}}^{2}  &  =2\left\Vert \left\langle
\omega\left(  \cdot,\cdot\right)  , c\right\rangle _{\mathbf{C}}\right\Vert
_{H^{\ast}\otimes H^{\ast}}^{2}, \label{e.3.33}%
\end{align}
and%
\begin{equation}
\left\Vert \omega\left(  \cdot,\cdot\right)  \right\Vert
_{H_{\operatorname{Re}}^{\ast}\otimes H_{\operatorname{Re}}^{\ast}%
\otimes\mathbf{C}_{\operatorname{Re}}}^{2}=4\left\Vert \omega\left(
\cdot,\cdot\right)  \right\Vert _{H^{\ast}\otimes H^{\ast}\otimes\mathbf{C}%
}^{2}. \label{e.3.34}%
\end{equation}

\end{lem}

\begin{proof}
A straightforward proof.
\end{proof}

\begin{ex}
[Finite-dimensional real Heisenberg group]\label{ex.3.18}Let $\mathbf{C}%
=\mathbb{R},$ $W=H=\left(  \mathbb{C}^{n}\right)  _{\operatorname{Re}}%
\cong\mathbb{R}^{2n},$ and $\omega\left(  w,z\right)  :=\operatorname{Im}%
\left\langle w,z\right\rangle $ be the standard symplectic form on
$\mathbb{R}^{2n},$ where $\left\langle w, z\right\rangle =w\cdot\bar{z}$ is
the usual inner product on $\mathbb{C}^{n}.$ Any element of the group
$\mathbf{H}_{\mathbb{R}}^{n}:=G\left(  \omega\right)  $ can be written as
$g=\left(  z,c\right)  $, where $z\in\mathbb{C}^{n}$ and $c\in\mathbb{R}.$ As
above, the Lie algebra, $\mathfrak{h}_{\mathbb{R}}^{n},$ of $\mathbf{H}%
_{\mathbb{R}}^{n}$ is, as a set, equal to $\mathbf{H}_{\mathbb{R}}^{n}$
itself. If $\{e_{j}\}_{j=1}^{n}$ is an orthonormal basis for $\mathbb{R}^{n}$
then $\left\{  e_{j},ie_{j}\right\}  _{j=1}^{n}$ is an orthonormal basis for
$H$ and (real) Hilbert Schmidt norm of $\omega$ is given by
\begin{equation}
\left\Vert \omega\right\Vert _{H^{\ast}\otimes H^{\ast}}^{2}=\sum_{j,k=1}%
^{n}\sum_{\varepsilon,\delta\in\left\{  1,i\right\}  }\left[
\operatorname{Im}\left\langle \varepsilon e_{j},\delta e_{k}\right\rangle
\right]  ^{2}=\sum_{j,k=1}^{n}2\delta_{j,k}=2n. \label{e.3.35}%
\end{equation}

\end{ex}

\begin{ex}
[Finite-dimensional complex Heisenberg group]\label{ex.3.19} Suppose that
$W=H=\mathbb{C}^{n}\times\mathbb{C}^{n},$ $\mathbf{C}=\mathbb{C},$ and
$\omega:W\times W\rightarrow\mathbb{C}$ is defined by%
\[
\omega\left(  \left(  w_{1},w_{2}\right)  , \left(  z_{1},z_{2}\right)
\right)  =w_{1}\cdot z_{2}-w_{2}\cdot z_{1}.
\]
Any element of the group $\mathbf{H}_{\mathbb{C}}^{n}:=G\left(  \omega\right)
$ can be written as $g=\left(  z_{1},z_{2},c\right)  $, where $z_{1},z_{2}%
\in\mathbb{C}^{n}$ and $c\in\mathbb{C}.$ As above, the Lie algebra,
$\mathfrak{h}_{\mathbb{C}}^{n},$ of $\mathbf{H}_{\mathbb{C}}^{n}$ is, as a
set, equal to $\mathbf{H}_{\mathbb{C}}^{n}$ itself. In this case $\left\{
\left(  e_{j},0\right)  , \left(  0,e_{j}\right)  \right\}  _{j=1}^{n}$ is a
complex orthonormal basis for $H.$ The (complex) Hilbert-Schmidt norm of the
symplectic form $\omega$ is given by
\[
\left\Vert \omega\right\Vert _{H^{\ast}\otimes H^{\ast}}^{2}=2\sum_{j=1}%
^{n}|\omega\left(  \left(  e_{j},0\right)  , \left(  0,e_{j}\right)  \right)
|^{2}=2n.
\]

\end{ex}

\begin{ex}
\label{ex.3.20}Let $\left(  K,\left\langle \cdot,\cdot\right\rangle \right)  $
be a complex Hilbert space and $Q$ be a strictly positive trace class operator
on $K.$ For $h,k\in K,$ let $\left\langle h,k\right\rangle _{Q}:=\left\langle
h,Qk\right\rangle $ and $\left\Vert h\right\Vert _{Q}:=\sqrt{\left\langle
h,h\right\rangle _{Q}}.$ Also let $\left(  K_{Q},\left\langle \cdot
,\cdot\right\rangle _{Q}\right)  $ denote the Hilbert space completion of
$\left(  K,\left\Vert \cdot\right\Vert _{Q}\right)  .$ Analogous to Example
\ref{ex.3.18}, let $H:=K_{\operatorname{Re}},$ $W:=\left(  K_{Q}\right)
_{\operatorname{Re}},$ and $\omega:W\times W\rightarrow\mathbb{R}=:\mathbf{C}$
be defined by%
\[
\omega\left(  w,z\right)  =\operatorname{Im}\left\langle w,z\right\rangle
_{Q}\text{ for all }w,z\in W.
\]
Then $G\left(  \omega\right)  =W\times\mathbb{R}$ is a real group and $\left(
W,H\right)  $ determines an abstract Wiener space (see for example
\cite[Exercise 17, p. 59]{Kuo75} and \cite[Example 3.9.7]{Bog98})). Let
$\left\{  e_{j}\right\}  _{j=1}^{\infty}$ be an orthonormal basis for $K$ so
that $\left\{  e_{j},ie_{j}\right\}  _{j=1}^{\infty}$ is an orthonormal basis
for $\left(  H,\operatorname{Re}\left\langle \cdot,\cdot\right\rangle
_{K}\right)  .$ Then the real Hilbert-Schmidt norm of $\omega$ is given by
\begin{align}
\left\Vert \omega\right\Vert _{H^{\ast}\otimes H^{\ast}}^{2}  &  =\sum
_{j,k=1}^{\infty}\sum_{\varepsilon,\delta\in\left\{  1,i\right\}  }\left[
\operatorname{Im}^{2}\left\langle \varepsilon e_{j},Q\delta e_{k}\right\rangle
\right] \nonumber\\
&  =2\sum_{j,k=1}^{\infty}\left[  \operatorname{Im}^{2}\left\langle
e_{j},Qe_{k}\right\rangle +\operatorname{Re}^{2}\left\langle e_{j}%
,Qe_{k}\right\rangle \right]  =2\sum_{j,k=1}^{\infty}\left\vert \left\langle
e_{j},Qe_{k}\right\rangle \right\vert ^{2}\nonumber\\
&  =2\sum_{k=1}^{\infty}\left\Vert Qe_{k}\right\Vert ^{2}=2\left\Vert
Q\right\Vert _{HS}^{2}=2\operatorname{tr}Q^{2}. \label{e.3.36}%
\end{align}

\end{ex}

\begin{ex}
\label{ex.3.21}Again let $\left(  K,\left\langle \cdot,\cdot\right\rangle
\right)  , $ $Q,$ and $\left(  K_{Q},\left\langle \cdot,\cdot\right\rangle
_{Q}\right)  $ be as in the previous example. Let us further assume that $K$
is equipped with a conjugation, $k\rightarrow\bar{k},$ which is isometric and
commutes with $Q.$ Analogously to Example \ref{ex.3.19}, let $W:=K_{Q}\times
K_{Q},$ $H=K\times K,$ and let $\omega:W\times W\rightarrow\mathbb{C}$ be
defined by%
\[
\omega\left(  \left(  w_{1},w_{2}\right)  , \left(  z_{1},z_{2}\right)
\right)  =\left\langle w_{1},\bar{z}_{2}\right\rangle _{Q}-\left\langle
w_{2},\bar{z}_{1}\right\rangle _{Q},
\]
which is skew symmetric because the conjugation commutes with $Q.$ If
$\left\{  e_{j}\right\}  _{j=1}^{\infty}$ is an orthonormal basis for $K,$
then $\left\{  \bar{e}_{j}\right\}  _{j=1}^{\infty}$ is also an orthonormal
basis for $K$ (because the conjugation is isometric) and $\left\{  \left(
e_{j},0\right)  , \left(  0,e_{j}\right)  \right\}  _{j=1}^{\infty}$ is a
orthonormal basis for $H.$ Hence, the (complex) Hilbert-Schmidt norm of
$\omega$ is given by
\begin{align}
\left\Vert \omega\right\Vert _{H^{\ast}\otimes H^{\ast}}^{2}  &  =\sum
_{j,k=1}^{\infty}\left(  |\omega\left(  \left(  e_{j},0\right)  , \left(
0,e_{k}\right)  \right)  |^{2}+|\omega\left(  \left(  0,e_{k}\right)  ,
\left(  e_{j},0\right)  \right)  |^{2}\right) \nonumber\\
&  =2\sum_{j,k=1}^{\infty}\left\vert \left\langle e_{j},Q\bar{e}%
_{k}\right\rangle \right\vert ^{2}=2\sum_{k=1}^{\infty}\left\Vert Q\bar{e}%
_{k}\right\Vert ^{2}=2\left\Vert Q\right\Vert _{HS}^{2}=2\operatorname{tr}%
Q^{2}. \label{e.3.37}%
\end{align}

\end{ex}

\begin{ex}
\label{ex.3.22}Suppose that $\left(  V,\left\langle \cdot,\cdot\right\rangle
_{V}\right)  $ is a $d$--dimensional $\mathbb{F}$--inner product space
($\mathbb{F}=\mathbb{R}$ or $\mathbb{C)}$, $\mathbf{C}$ is a
finite-dimensional $\mathbb{F}$ -- inner product space, $\alpha:V\times
V\rightarrow\mathbf{C}$ is an anti-symmetric bilinear form on $V,$ and
$\left\{  q_{j}\right\}  _{j=1}^{\infty}$ is a sequence of positive numbers
such that $\sum_{j=1}^{\infty}q_{j}<\infty.$ Let
\begin{align*}
W &  =\left\{  v\in V^{\mathbb{N}}:\sum_{j=1}^{\infty}q_{j}\left\Vert
v_{j}\right\Vert _{V}^{2}<\infty\right\}  \text{ and}\\
H &  =\left\{  v\in V^{\mathbb{N}}:\sum_{j=1}^{\infty}\left\Vert
v_{j}\right\Vert _{V}^{2}<\infty\right\}  \subset W,
\end{align*}
each of which are Hilbert spaces when equipped with the inner products
\begin{align*}
\left\langle v,w\right\rangle _{W} &  :=\sum_{j=1}^{\infty}q_{j}\left\langle
v_{j},w_{j}\right\rangle _{V}\text{ and}\\
\left\langle v,w\right\rangle _{H} &  :=\sum_{j=1}^{\infty}\left\langle
v_{j},w_{j}\right\rangle _{V}%
\end{align*}
respectively. Further let $\omega:W\times W\rightarrow\mathbf{C}$ be defined
by
\[
\omega\left(  v,w\right)  =\sum_{j=1}^{\infty}q_{j}\alpha\left(  v_{j}%
,w_{j}\right)  .
\]
Then $\left(  W_{\operatorname{Re}},H_{\operatorname{Re}}\right)  $ is an
abstract Wiener space (see for example \cite[Exercise 17, p. 59]{Kuo75} and
\cite[Example 3.9.7]{Bog98}) and, since%
\begin{align*}
\left\vert \omega\left(  v,w\right)  \right\vert  &  =\sum_{j=1}^{\infty}%
q_{j}\left\vert \alpha\left(  v_{j},w_{j}\right)  \right\vert \leqslant
\left\Vert \alpha\right\Vert _{0}\sum_{j=1}^{\infty}q_{j}\left\Vert
v_{j}\right\Vert _{V}\left\Vert w_{j}\right\Vert _{V}\\
&  \leqslant\left\Vert \alpha\right\Vert _{0}\left\Vert v\right\Vert
_{W}\left\Vert w\right\Vert _{W},
\end{align*}
we have $\left\Vert \omega\right\Vert _{0}\leqslant\left\Vert \alpha
\right\Vert _{0}.$ For $v\in V,$ let $v\left(  j\right)  :=\left(
0,\dots,0,v,0,0,\dots\right)  \in H$ where the $v$ is put in the
$j^{\text{th}}$ -- position. If $\left\{  u_{a}\right\}  _{a=1}^{d}$ is an
orthonormal basis for $V,$ then $\left\{  u_{a}\left(  j\right)
:a=1,\dots,d\right\}  _{j=1}^{\infty}$ is an orthonormal basis for $H.$
Therefore,%
\begin{align}
\left\Vert \left\langle \omega\left(  \cdot,\cdot\right)  ,c\right\rangle
_{\mathbf{C}}\right\Vert _{H^{\ast}\otimes H^{\ast}}^{2} &  =\sum
_{j,k=1}^{\infty}\sum_{a,b=1}^{d}\left\langle \omega\left(  u_{a}\left(
j\right)  ,u_{b}\left(  k\right)  \right)  ,c\right\rangle _{\mathbf{C}}%
^{2}\nonumber\\
&  =\sum_{j=1}^{\infty}\sum_{a,b=1}^{d}q_{j}^{2}\left\langle \alpha\left(
u_{a},u_{b}\right)  ,c\right\rangle _{\mathbf{C}}^{2}\nonumber\\
&  =\left(  \sum_{j=1}^{\infty}q_{j}^{2}\right)  \left\Vert \left\langle
\alpha\left(  \cdot,\cdot\right)  ,c\right\rangle _{\mathbf{C}}\right\Vert
_{V^{\ast}\otimes V^{\ast}}^{2}\text{ for all }c\in\mathbf{C}.\label{e.3.38}%
\end{align}

\end{ex}

\begin{ex}
\label{ex.3.23}Let $(V,\left\langle \cdot,\cdot\right\rangle ,\alpha)$ be as
in Example \ref{ex.3.22} with $\mathbb{F}=\mathbb{C},$%
\[
W=\left\{  \sigma\in C\left(  \left[  0,1\right]  ,V\right)  :\sigma\left(
0\right)  =0\right\}
\]
and $H$ be the associated Cameron -- Martin space,%
\[
H:=H\left(  V\right)  =\left\{  h\in W:\int_{0}^{1}\left\Vert h^{\prime
}\left(  s\right)  \right\Vert _{V}^{2}ds<\infty\right\}  ,
\]
wherein $\int_{0}^{1}\left\Vert h^{\prime}\left(  s\right)  \right\Vert
_{V\ }^{2}ds:=\infty$ if $h$ is not absolutely continuous. Further let $\eta$
be a complex measure on $\left[  0,1\right]  $ and%
\[
\omega\left(  \sigma_{1},\sigma_{2}\right)  :=\int_{0}^{1}\alpha\left(
\sigma_{1}\left(  s\right)  ,\sigma_{2}\left(  s\right)  \right)  d\eta\left(
s\right)  \text{ for all }\sigma_{1},\sigma_{2}\in W.
\]
Then $\left(  W,H,\omega\right)  $ satisfies all of the assumptions in
Notation \ref{n.3.1}. Let $\left\{  u_{a}\right\}  _{a=1}^{d}$ be the
orthonormal basis of $V$, $\{l_{j}\}_{j=1}^{\infty}$ be the orthonormal basis
of $H(\mathbb{R})$, then $\{l_{j}u_{a}:a=1,2,\dots,d\}_{j=1}^{\infty}$ is an
orthonormal basis of $H$ and (see \cite[Lemma 3.8]{Driver1996b})%
\begin{equation}
\sum_{j=1}^{\infty}l_{j}\left(  s\right)  l_{j}\left(  t\right)  =s\wedge
t\text{ for all }s,t\in\left[  0,1\right]  .\label{e.3.39}%
\end{equation}
If we let $\lambda$ be the total variation of $\eta,$ then $d\eta=\rho
d\lambda$, where $\rho=\frac{d\eta}{d\lambda}.$ Hence if $d\bar{\eta}\left(
t\right)  :=\bar{\rho}\left(  t\right)  d\lambda\left(  t\right)  $ and
$c\in\mathbf{C},$ then%
\begin{align}
&  \left\Vert \left\langle \omega\left(  \cdot,\cdot\right)  ,c\right\rangle
_{\mathbf{C}}\right\Vert _{2}^{2}\nonumber\\
&  \quad=\sum_{j,k=1}^{\infty}\sum_{a,b=1}^{d}\left\vert \left\langle
\omega\left(  l_{j}u_{a},l_{k}u_{b}\right)  ,c\right\rangle _{\mathbf{C}%
}\right\vert ^{2}\nonumber\\
&  \quad=\sum_{j,k=1}^{\infty}\sum_{a,b=1}^{d}\left\vert \int_{\left[
0,1\right]  }l_{j}\left(  s\right)  l_{k}\left(  s\right)  \rho\left(
s\right)  d\lambda\left(  s\right)  \,\left\langle \alpha\left(  u_{a}%
,u_{b}\right)  ,c\right\rangle _{\mathbf{C}}\right\vert ^{2}\nonumber\\
&  \quad=\left\Vert \left\langle \alpha\left(  \cdot,\cdot\right)
,c\right\rangle _{\mathbf{C}}\right\Vert _{2}^{2}\cdot\sum_{j,k=1}^{\infty
}\int_{\left[  0,1\right]  ^{2}}l_{j}\left(  s\right)  l_{k}\left(  s\right)
l_{j}\left(  t\right)  l_{k}\left(  t\right)  \rho\left(  s\right)  \bar{\rho
}\left(  t\right)  d\lambda\left(  s\right)  d\lambda\left(  t\right)
\nonumber\\
&  \quad=\left\Vert \left\langle \alpha\left(  \cdot,\cdot\right)
,c\right\rangle _{\mathbf{C}}\right\Vert _{2}^{2}\cdot\int_{\left[
0,1\right]  ^{2}}\left(  s\wedge t\right)  ^{2}\rho\left(  s\right)  \bar
{\rho}\left(  t\right)  d\lambda\left(  s\right)  d\lambda\left(  t\right)
\nonumber\\
&  \quad=\left\Vert \left\langle \alpha\left(  \cdot,\cdot\right)
,c\right\rangle _{\mathbf{C}}\right\Vert _{2}^{2}\cdot\int_{\left[
0,1\right]  ^{2}}\left(  s\wedge t\right)  ^{2}d\eta\left(  s\right)
d\bar{\eta}\left(  t\right)  ,\label{e.3.40}%
\end{align}
wherein we have used Equation \eqref{e.3.39} in the fourth equality above.
Summing this equation over $c$ in an orthonormal basis for $\mathbf{C}$ shows%
\begin{equation}
\left\Vert \omega\right\Vert _{2}^{2}=\left\Vert \alpha\right\Vert _{2}%
^{2}\cdot\int_{\left[  0,1\right]  ^{2}}\left(  s\wedge t\right)  ^{2}%
d\eta\left(  s\right)  d\bar{\eta}\left(  t\right)  .\label{e.3.41}%
\end{equation}

\end{ex}

\subsection{Finite Dimensional Projections and Cylinder Functions\label{s.3.4}%
}

Let $i:H\rightarrow W$ be the inclusion map, and $i^{\ast}:W^{\ast}\rightarrow
H^{\ast}$ be its transpose, i.e. $i^{\ast}\ell:=\ell\circ i$ for all $\ell\in
W^{\ast}.$ Also let
\[
H_{\ast}:=\left\{  h\in H:\left\langle \cdot,h\right\rangle _{H}%
\in\mathrm{Ran}(i^{\ast})\subset H^{\ast}\right\}
\]
or in other words, $h\in H$ is in $H_{\ast}$ iff $\left\langle \cdot
,h\right\rangle _{H}\in H^{\ast}$ extends to a continuous linear functional on
$W.$ (We will continue to denote the continuous extension of $\left\langle
\cdot,h\right\rangle _{H}$ to $W$ by $\left\langle \cdot,h\right\rangle
_{H}.)$ Because $H$ is a dense subspace of $W,$ $i^{\ast}$ is injective and
because $i$ is injective, $i^{\ast}$ has a dense range. Since $h\mapsto
\left\langle \cdot, h\right\rangle _{H}$ as a map from $H$ to $H^{\ast}$ is a
conjugate linear isometric isomorphism, it follows from the above comments
that $H_{\ast}\ni h\rightarrow\left\langle \cdot,h\right\rangle _{H}\in
W^{\ast}$ is a conjugate linear isomorphism too, and that $H_{\ast}$ is a
dense subspace of $H$.

Now suppose that $P:H\rightarrow H$ is a finite rank orthogonal projection
such that $PH\subset H_{\ast}.$ Let $\left\{  e_{j}\right\}  _{j=1}^{n}$ be an
orthonormal basis for $PH$ and $\ell_{j}=\left\langle \cdot,e_{j}\right\rangle
_{H}\in W^{\ast}.$ Then we may extend $P$ to a (unique) continuous operator
from $W$ $\rightarrow H$ (still denoted by $P)$ by letting%
\begin{equation}
Pw:=\sum_{j=1}^{n}\left\langle k,e_{j}\right\rangle _{H}e_{j}=\sum_{j=1}%
^{n}\ell_{j}\left(  w\right)  e_{j}\text{ for all }w\in W. \label{e.3.42}%
\end{equation}
For all $w\in W$ we have, $\left\Vert Pw\right\Vert _{H}\leqslant C_{2}\left(
P\right)  \left\Vert Pw\right\Vert _{W}$ and
\[
\left\Vert Pw\right\Vert _{W}\leqslant\left(  \sum_{i=1}^{n}\left\Vert
\left\langle \cdot,e_{i}\right\rangle _{H}\right\Vert _{W}\left\Vert
e_{i}\right\Vert _{W}\right)  \left\Vert w\right\Vert _{W}%
\]
and therefore there exists $C<\infty$ such that
\begin{equation}
\left\Vert Pw\right\Vert _{H}\leqslant C\left\Vert w\right\Vert _{W}\text{ for
all }w\in W. \label{e.3.43}%
\end{equation}

\begin{nota}
\label{n.3.24} Let $\operatorname*{Proj}\left(  W\right)  $ denote the
collection of finite rank projections on $W$ such that $PW\subset H_{\ast}$
and $P|_{H}:H\rightarrow H$ is an orthogonal projection, i.e. $P$ has the form
given in Equation \eqref{e.3.42}. Further, let $G_{P}:=PW\times\mathbf{C}$ (a
subgroup of $G_{CM})$ and
\[
\pi=\pi_{P}:G\rightarrow G_{P}%
\]
be defined by $\pi_{P}\left(  w,c\right)  :=\left(  Pw,c\right)  .$
\end{nota}

\begin{rem}
\label{r.3.25}The reader should be aware that $\pi_{P}:G\rightarrow
G_{P}\subset G_{CM}$ is \textbf{not }(for general $\omega$ and $P\in
\operatorname*{Proj}\left(  W\right)  )$ a group homomorphism. In fact we
have,
\begin{equation}
\pi_{P}\left[  \left(  w,c\right)  \cdot\left(  w^{\prime},c^{\prime}\right)
\right]  -\pi_{P}\left(  w,c\right)  \cdot\pi_{P}\left(  w^{\prime},c^{\prime
}\right)  =\Gamma_{P}\left(  w,w^{\prime}\right)  , \label{e.3.44}%
\end{equation}
where
\begin{equation}
\Gamma_{P}\left(  w,w^{\prime}\right)  =\frac{1}{2}\left(  0,\omega\left(
w,w^{\prime}\right)  -\omega\left(  Pw,Pw^{\prime}\right)  \right)  .
\label{e.3.45}%
\end{equation}
So unless $\omega$ is \textquotedblleft supported\textquotedblright\ on the
range of $P,$ $\pi_{P}$ is not a group homomorphism. Since, $\left(
w,b\right)  +\left(  0,c\right)  =\left(  w,b\right)  \cdot\left(  0,c\right)
$ for all $w\in W$ and $b,c\in\mathbf{C},$ we may also write Equation
\ref{e.3.44} as
\begin{equation}
\pi_{P}\left[  \left(  w,c\right)  \cdot\left(  w^{\prime},c^{\prime}\right)
\right]  =\pi_{P}\left(  w,c\right)  \cdot\pi_{P}\left(  w^{\prime},c^{\prime
}\right)  \cdot\Gamma_{P}\left(  w,w^{\prime}\right)  . \label{e.3.46}%
\end{equation}

\end{rem}

\begin{df}
\label{d.3.26} A function $f:G\rightarrow\mathbb{C}$ is said to be a
(\textbf{smooth)} \textbf{cylinder function} if it may be written as
$f=F\circ\pi_{P}$ for some $P\in\operatorname*{Proj}\left(  W\right)  $ and
some (smooth) function $F:G_{P}\mathbb{\rightarrow C}.$
\end{df}

\begin{nota}
\label{n.3.27}For $g=\left(  w,c\right)  \in G,$ let $\gamma\left(  g\right)
$ and $\chi\left(  g\right)  $ be the elements of $\mathfrak{g}_{CM}%
\otimes\mathfrak{g}_{CM}$ defined by
\begin{align*}
\gamma\left(  g\right)   &  :=\sum_{j=1}^{\infty}\left(  0,\omega\left(
w,e_{j}\right)  \right)  \otimes\left(  e_{j},0\right)  \text{ and}\\
\chi\left(  g\right)   &  :=\sum_{j=1}^{\infty}\left(  0,\omega\left(
w,e_{j}\right)  \right)  \otimes\left(  0,\omega\left(  w,e_{j}\right)
\right)
\end{align*}
where $\left\{  e_{j}\right\}  _{j=1}^{\infty}$ is any orthonormal basis for
$H.$ Both $\gamma$ and $\chi$ are well defined because of Lemma \ref{l.3.16}.
\end{nota}

\begin{nota}
[Left differentials]\label{n.3.28} Suppose $f:G\rightarrow\mathbb{C}$ is a
smooth cylinder function. For $g\in G$ and $h,h_{1},...,h_{n}\in\mathfrak{g}$,
$n\in\mathbb{N}$ let%
\begin{align}
&  \left(  D^{0}f\right)  \left(  g\right)  =f\left(  g\right)  \text{
and}\nonumber\\
&  \left(  D^{n}f\right)  \left(  g\right)  \left(  h_{1}\otimes...\otimes
h_{n}\right)  =\tilde{h}_{1}...\tilde{h}_{n}f\left(  g\right)  ,
\label{e.3.47}%
\end{align}
where $\tilde{h}f$ is given as in Equation \eqref{e.3.10} or Equation
\eqref{e.3.12}. We will write $Df$ for $D^{1}f.$
\end{nota}

\begin{prop}
\label{p.3.29}Let $\left\{  e_{j}\right\}  _{j=1}^{\infty}$ and $\left\{
f_{\ell}\right\}  _{\ell=1}^{d}$ be orthonormal bases for $H$ and $\mathbf{C}$
respectively. Then for any smooth cylinder function, $f:G\rightarrow
\mathbb{C},$%
\begin{equation}
Lf\left(  g\right)  :=\sum_{j=1}^{\infty}\left[  \widetilde{\left(
e_{j},0\right)  }^{2}f\right]  \left(  g\right)  +\sum_{\ell=1}^{d}\left[
\widetilde{\left(  0,f_{\ell}\right)  }^{2}f\right]  \left(  g\right)
\label{e.3.48}%
\end{equation}
is well defined. Moreover, if $f=F\circ\pi_{P}$, $\partial_{h}$ is as in
Notation \ref{n.3.5} for all $h\in\mathfrak{g}_{CM},$%
\begin{equation}
\Delta_{H}f\left(  g\right)  :=\sum_{j=1}^{\infty}\partial_{\left(
e_{j},0\right)  }^{2}f\left(  g\right)  =\left(  \Delta_{PH}F\right)  \left(
Pw,c\right)  \label{e.3.49}%
\end{equation}
and%
\begin{equation}
\Delta_{\mathbf{C}}f\left(  g\right)  :=\sum_{\ell=1}^{d}\left[
\partial_{\left(  0,f_{\ell}\right)  }^{2}f\right]  \left(  g\right)  =\left(
\Delta_{\mathbf{C}}F\right)  \left(  Pw,c\right)  , \label{e.3.50}%
\end{equation}
then%
\begin{equation}
Lf\left(  g\right)  =\left(  \Delta_{H}f+\Delta_{\mathbf{C}}f\right)  \left(
g\right)  +f^{\prime\prime}\left(  g\right)  \left(  \gamma\left(  g\right)
+\frac{1}{4}\chi\left(  g\right)  \right) . \label{e.3.51}%
\end{equation}

\end{prop}

\begin{proof}
The proof of the second equality in Eq. (\ref{e.3.49}) is straightforward and
will be left to the reader. Recall from Equation \eqref{e.3.12} that
\begin{equation}
\widetilde{\left(  e_{j},0\right)  }f\left(  g\right)  =f^{\prime}\left(
g\right)  \left(  e_{j},\frac{1}{2}\omega\left(  w,e_{j}\right)  \right)  .
\label{e.3.52}%
\end{equation}
Applying $\widetilde{\left(  e_{j},0\right)  }$ to both sides of Equation
\eqref{e.3.52} gives%
\begin{align}
\widetilde{\left(  e_{j},0\right)  }^{2}f\left(  g\right)   &  =f^{\prime
\prime}\left(  g\right)  \left(  \left(  e_{j},\frac{1}{2}\omega\left(
w,e_{j}\right)  \right)  \otimes\left(  e_{j},\frac{1}{2}\omega\left(
w,e_{j}\right)  \right)  \right) \label{e.3.53}\\
&  =f^{\prime\prime}\left(  g\right)  \left(  \left(  e_{j},0\right)
\otimes\left(  e_{j},0\right)  \right)  +f^{\prime\prime}\left(  g\right)
\left(  \left(  0,\omega\left(  w,e_{j}\right)  \right)  \otimes\left(
e_{j},0\right)  \right) \nonumber\\
&  +\frac{1}{4}f^{\prime\prime}\left(  g\right)  \left(  \left(
0,\omega\left(  w,e_{j}\right)  \right)  \otimes\left(  0,\omega\left(
w,e_{j}\right)  \right)  \right)  \label{e.3.54}%
\end{align}
wherein we have used,
\[
\partial_{e_{j}}\omega\left(  \cdot,e_{j}\right)  =\omega\left(  e_{j}%
,e_{j}\right)  =0.
\]

Summing Equation \eqref{e.3.54} on $j$ shows,%
\begin{align*}
\sum_{j=1}^{\infty}\left[  \widetilde{\left(  e_{j},0\right)  }^{2}f\right]
\left(  g\right)   &  =\sum_{j=1}^{\infty}f^{\prime\prime}\left(  g\right)
\left(  \left(  e_{j},0\right)  \otimes\left(  e_{j},0\right)  \right)
+f^{\prime\prime}\left(  g\right)  \left(  \gamma\left(  g\right)  +\frac
{1}{4}\chi\left(  g\right)  \right) \\
&  =\sum_{j=1}^{\infty}\partial_{\left(  e_{j},0\right)  }^{2}f\left(
g\right)  +f^{\prime\prime}\left(  g\right)  \left(  \gamma\left(  g\right)
+\frac{1}{4}\chi\left(  g\right)  \right)  .
\end{align*}
The formula in Equation \eqref{e.3.51} for $Lf$ is now easily verified and
this shows that $Lf$ is independent of the choice of orthonormal bases for $H$
and $\mathbf{C}$ appearing in Equation \eqref{e.3.48}.
\end{proof}

\section{Brownian Motion and Heat Kernel Measures\label{s.4}}

For the Hilbert space stochastic calculus background needed for this section,
see M\'{e}tivier \cite{Metivier82}. For the background on It\^{o} integral
relative to an abstract Wiener space--valued Brownian motion, see Kuo
\cite[pages 188-207]{Kuo75} (especially Theorem 5.1), Kusuoka and Stroock
\cite[p. 5]{KS1991}, and the appendix in \cite{Driver1997b}.

Suppose now that $\left(  B\left(  t\right)  ,B_{0}\left(  t\right)  \right)
$ is a smooth curve in $\mathfrak{g}_{CM}$ with $\left(  B\left(  0\right)
,B_{0}\left(  0\right)  \right)  =\left(  0,0\right)  $ and consider solving,
for $g\left(  t\right)  =\left(  w\left(  t\right)  ,c\left(  t\right)
\right)  \in G_{CM},$ the differential equation%
\begin{equation}
\left(  \dot{w}\left(  t\right)  ,\dot{c}\left(  t\right)  \right)  =\dot
{g}\left(  t\right)  =l_{g\left(  t\right)  \ast}\left(  \dot{B}\left(
t\right)  ,\dot{B}_{0}\left(  t\right)  \right)  \text{ with }g\left(
0\right)  =\left(  0,0\right)  . \label{e.4.1}%
\end{equation}
By Equation \eqref{e.3.11}, it follows that
\[
\left(  \dot{w}\left(  t\right)  ,\dot{c}\left(  t\right)  \right)  =\left(
\dot{B}\left(  t\right)  ,\dot{B}_{0}\left(  t\right)  +\frac{1}{2}%
\omega\left(  w\left(  t\right)  ,\dot{B}\left(  t\right)  \right)  \right)
\]
and therefore the solution to Equation \eqref{e.4.1} is given by%
\begin{equation}
g\left(  t\right)  =\left(  w\left(  t\right)  ,c\left(  t\right)  \right)
=\left(  B\left(  t\right)  ,B_{0}\left(  t\right)  +\frac{1}{2}\int_{0}%
^{t}\omega\left(  B\left(  \tau\right)  ,\dot{B}\left(  \tau\right)  \right)
d\tau\right)  . \label{e.4.2}%
\end{equation}
Below in subsection \ref{s.4.2}, we will replace $B$ and $B_{0}$ by Brownian
motions and use this to define a Brownian motion on $G.$

\subsection{A quadratic integral\label{s.4.1}}

Let $\left\{  \left(  B\left(  t\right)  ,B_{0}\left(  t\right)  \right)
\right\}  _{t\geqslant0}$ be a \textbf{Brownian motion} on $\mathfrak{g}$ with
variance determined by
\begin{multline*}
\mathbb{E}\left[  \left\langle \left(  B\left(  s\right)  ,B_{0}\left(
s\right)  \right)  ,\left(  A,a\right)  \right\rangle _{\mathfrak{g}_{CM}%
}\left\langle \left(  B\left(  t\right)  ,B_{0}\left(  t\right)  \right)
,\left(  C,c\right)  \right\rangle _{\mathfrak{g}_{CM}}\right] \\
=\operatorname{Re}\left\langle \left(  A,a\right)  ,\left(  C,c\right)
\right\rangle _{\mathfrak{g}_{CM}}\min\left(  s,t\right)
\end{multline*}
for all $s,t\in\lbrack0,\infty),$ $A,C\in H_{\ast}$ and $a,c\in\mathbf{C}.$
Also let $\left\{  e_{j}\right\}  _{j=1}^{\infty}\subset H_{\ast}$ be an
orthonormal basis for $H.$ For $n\in\mathbb{N},$ define $P_{n}\in
\operatorname*{Proj}\left(  W\right)  $ as in Notation \ref{n.3.24}, i.e.
\begin{equation}
P_{n}\left(  w\right)  =\sum_{j=1}^{n}\left\langle w,e_{j}\right\rangle
_{H}e_{j}=\sum_{j=1}^{n}\ell_{j}\left(  w\right)  e_{j}\text{ for all }w\in W.
\label{e.4.3}%
\end{equation}

\begin{prop}
\label{p.4.1}For each $n,$ let $M_{t}^{n}:=\int_{0}^{t}\omega\left(  B\left(
\tau\right)  ,dP_{n}B\left(  \tau\right)  \right)  .$ Then

\begin{enumerate}
\item $\left\{  M_{t}^{n}\right\}  _{t\geqslant0}$ is an $L^{2}$--martingale
and there exists an $L^{2}$--martingale, $\left\{  M_{t}\right\}
_{t\geqslant0}$ with values in $\mathbf{C}$ such that
\begin{equation}
\lim_{n\rightarrow\infty}\mathbb{E}\left[  \max_{t\leqslant T}\left\Vert
M_{t}-M_{t}^{n}\right\Vert _{\mathbf{C}}^{2}\right]  =0\text{ for all
}T<\infty. \label{e.4.4}%
\end{equation}

\item The quadratic variation of $M$ is given by;
\begin{equation}
\left\langle M\right\rangle _{t}=\int_{0}^{t}\left\Vert \omega\left(  B\left(
\tau\right)  ,\cdot\right)  \right\Vert _{H^{\ast}\otimes\mathbf{C}}^{2}d\tau.
\label{e.4.5}%
\end{equation}

\item The square integrable martingale, $M_{t},$ is well defined independent
of the choice of the orthonormal basis, $\left\{  e_{j}\right\}
_{j=1}^{\infty}$ and hence will be denoted by $\int_{0}^{t}\omega\left(
B\left(  \tau\right)  , dB\left(  \tau\right)  \right)  .$

\item For each $p\in\lbrack1,\infty),$ $\left\{  M_{t}\right\}  _{t\geqslant
0}$ is $L^{p}$--integrable and there exists $c_{p}<\infty$ such that
\[
\mathbb{E}\left(  \sup_{0\leqslant t\leqslant T}\left\Vert M_{t}\right\Vert
_{\mathbf{C}}^{p}\right)  \leqslant c_{p}T^{p}<\infty\text{ for all
}0\leqslant T<\infty.
\]
(This estimate will be considerably generalized in Proposition \ref{p.4.13} below.)
\end{enumerate}
\end{prop}

\begin{proof}
1. For $P\in\operatorname*{Proj}\left(  W\right)  $ let $M_{t}^{P}:=\int
_{0}^{t}\omega\left(  B\left(  \tau\right)  ,dPB\left(  \tau\right)  \right)
.$ Let $P,Q\in\operatorname*{Proj}\left(  W\right)  $ and choose an
orthonormal basis, $\left\{  v_{l}\right\}  _{l=1}^{N}$ for
$\operatorname*{Ran}\left(  P\right)  +\operatorname*{Ran}\left(  Q\right)  .$
We then have%
\begin{align}
\mathbb{E}\left[  \left\Vert M_{T}^{P}-M_{T}^{Q}\right\Vert _{\mathbf{C}}%
^{2}\right]   &  =\mathbb{E}\int_{0}^{T}\sum_{l=1}^{N}\left\Vert \omega\left(
B\left(  \tau\right)  ,\left(  P-Q\right)  v_{l}\right)  \right\Vert
_{\mathbf{C}}^{2}d\tau\nonumber\\
&  =\mathbb{E}\int_{0}^{T}\sum_{l=1}^{\infty}\left\Vert \omega\left(  B\left(
\tau\right)  ,\left(  P-Q\right)  e_{l}\right)  \right\Vert _{\mathbf{C}}%
^{2}d\tau\label{e.4.6}\\
&  =\int_{0}^{T}\sum_{l=1}^{\infty}\sum_{k=1}^{\infty}\left\Vert \omega\left(
e_{k},\left(  P-Q\right)  e_{l}\right)  \right\Vert _{\mathbf{C}}^{2}\tau
d\tau\nonumber\\
&  =\frac{T^{2}}{2}\sum_{l=1}^{\infty}\sum_{k=1}^{\infty}\left\Vert
\omega\left(  e_{k},\left(  P-Q\right)  e_{l}\right)  \right\Vert
_{\mathbf{C}}^{2}. \label{e.4.7}%
\end{align}
Taking $P=P_{n}$ and $Q=P_{m}$ with $m\leqslant n$ in Eq. (\ref{e.4.7}) allows
us to conclude that
\[
\mathbb{E}\left[  \left\Vert M_{T}^{n}-M_{T}^{m}\right\Vert _{\mathbf{C}}%
^{2}\right]  =\frac{T^{2}}{2}\sum_{j=m+1}^{n}\sum_{l=1}^{\infty}\left\Vert
\omega\left(  e_{l},e_{j}\right)  \right\Vert _{\mathbf{C}}^{2}\rightarrow
0\text{ as }m,n\rightarrow\infty
\]
because $\left\Vert \omega\right\Vert _{2}^{2}<\infty$ by Proposition
\ref{p.3.14}. Since the space of continuous $L^{2}$--martingales on $\left[
0,T\right]  $ is complete in the norm, $N\rightarrow\mathbb{E}\left\Vert
N_{T}\right\Vert _{\mathbf{C}}^{2}$ and, by Doob's maximal inequality
(\cite[Proposition 7.16]{Kall}), there exists $c<\infty$ such that
\[
\mathbb{E}\left[  \max_{t\leqslant T}\left\Vert N_{t}\right\Vert _{\mathbf{C}%
}^{p}\right]  \leqslant c\mathbb{E}\left\Vert N_{T}\right\Vert _{\mathbf{C}%
}^{p},
\]
it follows that there exists a square integrable $\mathbf{C}$--valued
martingale, $\left\{  M_{t}\right\}  _{t\geqslant0},$ such that Eq.
(\ref{e.4.4}) holds.

2. Since the quadratic variation of $M^{n}$ is given by%
\[
\left\langle M^{n}\right\rangle _{t}=\int_{0}^{t}\left\Vert \omega\left(
B\left(  \tau\right)  ,dP_{n}B\left(  \tau\right)  \right)  \right\Vert
_{\mathbf{C}}^{2}=\int_{0}^{t}\sum_{l=1}^{n}\left\Vert \omega\left(  B\left(
\tau\right)  , e_{l}\right)  \right\Vert _{\mathbf{C}}^{2}d\tau
\]
and%
\begin{align*}
\mathbb{E}\left[  \left\vert \left\langle M\right\rangle _{t}-\left\langle
M^{n}\right\rangle _{t}\right\vert \right]   &  \leqslant\sqrt{\mathbb{E}%
\left[  \left\langle M-M^{n}\right\rangle _{t}\right]  \cdot\mathbb{E}\left[
\left\langle M+M^{n}\right\rangle _{t}\right]  }\\
&  =\sqrt{\mathbb{E}\left\Vert M_{t}-M_{t}^{n}\right\Vert _{\mathbf{C}}%
^{2}\cdot\mathbb{E}\left\Vert M_{t}+M_{t}^{n}\right\Vert _{\mathbf{C}}^{2}%
}\rightarrow0\text{ as }n\rightarrow\infty,
\end{align*}
Eq. (\ref{e.4.5}) easily follows.

3. Suppose now that $\left\{  e_{j}^{\prime}\right\}  _{j=1}^{\infty}\subset
H_{\ast}$ is another orthonormal basis for $H$ and $P_{n}^{\prime
}:W\rightarrow H_{\ast}$ are the corresponding orthogonal projections. Taking
$P=P_{n}$ and $P^{\prime}=P_{n}^{\prime}$ in Eq. (\ref{e.4.7}) gives,%
\begin{equation}
\mathbb{E}\left\Vert M_{T}^{P_{n}}-M_{T}^{P_{n}^{\prime}}\right\Vert
_{\mathbf{C}}^{2}=\frac{T^{2}}{2}\sum_{l=1}^{\infty}\sum_{k=1}^{\infty
}\left\Vert \omega\left(  e_{k},\left(  P_{n}-P_{n}^{\prime}\right)
e_{l}\right)  \right\Vert _{\mathbf{C}}^{2}. \label{e.4.8}%
\end{equation}
Since%
\begin{align*}
\sum_{l=1}^{\infty}\sum_{k=1}^{\infty}\left\Vert \omega\left(  e_{k}%
,P_{n}^{\prime}e_{l}\right)  -\omega\left(  e_{k},e_{l}\right)  \right\Vert
_{\mathbf{C}}^{2}  &  =\sum_{l=1}^{\infty}\sum_{k=1}^{\infty}\left\Vert
\omega\left(  e_{k},\left(  P_{n}^{\prime}-I\right)  e_{l}\right)  \right\Vert
_{\mathbf{C}}^{2}\\
&  =\sum_{l=1}^{\infty}\sum_{k=1}^{\infty}\left\Vert \omega\left(
e_{k}^{\prime},\left(  P_{n}^{\prime}-I\right)  e_{l}^{\prime}\right)
\right\Vert _{\mathbf{C}}^{2}\\
&  =\sum_{l=n+1}^{\infty}\sum_{k=1}^{\infty}\left\Vert \omega\left(
e_{k}^{\prime},e_{l}^{\prime}\right)  \right\Vert _{\mathbf{C}}^{2}%
\rightarrow0\text{ as }n\rightarrow\infty
\end{align*}
and similarly but more easily, $\sum_{l=1}^{\infty}\sum_{k=1}^{\infty
}\left\Vert \omega\left(  e_{k},P_{n}e_{l}\right)  -\omega\left(  e_{k}%
,e_{l}\right)  \right\Vert _{\mathbf{C}}^{2}\rightarrow0$ as $n\rightarrow
\infty,$ we may pass to the limit in Eq. (\ref{e.4.8}) to learn that
$\lim_{n\rightarrow\infty}\mathbb{E}\left\Vert M_{T}^{P_{n}}-M_{T}%
^{P_{n}^{\prime}}\right\Vert _{\mathbf{C}}^{2}=0.$

4. By Jensen's inequality
\begin{align*}
\left(  \int_{0}^{T}\left\Vert \omega\left(  B\left(  s\right)  ,\cdot\right)
\right\Vert _{H^{\ast}\otimes\mathbf{C}}^{2}ds\right)  ^{p/2}  &
=T^{p/2}\left(  \int_{0}^{T}\left\Vert \omega\left(  B\left(  s\right)
,\cdot\right)  \right\Vert _{H^{\ast}\otimes\mathbf{C}}^{2}\frac{ds}%
{T}\right)  ^{p/2}\\
&  \leqslant T^{p/2}\int_{0}^{T}\left\Vert \omega\left(  B\left(  s\right)
,\cdot\right)  \right\Vert _{H^{\ast}\otimes\mathbf{C}}^{p}\frac{ds}{T}\\
&  =T^{\frac{p}{2}-1}\int_{0}^{T}\left\Vert \omega\left(  B\left(  s\right)
,\cdot\right)  \right\Vert _{H^{\ast}\otimes\mathbf{C}}^{p}ds.
\end{align*}

Combining this estimate with Equation \eqref{e.3.28} and then applying either
Skorohod's or Fernique's inequality (see Equations \eqref{e.2.3} or
\eqref{e.2.4}) shows
\begin{align}
\mathbb{E}\left[  \left\langle M\right\rangle _{T}^{p/2}\right]   &  \leqslant
T^{\frac{p}{2}-1}\int_{0}^{T}\mathbb{E}\left\Vert \omega\left(  B\left(
s\right)  ,\cdot\right)  \right\Vert _{H^{\ast}\otimes\mathbf{C}}%
^{p}ds\nonumber\\
&  \leqslant T^{\frac{p}{2}-1}\int_{0}^{T}C_{2}^{p/2}\left\Vert \omega
\right\Vert _{0}^{p}\left\Vert B\left(  s\right)  \right\Vert _{W}%
^{p}ds\nonumber\\
&  \leqslant T^{\frac{p}{2}-1}C_{2}^{p/2}\left\Vert \omega\right\Vert _{0}%
^{p}\int_{W}\left\Vert y\right\Vert _{W}^{p}d\mu\left(  y\right)  \int_{0}%
^{T}s^{p/2}ds\nonumber\\
&  =T^{\frac{p}{2}-1}C_{2}^{p/2}\left\Vert \omega\right\Vert _{0}^{p}%
C_{p}\frac{T^{p/2+1}}{p/2+1}=c_{p}^{\prime}T^{p}. \label{e.4.9}%
\end{align}
As a consequence of the Burkholder-Davis-Gundy inequality (see for example
\cite[Corollary 6.3.1a on p.344]{StroockBook}, \cite[Appendix A.2]{Nualart95},
or \cite[p. 212]{Metivier82} and \cite[Theorem 17.7]{Kall} for the real case),
for any $p\geqslant2$ there exists $c_{p}^{\prime\prime}<\infty$ such that
\[
\mathbb{E}\left(  \sup_{0\leqslant t\leqslant T}\left\Vert M_{t}\right\Vert
_{\mathbf{C}}\right)  ^{p}\leqslant c_{p}^{\prime\prime}\mathbb{E}\left[
\left\langle M\right\rangle _{T}^{p/2}\right]  =c_{p}^{\prime\prime}%
c_{p}^{\prime}T^{p}=:c_{p}T^{p}.
\]

\end{proof}

\subsection{Brownian Motion on $G\left(  \omega\right)  $\label{s.4.2}}

Motivated by Eq. (\ref{e.4.2}) we have the following definition.

\begin{df}
\label{d.4.2}Let $\left(  B\left(  t\right)  ,B_{0}\left(  t\right)  \right)
$ be a $\mathfrak{g}$ -- valued Brownian motion as in subsection \ref{s.4.1}.
A \textbf{Brownian motion }on $G$ is the continuous $G$--valued process
defined by%
\begin{equation}
g\left(  t\right)  =\left(  B\left(  t\right)  ,B_{0}\left(  t\right)
+\frac{1}{2}\int_{0}^{t}\omega\left(  B\left(  \tau\right)  ,dB\left(
\tau\right)  \right)  \right)  . \label{e.4.10}%
\end{equation}
Further, for $T>0,$ let $\nu_{T}=\operatorname{Law}\left(  g\left(  T\right)
\right)  $ be a probability measure on $G.$ We refer to $\nu_{T}$ as the
\textbf{time }$T$ \textbf{heat kernel measure on }$G.$
\end{df}

\begin{rem}
\label{r.4.3} An alert reader may complain that we should use the Stratonovich
integral in Eq. (\ref{e.4.10}) rather than the It\^{o} integral. However,
these two integrals are equal since $\omega$ is a skew symmetric form
\begin{align*}
\int_{0}^{t}\omega\left(  B\left(  \tau\right)  ,\circ dB\left(  \tau\right)
\right)   &  =\int_{0}^{t}\omega\left(  B\left(  \tau\right)  ,dB\left(
\tau\right)  \right)  +\frac{1}{2}\int_{0}^{t}\omega\left(  dB\left(
\tau\right)  ,dB\left(  \tau\right)  \right) \\
&  =\int_{0}^{t}\omega\left(  B\left(  \tau\right)  ,dB\left(  \tau\right)
\right)  .
\end{align*}

\end{rem}

\begin{thm}
[The generator of $g\left(  t\right)  $]\label{t.4.4}The generator of
$g\left(  t\right)  $ is the operator $L$ defined in Proposition \ref{p.3.29}.
More precisely, if $f:G\rightarrow\mathbb{C}$ is a smooth cylinder function,
then
\begin{equation}
d\left[  f\left(  g\left(  t\right)  \right)  \right]  =f^{\prime}\left(
g\left(  t\right)  \right)  dg\left(  t\right)  +\frac{1}{2}Lf(g\left(
t\right)  dt \label{e.4.11}%
\end{equation}
where $L$ is given in Proposition \ref{p.3.29}, $f^{\prime}$ is defined as in
Notation \ref{n.3.5} and
\[
dg\left(  t\right)  =\left(  dB\left(  t\right)  ,dB_{0}\left(  t\right)
+\frac{1}{2}\omega\left(  B\left(  t\right)  ,dB\left(  t\right)  \right)
\right)  .
\]

\end{thm}

\begin{proof}
Let us begin by observing that%
\begin{align}
dg\left(  t\right)  \otimes dg\left(  t\right)   &  =\left(  dB\left(
t\right)  ,dB_{0}\left(  t\right)  +\frac{1}{2}\omega\left(  B\left(
t\right)  ,dB\left(  t\right)  \right)  \right)  ^{\otimes2}\nonumber\\
&  =\left[  \left(  dB\left(  t\right)  ,\frac{1}{2}\omega\left(  B\left(
t\right)  ,dB\left(  t\right)  \right)  \right)  +\left(  0,dB_{0}\left(
t\right)  \right)  \right]  ^{\otimes2}\nonumber\\
&  =\sum_{j=1}^{\infty}\left(  e_{j},\frac{1}{2}\omega\left(  B\left(
t\right)  ,e_{j}\right)  \right)  ^{\otimes2}dt+\sum_{\ell=1}^{d}\left(
0,f_{\ell}\right)  ^{\otimes2}dt\label{e.4.12}%
\end{align}
where $\left\{  f_{\ell}\right\}  _{\ell=1}^{d}$ is an orthonormal basis for
$\mathbf{C}$ and $\left\{  e_{j}\right\}  _{j=1}^{\infty}$ is an orthonormal
basis for $H.$ Hence, as a consequence of It\^{o}'s formula, we have%
\begin{align*}
d\left[  f\left(  g\left(  t\right)  \right)  \right]  = &  f^{\prime}\left(
g\left(  t\right)  \right)  dg\left(  t\right)  +\frac{1}{2}f^{\prime\prime
}\left(  g\left(  t\right)  \right)  \left(  dg\left(  t\right)  \otimes
dg\left(  t\right)  \right)  \\
= &  f^{\prime}\left(  g\left(  t\right)  \right)  dg\left(  t\right)  \\
&  +\frac{1}{2}f^{\prime\prime}\left(  g\left(  t\right)  \right)  \sum
_{j=1}^{\infty}\left(  e_{j},\frac{1}{2}\omega\left(  B\left(  t\right)
,e_{j}\right)  \right)  ^{\otimes2}dt+\frac{1}{2}f^{\prime\prime}\left(
g\left(  t\right)  \right)  \sum_{\ell=1}^{d}\left(  0,f_{\ell}\right)
^{\otimes2}dt\\
= &  f^{\prime}\left(  g\left(  t\right)  \right)  dg\left(  t\right)
+\frac{1}{2}\sum_{j=1}^{\infty}\left(  \widetilde{\left(  e_{j},0\right)
}^{2}f\right)  \left(  g\left(  t\right)  \right)  dt\\
&  +\frac{1}{2}\sum_{\ell=1}^{d}\left(  \widetilde{\left(  0,f_{\ell}\right)
}^{2}f\right)  \left(  g\left(  t\right)  \right)  dt\\
= &  f^{\prime}\left(  g\left(  t\right)  \right)  \left(  dg\left(  t\right)
\right)  +\frac{1}{2}Lf\left(  g\left(  t\right)  \right)  dt.
\end{align*}

\end{proof}

For the next corollary, let $P\in\operatorname*{Proj}\left(  W\right)  $ as in
Equation \eqref{e.3.42}, $F\in C^{2}\left(  PH\times\mathbf{C},\mathbb{C}%
\right)  ,$ and $f=F\circ\pi_{P}:G\rightarrow\mathbb{C}$ be a cylinder
function where $P\in\operatorname*{Proj}\left(  W\right)  .$ We will further
suppose there exist $0<K,p<\infty$ such that%
\begin{equation}
\left\vert F\left(  h,c\right)  \right\vert +\left\Vert F^{\prime}\left(
h,c\right)  \right\Vert +\left\Vert F^{\prime\prime}\left(  h,c\right)
\right\Vert \leqslant K\left(  1+\left\Vert h\right\Vert _{H}+\left\Vert
c\right\Vert _{\mathbf{C}}\right)  ^{p} \label{e.4.13}%
\end{equation}
for any $h\in PH$ and $c\in\mathbf{C}.$ Further let $\left\{  f_{\ell
}\right\}  _{\ell=1}^{d}$ be an orthonormal basis for $\mathbf{C}$ and extend
$\left\{  e_{j}\right\}  _{j=1}^{n}$\textbf{ }to an orthonormal basis,
$\left\{  e_{j}\right\}  _{j=1}^{\infty},$ for $H.$

\begin{cor}
\label{c.4.5}If $f:G\rightarrow\mathbb{C}$ is a cylinder function as above,
then%
\begin{equation}
\mathbb{E}\left[  f\left(  g\left(  T\right)  \right)  \right]  =f\left(
\mathbf{e}\right)  +\frac{1}{2}\int_{0}^{T}\mathbb{E}\left[  \left(
Lf\right)  \left(  g\left(  t\right)  \right)  \right]  dt, \label{e.4.14}%
\end{equation}
i.e.%
\begin{equation}
\nu_{T}\left(  f\right)  =f\left(  \mathbf{e}\right)  +\frac{1}{2}\int_{0}%
^{T}\mathbb{\nu}_{t}\left(  Lf\right)  dt. \label{e.4.15}%
\end{equation}
In other words, $\nu_{t}$ weakly solves the heat equation
\[
\partial_{t}\nu_{t}=\frac{1}{2}L\nu_{t}\text{ with }\lim_{t\downarrow0}\nu
_{t}=\delta_{\mathbf{e}}.
\]

\end{cor}

\begin{proof}
Integrating Equation \eqref{e.4.11} shows
\begin{equation}
f\left(  g\left(  T\right)  \right)  =f\left(  \mathbf{e}\right)  +N_{T}%
+\frac{1}{2}\int_{0}^{T}Lf\left(  g\left(  \tau\right)  \right)
d\tau\label{e.4.16}%
\end{equation}
where%
\[
N_{t}:=\int_{0}^{t}f^{\prime}\left(  g\left(  \tau\right)  \right)  dg\left(
\tau\right)  =\int_{0}^{t}f^{\prime}\left(  g\left(  \tau\right)  \right)
\left(  dB\left(  \tau\right)  ,dB_{0}\left(  \tau\right)  +\frac{1}%
{2}dM_{\tau}\right)
\]
and $M_{t}=\int_{0}^{t}\omega\left(  B\left(  \tau\right)  ,dB\left(
\tau\right)  \right)  .$ Using Eqs. \eqref{e.4.12} and \eqref{e.4.13} there
exists $C=C\left(  P,\left\Vert \omega_{0}\right\Vert \right)  <\infty$ such
that
\begin{align*}
d\left\langle N\right\rangle _{t}  &  :=\left\vert dN_{t}\right\vert
^{2}=\left\langle f^{\prime}\left(  g_{t}\right)  \otimes f^{\prime}\left(
g_{t}\right)  ,dg_{t}\otimes dg_{t}\right\rangle \\
&  =\sum_{j=1}^{\infty}\left\vert f^{\prime}\left(  g\left(  t\right)
\right)  \left(  e_{j},\frac{1}{2}\omega\left(  B\left(  t\right)
,e_{j}\right)  \right)  \right\vert ^{2}dt+\sum_{\ell=1}^{d}\left\vert
f^{\prime}\left(  g\left(  t\right)  \right)  \left(  0,f_{\ell}\right)
\right\vert ^{2}dt\\
&  =\sum_{j=1}^{n}\left\vert f^{\prime}\left(  g\left(  t\right)  \right)
\left(  e_{j},\frac{1}{2}\omega\left(  B\left(  t\right)  ,e_{j}\right)
\right)  \right\vert ^{2}dt+\sum_{\ell=1}^{d}\left\vert f^{\prime}\left(
g\left(  t\right)  \right)  \left(  0,f_{\ell}\right)  \right\vert ^{2}dt\\
&  \leqslant C_{1}\left(  P,\left\Vert \omega_{0}\right\Vert \right)  \left(
1+\left\Vert PB\left(  t\right)  \right\Vert _{H}+\left\Vert B_{0}\left(
t\right)  \right\Vert _{\mathbf{C}}\right)  ^{p}\left(  \left\Vert B\left(
t\right)  \right\Vert _{W}^{2}+1\right)  dt\\
&  \leqslant C\left(  1+\left\Vert B\left(  t\right)  \right\Vert
_{W}+\left\Vert B_{0}\left(  t\right)  \right\Vert _{\mathbf{C}}\right)
^{p+2}dt,
\end{align*}
wherein we have used Equation \eqref{e.3.43} for the last inequality. From
this inequality and either of Equations \eqref{e.2.3} or \eqref{e.2.4}, we
find%
\[
\mathbb{E}\left[  \left\langle N\right\rangle _{T}\right]  \leqslant C\int
_{0}^{T}\mathbb{E}\left(  1+\left\Vert B\left(  t\right)  \right\Vert
_{W}+\left\Vert B_{0}\left(  t\right)  \right\Vert _{\mathbf{C}}\right)
^{p+2}dt<\infty
\]
and hence that $N_{t}$ is a square integrable martingale. Therefore we may
take the expectation of Equation \eqref{e.4.16} which implies Equation \eqref{e.4.14}.
\end{proof}

\subsection{Finite Dimensional Approximations\label{s.4.3}}

\begin{prop}
\label{p.4.6}Let $\left\{  P_{n}\right\}  _{n=1}^{\infty}\subset
\operatorname*{Proj}\left(  W\right)  $ be as in Eq. \eqref{e.4.3} and
\begin{equation}
B_{n}\left(  t\right)  :=P_{n}B\left(  t\right)  \in P_{n}H\subset H\subset W.
\label{e.4.17}%
\end{equation}
Then%
\begin{equation}
\lim_{n\rightarrow\infty}\mathbb{E}\left[  \max_{0\leqslant t\leqslant
T}\left\Vert B\left(  t\right)  -B_{n}\left(  t\right)  \right\Vert _{W}%
^{p}\right]  =0\text{ for all }p\in\lbrack1,\infty), \label{e.4.18}%
\end{equation}
and%
\begin{equation}
\lim_{n\rightarrow\infty}\max_{0\leqslant t\leqslant T}\left\Vert B\left(
t\right)  -B_{n}\left(  t\right)  \right\Vert _{W}=0\text{ a.s.}
\label{e.4.19}%
\end{equation}

\end{prop}

\begin{proof}
Let $\left\{  w_{k}\right\}  _{k=1}^{\infty}\subset W$ be a countable dense
set and for each $k\in\mathbb{N},$ choose $\varphi_{k}\in W^{\ast}$ such that
$\left\Vert \varphi_{k}\right\Vert _{W^{\ast}}=1$ and $\varphi_{k}\left(
w_{k}\right)  =\left\Vert w_{k}\right\Vert _{W}.$ We then have,%
\[
\left\Vert w\right\Vert _{W}=\sup_{k}\left\vert \varphi_{k}\left(  w\right)
\right\vert =\sup\operatorname{Re}\varphi_{k}\left(  w\right)  \text{ for all
}w\in W.
\]
By \cite[Theorem 3.5.7]{Bog98} with $A=I$, if $\varepsilon_{n}\left(
t\right)  :=B\left(  t\right)  -B_{n}\left(  t\right)  ,$ then
\begin{equation}
\lim_{n\rightarrow\infty}\mathbb{E}\left\Vert \varepsilon_{n}\left(  T\right)
\right\Vert _{W}^{p}=0\text{ for all }p\in\lbrack1,\infty). \label{e.4.20}%
\end{equation}
Since $\left\{  \varphi_{k}\left(  \varepsilon_{n}\left(  t\right)  \right)
\right\}  _{t\geqslant0}$ is (up to a multiplicative factor) a standard
Brownian motion, $\left\{  \left\vert \varphi_{k}\left(  \varepsilon
_{n}\left(  t\right)  \right)  \right\vert \right\}  _{t\geqslant0}$ is a
submartingale for each $k\in\mathbb{N}$ and therefore so is $\left\{
\left\Vert \varepsilon_{n}\left(  t\right)  \right\Vert =\sup_{k}\left\vert
\varphi_{k}\left(  \varepsilon_{n}\left(  t\right)  \right)  \right\vert
\right\}  _{t\geqslant0}.$ Hence, according to Doob's inequality, for each
$p\in\lbrack1,\infty)$ there exists $C_{p}<\infty$ such that
\begin{equation}
\mathbb{E}\left\vert \max_{t\leqslant T}\left\Vert \varepsilon_{n}\left(
t\right)  \right\Vert _{W}\right\vert ^{p}\leqslant C_{p}\mathbb{E}\left\Vert
\varepsilon_{n}\left(  T\right)  \right\Vert _{W}^{p}. \label{e.4.21}%
\end{equation}
Combining Equation \eqref{e.4.21} with Equation \eqref{e.4.20} proves Equation
\eqref{e.4.18}. Equation \eqref{e.4.19} now follows from Equation
\eqref{e.4.18} and \cite[Proposition 2.11]{DaPratoBook1992}. To apply this
proposition, let $E$ be the Banach space, $C\left(  \left[  0,T\right]  ,
W\right)  $ equipped with the sup-norm, and let $\xi_{k}:=\ell_{k}\left(
B\left(  \cdot\right)  \right)  e_{k}\in E$ for all $k\in\mathbb{N}.$
\end{proof}

\begin{lem}
[Finite Dimensional Approximations to $g\left(  t\right)  $]\label{l.4.7}For
$P\in\operatorname*{Proj}\left(  W\right)  ,$ $Q:=I_{W}-P,$ let $g_{P}\left(
t\right)  $ be the Brownian motion on $G_{P}$ defined by%
\[
g_{P}\left(  t\right)  :=\left(  PB\left(  t\right)  ,B_{0}\left(  t\right)
+\frac{1}{2}\int_{0}^{t}\omega\left(  PB\left(  \tau\right)  ,PdB\left(
\tau\right)  \right)  \right)  .
\]
Then%
\begin{equation}
g\left(  t\right)  =g_{P}\left(  t\right)  \left(  QB\left(  t\right)
,\frac{1}{2}\int_{0}^{t}\left[  2\omega\left(  QB\left(  \tau\right)
,PdB\left(  \tau\right)  \right)  +\omega\left(  QB\left(  \tau\right)
,QdB\left(  \tau\right)  \right)  \right]  \right)  , \label{e.4.22}%
\end{equation}
and
\begin{equation}
g_{P}\left(  t\right)  ^{-1}\pi_{P}\left(  g\left(  t\right)  \right)
=\frac{1}{2}\left(  0,\int_{0}^{t}\left[  \omega\left(  B\left(  \tau\right)
,dB\left(  \tau\right)  \right)  -\omega\left(  PB\left(  \tau\right)
,PdB\left(  \tau\right)  \right)  \right]  \right)  . \label{e.4.23}%
\end{equation}
Also, if $\left\{  P_{n}\right\}  _{n=1}^{\infty}\subset\operatorname*{Proj}%
\left(  W\right)  $ are as in Eq. \eqref{e.4.3} and%
\begin{equation}
g_{n}\left(  t\right)  =g_{P_{n}}\left(  t\right)  =\left(  P_{n}B\left(
t\right)  ,B_{0}\left(  t\right)  +\frac{1}{2}\int_{0}^{t}\omega\left(
P_{n}B\left(  \tau\right)  ,dP_{n}B\left(  \tau\right)  \right)  \right)  ,
\label{e.4.24}%
\end{equation}
then
\begin{equation}
\lim_{n\rightarrow\infty}\mathbb{E}\left[  \max_{0\leqslant t\leqslant
T}\left\Vert g\left(  t\right)  -g_{n}\left(  t\right)  \right\Vert
_{\mathfrak{g}}^{p}\right]  =0 \label{e.4.25}%
\end{equation}
for all $1\leqslant p<\infty.$
\end{lem}

\begin{proof}
A simple computation shows%
\begin{align*}
l_{g_{P}\left(  t\right)  ^{-1}\ast}\circ dg_{P}\left(  t\right)   &  =\left(
\begin{array}
[c]{c}%
dPB(t),dB_{0}(t)+\frac{1}{2}\omega\left(  PB(t),PdB(t)\right)  \\
+\frac{1}{2}\omega\left(  -PB(t),PdB(t)\right)
\end{array}
\right)  \\
&  =\left(  dPB(t),dB_{0}(t)\right)  =d\left(  PB(t),B_{0}(t)\right)  .
\end{align*}
Hence it follows that $g_{P}$ solves the stochastic differential equation,%
\[
dg_{P}\left(  t\right)  =l_{g_{P}\left(  t\right)  \ast}\circ d\left(
PB(t),B_{0}(t)\right)  \text{ with }g_{P}\left(  0\right)  =0
\]
and therefore $g_{P}$ is a $G_{P}$--valued Brownian motion. The proof of the
equalities in Equations \eqref{e.4.22} and \eqref{e.4.23} follows by
elementary manipulations which are left to the reader.

In light of Equation \eqref{e.4.18} of Proposition \ref{p.4.6}, to prove the
last assertion we must show
\begin{equation}
\lim_{n\rightarrow\infty}\mathbb{E}\left[  \max_{0\leqslant t\leqslant
T}\left\vert M_{t}\left(  n\right)  \right\vert ^{p}\right]  =0,
\label{e.4.26}%
\end{equation}
where $M_{t}\left(  n\right)  $ is the local martingale defined by
\[
M_{t}\left(  n\right)  :=\int_{0}^{t}\left[  \omega\left(  B\left(
\tau\right)  , dB\left(  \tau\right)  \right)  -\omega\left(  B_{n}\left(
\tau\right)  , dB_{n}\left(  \tau\right)  \right)  \right]  .
\]
Since%
\begin{align*}
\left\langle M\left(  n\right)  \right\rangle _{T}  &  =\sum_{j=1}^{\infty
}\int_{0}^{T}\left\Vert \omega\left(  B\left(  \tau\right)  , e_{j}\right)
\right\Vert _{\mathbf{C}}^{2}d\tau+\sum_{j=1}^{n}\int_{0}^{T}\left\Vert
\omega\left(  B_{n}\left(  \tau\right)  , e_{j}\right)  \right\Vert
_{\mathbf{C}}^{2}d\tau\\
&  -2\sum_{j=1}^{n}\int_{0}^{T}\left\langle \omega\left(  B\left(
\tau\right)  , e_{j}\right)  , \omega\left(  B_{n}\left(  \tau\right)  ,
e_{j}\right)  \right\rangle _{\mathbf{C}}d\tau
\end{align*}
and
\[
\frac{2}{T^{2}}\mathbb{E}\left[  \left\langle M\left(  n\right)  \right\rangle
_{T}\right]  =\sum_{j,k=1}^{\infty}\left\Vert \omega\left(  e_{k}%
,e_{j}\right)  \right\Vert _{\mathbf{C}}^{2}-\sum_{k=1}^{n}\sum_{j=1}%
^{n}\left\Vert \omega\left(  e_{k},e_{j}\right)  \right\Vert _{\mathbf{C}}%
^{2}\rightarrow0
\]
as $n\rightarrow\infty,$ it follows by the Burkholder-Davis-Gundy inequalities
that $M\left(  n\right)  $ is a martingale and Equation \eqref{e.4.26} holds
for $p=2$ and hence for $p\in\left[  1,2\right]  .$

By Doob's maximal inequality (\cite[Proposition 7.16]{Kall}), to prove
Equation \eqref{e.4.26} for $p\geqslant2,$ it suffices to show $\lim
_{n\rightarrow\infty}\mathbb{E}\left[  \left\vert M_{T}\left(  n\right)
\right\vert ^{p}\right]  =0.$ However, $M_{T}\left(  n\right)  $ has It\^{o}'s
chaos expansion terminating at degree two and hence by a theorem of Nelson
(see \cite[Lemma 2 on p. 415]{Nelson73b} and \cite[pp. 216-217]{Nelson73c})
for each $j\in\mathbb{N}$ there exists $c_{j}<\infty$ such that $\mathbb{E}%
\left[  M_{T}^{2j}\left(  n\right)  \right]  \leqslant c_{j}\left[
\mathbb{E}M_{T}^{2}\left(  n\right)  \right]  ^{j}.$ (This result also follows
from Nelson's hypercontractivity for the Ornstein-Uhlenbeck operator.) This
clearly suffices to complete the proof of the theorem.
\end{proof}

\begin{lem}
\label{l.4.8}For all $P\in\operatorname*{Proj}\left(  W\right)  $ and $t>0,$
let $\nu_{t}^{P}:=\operatorname*{Law}\left(  g_{P}\left(  t\right)  \right)
.$ Then $\nu_{t}^{P}\left(  dx\right)  =p_{t}^{P}\left(  e,x\right)  dx,$
where $dx$ is the Riemannian volume measure (equal to a Haar measure)
$p_{t}^{P}\left(  x,y\right)  $ is the heat kernel on $G_{P}.$
\end{lem}

\begin{proof}
An application of Corollary \ref{c.4.5} with $G$ replaced by $G_{P}$ implies
that $\nu_{t}^{P}=\operatorname*{Law}\left(  g_{P}\left(  t\right)  \right)  $
is a weak solution to the heat equation on $G_{P}.$ The result now follows as
an application of \cite[Theorem 2.6 ]{Driver2003b}.
\end{proof}

\begin{cor}
\label{c.4.9}For any $T>0,$ the heat kernel measure, $\nu_{T},$ is invariant
under the inversion map, $g\mapsto g^{-1}$ for any $g \in G$.
\end{cor}

\begin{proof}
It is well known (see for example \cite[Proposition 3.1]{Driver1997c}) that
heat kernel measures based at the identity of a finite-dimensional Lie group
are invariant under inversion. Now suppose that $f:G\rightarrow\mathbb{R}$ is
a bounded continuous function. By passing to a subsequence if necessary, we
may assume that the sequence of $G$--valued random variables, $\left\{
g_{n}\left(  T\right)  \right\}  _{t\geqslant0},$ in Lemma \ref{l.4.7}
converges almost surely to $g\left(  T\right)  .$ Therefore by the dominated
convergence theorem,%
\[
\mathbb{E}f\left(  g\left(  T\right)  ^{-1}\right)  =\lim_{n\rightarrow\infty
}\mathbb{E}f\left(  g_{n}\left(  T\right)  ^{-1}\right)  =\lim_{n\rightarrow
\infty}\mathbb{E}f\left(  g_{n}\left(  T\right)  \right)  =\mathbb{E}f\left(
g\left(  T\right)  \right)  .
\]
This completes the proof because $\nu_{T}$ is the law of $g\left(  T\right)
.$
\end{proof}

We are now going to give exponential bounds which are much stronger than the
moment estimates in Equation \eqref{e.4.9} of Proposition \ref{p.4.1}. Before
doing so we need to recall the following result of Cameron--Martin and Kac,
\cite{Cameron-Martin1945,Kac1951}.

\begin{lem}
[Cameron--Martin and Kac]\label{l.4.10}Let $\left\{  b_{s}\right\}
_{s\geqslant0}$ be a one dimensional Brownian motion. Then for any $T>0$ and
$\lambda\in\lbrack0,\frac{\pi}{2T});$%
\begin{equation}
\mathbb{E}\left[  \exp\left(  \frac{\lambda^{2}}{2}\int_{0}^{T}b_{s}%
^{2}ds\right)  \right]  =\left[  \cos\left(  \lambda T\right)  \right]
^{-1/2}<\infty. \label{e.4.27}%
\end{equation}

\end{lem}

\begin{proof}
When $T=1,$ simply follow the proof of \cite[Equation (6.9) on p. 472]{IW89}
with $\lambda$ replaced by $-\lambda^{2}.$ For general $T>0,$ by a change of
variables and a Brownian motion scaling we have%
\[
\int_{0}^{T}b_{s}^{2}ds=T\int_{0}^{1}b_{tT}^{2}dt\overset{d}{=}T^{2}\int
_{0}^{1}b_{t}^{2}dt.
\]
Therefore,%
\begin{align}
\mathbb{E}\left[  \exp\left(  \frac{\lambda^{2}}{2}\int_{0}^{T}b_{s}%
^{2}ds\right)  \right]   &  =\mathbb{E}\left[  \exp\left(  \frac{\lambda
^{2}T^{2}}{2}\int_{0}^{1}b_{s}^{2}ds\right)  \right] \label{e.4.28}\\
&  =\cos^{-1/2}\left(  \sqrt{\lambda^{2}T^{2}}\right) \nonumber
\end{align}
provided that $\lambda\in\lbrack0,\frac{\pi}{2T}).$
\end{proof}

\begin{rem}
\label{r.4.11}For our purposes below, all we really need later from Lemma
\ref{l.4.10} is the qualitative statement that for $\lambda T>0$ sufficiently
small
\begin{equation}
\mathbb{E}\left[  \exp\left(  \frac{\lambda^{2}}{2}\int_{0}^{T}b_{s}%
^{2}ds\right)  \right]  =1+\frac{\lambda^{2}T^{2}}{4}+O\left(  \lambda
^{4}T^{4}\right)  . \label{e.4.29}%
\end{equation}
Instead of using Lemma \ref{l.4.10} we can derive this statement as an easy
consequence of the scaling identity in Equation \eqref{e.4.28} along with the
analyticity (use Fernique's theorem) of the function,%
\[
F\left(  z\right)  :=\mathbb{E}\left[  \exp\left(  z\int_{0}^{1}b_{s}%
^{2}ds\right)  \right]  \text{ for }\left\vert z\right\vert \text{ small.}%
\]

\end{rem}

\begin{prop}
\label{p.4.12}If $\left\{  N_{t}\right\}  _{t\geqslant0}$ is a continuous
local martingale such that $N_{0}=0.$ Then%
\begin{equation}
\mathbb{E}e^{\left\vert N_{t}\right\vert }\leqslant2\sqrt{\mathbb{E}\left[
e^{2\left\langle N\right\rangle _{t}}\right]  }. \label{e.4.30}%
\end{equation}

\end{prop}

\begin{proof}
By It\^{o}'s formula, we know that%
\[
Z_{t}:=e^{2N_{t}-\left\langle 2N\right\rangle _{t}/2}=e^{2N_{t}-2\left\langle
N\right\rangle _{t}}%
\]
is a non-negative local martingale. If $\left\{  \sigma_{n}\right\}
_{n=1}^{\infty}$ is a localizing sequence of stopping times for $Z,$ then, by
Fatou's lemma,
\[
\mathbb{E}\left[  Z_{t}|\mathcal{B}_{s}\right]  \leqslant\liminf
_{n\rightarrow\infty}\mathbb{E}\left[  Z_{t}^{\sigma_{n}}|\mathcal{B}%
_{s}\right]  =\liminf_{n\rightarrow\infty}Z_{s}^{\sigma_{n}}=Z_{s}.
\]
This shows that $Z$ is a supermartingale and in particular that $\mathbb{E}%
\left[  Z_{t}\right]  \leqslant\mathbb{E}Z_{0}=1.$

By the Cauchy-Schwarz inequality we find%
\begin{align}
\mathbb{E}\left[  e^{N_{t}}\right]   &  =\mathbb{E}\left[  e^{N_{t}%
-\left\langle N\right\rangle _{t}}e^{\left\langle N\right\rangle _{t}}\right]
\nonumber\\
&  \leqslant\sqrt{\mathbb{E}\left[  e^{2N_{t}-2\left\langle N\right\rangle
_{t}}\right]  \cdot\mathbb{E}\left[  e^{2\left\langle N\right\rangle _{t}%
}\right]  }=\sqrt{\mathbb{E}\left[  Z_{t}\right]  \cdot\mathbb{E}\left[
e^{2\left\langle N\right\rangle _{t}}\right]  }\nonumber\\
&  \leqslant\sqrt{\mathbb{E}\left[  e^{2\left\langle N\right\rangle _{t}%
}\right]  } \label{e.4.31}%
\end{align}

Applying this inequality with $N$ replaced by $-N$ and using $e^{\left\vert
x\right\vert }\leqslant e^{x}+e^{-x}$ easily give Equation \eqref{e.4.30}.
\end{proof}

\begin{prop}
\label{p.4.13}Let $n\in\mathbb{N},$ $T>0,$ $d=\dim_{\mathbb{R}}\mathbf{C},$%
\begin{equation}
\gamma:=\sup\left\{  \sum_{j=1}^{\infty}\left\vert \left\langle \omega\left(
h,e_{j}\right)  ,c\right\rangle _{\mathbf{C}}\right\vert ^{2}:\left\Vert
h\right\Vert _{H}=\left\Vert c\right\Vert _{\mathbf{C}}=1\right\}
\leqslant\left\Vert \omega\right\Vert _{2}^{2}<\infty\label{e.4.32}%
\end{equation}
and for $P\in\operatorname*{Proj}\left(  W\right)  $ let $B_{P}\left(
t\right)  :=PB\left(  t\right)  .$ Then for all
\begin{equation}
0\leqslant\lambda<\frac{\pi}{4dT\sqrt{\gamma}}, \label{e.4.33}%
\end{equation}%
\begin{equation}
\sup_{P\in\operatorname*{Proj}\left(  W\right)  }\mathbb{E}\left[  \exp\left(
\lambda\left\Vert \int_{0}^{t}\omega\left(  B_{P}\left(  \tau\right)
,dB_{P}\left(  \tau\right)  \right)  \right\Vert _{\mathbf{C}}\right)
\right]  <\infty\label{e.4.34}%
\end{equation}
and%
\begin{equation}
\mathbb{E}\left[  \exp\left(  \lambda\left\Vert \int_{0}^{t}\omega\left(
B\left(  \tau\right)  ,dB\left(  \tau\right)  \right)  \right\Vert
_{\mathbf{C}}\right)  \right]  <\infty. \label{e.4.35}%
\end{equation}

\end{prop}

\begin{proof}
Equation \eqref{e.4.35} follows by choosing $\left\{  P_{n}\right\}
_{n=1}^{\infty}\subset\operatorname*{Proj}\left(  W\right)  $ as in Eq.
\eqref{e.4.3} and then using Fatou's lemma in conjunction with the estimate in
Equation \eqref{e.4.34}. So we need only to concentrate on proving Equation \eqref{e.4.34}.

Fix a $P\in\operatorname*{Proj}\left(  W\right)  $ as in Equation
\eqref{e.3.42} and let%
\[
M_{t}^{P}:=\int_{0}^{t}\omega\left(  B_{P}\left(  \tau\right)  ,dB_{P}\left(
\tau\right)  \right)  .
\]
If $\left\{  f_{\ell}\right\}  _{\ell=1}^{d}$ is an orthonormal basis for
$\mathbf{C},$ then%
\[
\left\Vert M_{t}^{P}\right\Vert _{\mathbf{C}}\leqslant\sum_{\ell=1}%
^{d}\left\vert \left\langle M_{t}^{P},f_{\ell}\right\rangle _{\mathbf{C}%
}\right\vert ,
\]
and it follows by H\"{o}lder's inequality and the martingale estimate in
Proposition \ref{p.4.12} that%
\begin{align}
\mathbb{E}\left[  e^{\lambda\left\Vert M_{t}^{P}\right\Vert _{\mathbf{C}}%
}\right]  \leqslant &  \mathbb{E}\left[  e^{\lambda\sum_{\ell=1}^{d}\left\vert
\left\langle M_{t}^{P},f_{\ell}\right\rangle _{\mathbf{C}}\right\vert
}\right]  \leqslant\prod_{\ell=1}^{d}\left(  \mathbb{E}\left[  e^{\lambda
d\left\vert \left\langle M_{t}^{P},f_{\ell}\right\rangle _{\mathbf{C}%
}\right\vert }\right]  \right)  ^{1/d}\nonumber\\
\leqslant &  \prod_{\ell=1}^{d}\left(  2\sqrt{\mathbb{E}\left[  e^{2\lambda
^{2}d^{2}\left\langle \left\langle M_{\cdot}^{P},f_{\ell}\right\rangle
_{\mathbf{C}}\right\rangle _{t}}\right]  }\right)  ^{1/d}\nonumber\\
&  =2\prod_{\ell=1}^{d}\left(  \mathbb{E}\left[  e^{2\lambda^{2}%
d^{2}\left\langle \left\langle M_{\cdot}^{P},f_{\ell}\right\rangle
_{\mathbf{C}}\right\rangle _{t}}\right]  \right)  ^{1/2d}.\label{e.4.36}%
\end{align}

We will now evaluate each term in the product in Eq. (\ref{e.4.36}). So let
$c:=f_{\ell}$ and $N_{t}:=\left\langle M_{t}^{P},c\right\rangle _{\mathbf{C}%
},$ and $Q_{P}:H\rightarrow H$ and $Q:H\rightarrow H$ be the unique
non-negative symmetric operators such that, for all $h\in H,$%
\[
\sum_{j=1}^{n}\left\vert \left\langle \omega\left(  Ph,e_{j}\right)
,c\right\rangle _{\mathbf{C}}\right\vert ^{2}=\left\langle Q_{P}%
h,h\right\rangle _{H}\text{~}\text{ for all }h\in H
\]
and
\[
\sum_{j=1}^{\infty}\left\vert \left\langle \omega\left(  h,e_{j}\right)
,c\right\rangle _{\mathbf{C}}\right\vert ^{2}=\left\langle Qh,h\right\rangle
_{H}\text{ ~}\text{ for all }h\in H.
\]
Also let $\left\{  q_{l}\left(  P\right)  \right\}  _{l=1}^{\infty}$ be the
eigenvalues listed in decreasing order (counted with multiplicities) for
$Q_{P}$ and observe that
\begin{equation}
q_{1}\left(  P\right)  =\sup_{h\neq0}\frac{\left\langle Q_{P}h,h\right\rangle
}{\left\Vert h\right\Vert _{H}^{2}}\leqslant\sup_{h\neq0}\frac{\left\langle
QPh,Ph\right\rangle }{\left\Vert h\right\Vert _{H}^{2}}\leqslant\sup_{h\neq
0}\frac{\left\langle Qh,h\right\rangle }{\left\Vert h\right\Vert _{H}^{2}%
}\leqslant\gamma. \label{e.4.37}%
\end{equation}
With this notation, the quadratic variation of $N$ is given by%
\begin{equation}
\left\langle N\right\rangle _{T}=\int_{0}^{T}\sum_{j=1}^{n}\left\vert
\left\langle \omega\left(  B_{P}\left(  t\right)  ,e_{j}\right)
,c\right\rangle _{\mathbf{C}}\right\vert ^{2}dt=\int_{0}^{T}\left\langle
Q_{P}B_{P}\left(  t\right)  ,B_{P}\left(  t\right)  \right\rangle _{H}dt.
\label{e.4.38}%
\end{equation}
Moreover, by expanding $B_{P}\left(  \tau\right)  $ in an orthonormal basis of
eigenvectors of $Q_{P}|_{PH}$ it follows that
\begin{equation}
\left\langle N\right\rangle _{T}=\sum_{l=1}^{n}q_{l}\left(  P\right)  \int
_{0}^{T}b_{l}^{2}\left(  \tau\right)  d\tau\label{e.4.39}%
\end{equation}
where $\left\{  b_{l}\right\}  _{l=1}^{n}$ is a sequence of independent
Brownian motions. Hence it follows that
\begin{align}
\mathbb{E}\left[  e^{2\lambda^{2}d^{2}\left\langle \left\langle M_{\cdot}%
^{P},f_{\ell}\right\rangle _{\mathbf{C}}\right\rangle _{T}}\right]   &
=\mathbb{E}\left[  e^{2\lambda^{2}d^{2}\left\langle N\right\rangle _{T}%
}\right] \nonumber\\
&  =\prod_{l=1}^{n}\mathbb{E}\left[  \exp\left(  2\lambda^{2}d^{2}q_{l}\left(
P\right)  \int_{0}^{T}b_{l}^{2}\left(  \tau\right)  d\tau\right)  \right]  .
\label{e.4.40}%
\end{align}

If Eq. (\ref{e.4.33}) holds then (using Eq. (\ref{e.4.37}))%
\[
2\lambda d\sqrt{q_{1}\left(  P\right)  }=\sqrt{4\lambda^{2}d^{2}q_{1}\left(
P\right)  }\leqslant2\lambda d\sqrt{\gamma}<\pi/2T
\]
and we may apply Lemma \ref{l.4.10} to find
\begin{align}
\mathbb{E}\left[  \exp\left(  2\lambda^{2}d^{2}q_{l}\left(  P\right)  \int
_{0}^{T}b_{l}^{2}\left(  \tau\right)  d\tau\right)  \right]   &  =\frac
{1}{\sqrt{\cos\left(  2\lambda d\sqrt{q_{l}\left(  P\right)  }T\right)  }%
}\nonumber\\
&  =\exp\left(  -\frac{1}{2}\ln\cos\left(  2\lambda d\sqrt{q_{l}\left(
P\right)  }T\right)  \right)  . \label{e.4.41}%
\end{align}
Moreover, a simple calculus exercise shows for any $k\in\left(  0,\pi
/2\right)  $ there exists $c\left(  k\right)  <\infty$ such that $-\frac{1}%
{2}\ln\cos\left(  x\right)  \leqslant c\left(  k\right)  x^{2}$ for
$0\leqslant x\leqslant k.$ Taking $k=2\lambda d\sqrt{\gamma}T$ we may apply
this estimate to Eq. (\ref{e.4.41}) and combine the result with Eq.
(\ref{e.4.40}) to find
\[
\mathbb{E}\left[  e^{2\lambda^{2}d^{2}\left\langle \left\langle M_{\cdot}%
^{P},f_{\ell}\right\rangle _{\mathbf{C}}\right\rangle _{T}}\right]
\leqslant\prod_{l=1}^{n}\exp\left(  c\left(  k\right)  4\lambda^{2}d^{2}%
T^{2}q_{l}\left(  P\right)  \right)  =\exp\left(  c\left(  k\right)
4\lambda^{2}d^{2}T^{2}\operatorname*{tr}\left(  Q_{P}\right)  \right)  .
\]
Since $Q_{P}\leqslant PQP,$ we have
\[
\operatorname*{tr}Q_{P}\leqslant\operatorname{tr}Q=\sum_{l=1}^{\infty
}\left\langle Qe_{l},e_{l}\right\rangle _{H}=\sum_{j,l=1}^{\infty}\left\vert
\left\langle \omega\left(  e_{l},e_{j}\right)  ,c\right\rangle _{\mathbf{C}%
}\right\vert ^{2}=\left\Vert \left\langle \omega\left(  \cdot,\cdot\right)
,c\right\rangle _{\mathbf{C}}\right\Vert _{2}^{2}<\infty.
\]
Combining the last two equations (recalling that $c=f_{\ell})$ then shows,
\begin{equation}
\mathbb{E}\left[  e^{2\lambda^{2}d^{2}\left\langle \left\langle M_{\cdot}%
^{P},f_{\ell}\right\rangle _{\mathbf{C}}\right\rangle _{T}}\right]
\leqslant\exp\left(  c\left(  k\right)  4\lambda^{2}d^{2}T^{2}\left\Vert
\left\langle \omega\left(  \cdot,\cdot\right)  ,f_{\ell}\right\rangle
_{\mathbf{C}}\right\Vert _{2}^{2}\right)  . \label{e.4.42}%
\end{equation}
Using this estimate back in Eq. (\ref{e.4.36}) gives,%
\begin{align}
&  \mathbb{E}\left[  e^{\lambda\left\Vert M_{t}^{P}\right\Vert _{\mathbf{C}}%
}\right]  \leqslant2\exp\left(  c\left(  k\right)  2\lambda^{2}dT^{2}%
\sum_{\ell=1}^{d}\left\Vert \left\langle \omega\left(  \cdot,\cdot\right)
,f_{\ell}\right\rangle _{\mathbf{C}}\right\Vert _{2}^{2}\right) \label{e.4.43}%
\\
&  =2\exp\left(  c\left(  k\right)  2\lambda^{2}dT^{2}\left\Vert
\omega\right\Vert _{2}^{2}\right) \nonumber
\end{align}
which completes the proof as this last estimate is independent of
$P\in\operatorname*{Proj}\left(  W\right)  .$
\end{proof}

\begin{prop}
\label{p.4.14}Suppose that $\nu$ and $\mu$ are Gaussian measures on $W$ such
$q_{\nu}\left(  f\right)  :=\nu\left(  f^{2}\right)  \leqslant q_{\mu}\left(
f\right)  :=\mu\left(  f^{2}\right)  $ for all $f\in W_{\mathbb{R}}^{\ast}.$
If $g:[0,\infty)\rightarrow\lbrack0,\infty)$ is a non-negative,
non-decreasing, $C^{1}$ -- function, then
\[
\int_{W}g\left(  \left\Vert w\right\Vert \right)  d\nu\left(  w\right)
\leqslant\int_{W}g\left(  \left\Vert w\right\Vert \right)  d\mu\left(
w\right)  .
\]

\end{prop}

\begin{proof}
Theorem 3.3.6 in \cite[p. 107]{Bog98} states that if $q_{\nu}\leqslant q_{\mu
}$ then $\mu\left(  A\right)  \leqslant\nu\left(  A\right)  $ for every Borel
set $A$ which is convex and balanced. In particular, since $B_{t}:=\left\{
w\in W:\left\Vert w\right\Vert <t\right\}  $ is convex and balanced, it
follows that $\mu\left(  B_{t}\right)  \leqslant\nu\left(  B_{t}\right)  $ or
equivalently that $1-\nu\left(  B_{t}\right)  \leqslant1-\mu\left(
B_{t}\right)  $ for all $t\geqslant0.$ Since%
\begin{align}
\int_{W}g\left(  \left\Vert w\right\Vert \right)  d\nu\left(  w\right)   &
=\int_{W}\left[  g\left(  0\right)  +\int_{0}^{\infty}1_{t\leqslant\left\Vert
w\right\Vert }g^{\prime}\left(  t\right)  dt\right]  d\nu\left(  w\right)
\nonumber\\
&  =g\left(  0\right)  +\int_{0}^{\infty}\left(  g^{\prime}\left(  t\right)
\int_{W}1_{t\leqslant\left\Vert w\right\Vert }d\nu\left(  w\right)  \right)
dt\nonumber\\
&  =g\left(  0\right)  +\int_{0}^{\infty}g^{\prime}\left(  t\right)  \left[
1-\nu\left(  B_{t}\right)  \right]  dt \label{e.4.44}%
\end{align}
with the same formula holding when $\nu$ is replaced by $\mu,$ it follows that%
\begin{align*}
\int_{W}g\left(  \left\Vert w\right\Vert \right)  d\nu\left(  w\right)   &
=g\left(  0\right)  +\int_{0}^{\infty}g^{\prime}\left(  t\right)  \left[
1-\nu\left(  B_{t}\right)  \right]  dt\\
&  \leqslant g\left(  0\right)  +\int_{0}^{\infty}dtg^{\prime}\left(
t\right)  \left[  1-\mu\left(  B_{t}\right)  \right]  =\int_{W}g\left(
\left\Vert w\right\Vert \right)  d\mu\left(  w\right)  .
\end{align*}

\end{proof}

\begin{df}
\label{d.4.15}Let $\rho^{2}:G\rightarrow\lbrack0,\infty)$ be defined as
\[
\rho^{2}\left(  w,c\right)  :=\left\Vert w\right\Vert _{W}^{2}+\left\Vert
c\right\Vert _{\mathbf{C}}.
\]

\end{df}

In analogy to Gross' theory of measurable semi-norms (see e.g. Definition 5 in
\cite{Gross1963}) in the abstract Wiener space setting and in light of Theorem
\ref{t.3.12}, we view $\rho$ as a \textquotedblleft
measurable\textquotedblright\ extension of $d_{G_{CM}}$

\begin{thm}
[Integrated Gaussian heat kernel bounds]\label{t.4.16}There exists a
$\delta>0$ such that for all $\varepsilon\in\left(  0,\delta\right)  ,$ $T>0,$
$p\in\lbrack1,\infty),$
\begin{equation}
\sup_{P\in\operatorname*{Proj}\left(  W\right)  }\mathbb{E}\left[
e^{\frac{\varepsilon}{T}\rho^{2}\left(  g_{P}\left(  T\right)  \right)
}\right]  <\infty\text{ and }\int_{G}e^{\frac{\varepsilon}{T}\rho^{2}\left(
g\right)  }d\nu_{T}\left(  g\right)  <\infty\label{e.4.45}%
\end{equation}
whenever $\varepsilon<\delta.$
\end{thm}

\begin{proof}
Let $\varepsilon^{\prime}:=\varepsilon/T.$ For $P\in\operatorname*{Proj}%
\left(  W\right)  ,$
\[
\rho^{2}\left(  g_{P}\left(  T\right)  \right)  \leqslant\left\Vert
B_{P}\left(  T\right)  \right\Vert _{W}^{2}+\left\Vert B_{0}\left(  T\right)
\right\Vert _{\mathbf{C}}+\frac{1}{2}\left\Vert N_{P}\left(  T\right)
\right\Vert _{\mathbf{C}},
\]
where $N_{P}\left(  T\right)  :=\int_{0}^{T}\omega\left(  B_{P}\left(
t\right)  ,dB_{P}\left(  t\right)  \right)  $ and therefore,%
\[
\mathbb{E}\left[  e^{\varepsilon^{\prime}\rho^{2}\left(  g_{P}\left(
T\right)  \right)  }\right]  \leqslant\mathbb{E}\left[  e^{\varepsilon
^{\prime}\left[  \left\Vert B_{P}\left(  T\right)  \right\Vert _{W}^{2}%
+\frac{1}{2}\left\Vert N_{P}\left(  T\right)  \right\Vert _{\mathbf{C}%
}\right]  }\right]  \cdot\mathbb{E}\left[  e^{\varepsilon^{\prime}\left\Vert
B_{0}\left(  T\right)  \right\Vert _{\mathbf{C}}}\right]  .
\]
Moreover, by H\"{o}lder's inequality we have,%
\begin{align*}
\mathbb{E}\left[  e^{\varepsilon^{\prime}\rho^{2}\left(  g_{P}\left(
T\right)  \right)  }\right]   &  \leqslant\mathbb{E}\left[  e^{\varepsilon
^{\prime}\left\Vert B_{0}\left(  T\right)  \right\Vert _{\mathbf{C}}}\right]
\sqrt{\mathbb{E}\left[  e^{2\varepsilon^{\prime}\left\Vert B_{P}\left(
T\right)  \right\Vert _{W}^{2}}\right]  \cdot\mathbb{E}\left[  e^{\varepsilon
^{\prime}\left\Vert N_{P}\left(  T\right)  \right\Vert _{\mathbf{C}}}\right]
}\\
&  \leqslant\mathbb{E}\left[  e^{\varepsilon^{\prime}\left\Vert B_{0}\left(
T\right)  \right\Vert _{\mathbf{C}}}\right]  \sqrt{\mathbb{E}\left[
e^{2\varepsilon^{\prime}\left\Vert B\left(  T\right)  \right\Vert _{W}^{2}%
}\right]  \cdot\sup_{P^{\prime}\in\operatorname*{Proj}\left(  W\right)
}\mathbb{E}\left[  e^{\varepsilon^{\prime}\left\Vert N_{P^{\prime}}\left(
T\right)  \right\Vert _{\mathbf{C}}}\right]  }.
\end{align*}
wherein we have made use of Proposition \ref{p.4.14} to conclude that
\[
\mathbb{E}\left[  e^{2\varepsilon^{\prime}\left\Vert B_{P}\left(  T\right)
\right\Vert _{W}^{2}}\right]  \leqslant\mathbb{E}\left[  e^{2\varepsilon
^{\prime}\left\Vert B\left(  T\right)  \right\Vert _{W}^{2}}\right]
=\mathbb{E}\left[  e^{2\varepsilon^{\prime}T\left\Vert B\left(  1\right)
\right\Vert _{W}^{2}}\right]
\]
which is finite by Fernique's theorem provided $2\varepsilon=2\varepsilon
^{\prime}T<\delta^{\prime}$ for some $\delta^{\prime}>0.$ Similarly by
Proposition \ref{p.4.13},%
\[
\sup_{P^{\prime}\in\operatorname*{Proj}\left(  W\right)  }\mathbb{E}\left[
e^{\varepsilon^{\prime}\left\Vert N_{P^{\prime}}\left(  T\right)  \right\Vert
_{\mathbf{C}}}\right]  <\infty
\]
provided $\varepsilon=\varepsilon^{\prime}T<\frac{\pi}{4\sqrt{\gamma}}.$ The
assertion in Equation \eqref{e.4.45} now follows from these observations and
the fact that $\mathbb{E}\left[  e^{\varepsilon^{\prime}\left\Vert
B_{0}\left(  T\right)  \right\Vert _{\mathbf{C}}}\right]  <\infty$ for all
$\varepsilon^{\prime}>0.$
\end{proof}

\section{Path space quasi-invariance\label{s.5}}

\begin{nota}
\label{n.5.1}Let $W_{T}\left(  G\right)  $ denote the collection of continuous
paths, $g:\left[  0,T\right]  \rightarrow G$ such that $g\left(  0\right)
=\mathbf{e}.$ Moreover, if $V$ is a separable Hilbert space, let
$\mathcal{H}_{T}\left(  V\right)  $ denote the collection of absolutely
continuous functions (see \cite[pages 106-107]{Diestel-Uhl77}), $h:\left[
0,T\right]  \rightarrow V$ such that $h\left(  0\right)  =0$ and
\[
\left\Vert h\right\Vert _{\mathcal{H}_{T}\left(  V\right)  }:=\left(  \int
_{0}^{T}\left\Vert \dot{h}\left(  t\right)  \right\Vert _{V}^{2}dt\right)
^{1/2}<\infty.
\]
By polarization, we endow $\mathcal{H}_{T}\left(  V\right)  $ with the inner
product%
\[
\left\langle h,k\right\rangle _{\mathcal{H}_{T}\left(  V\right)  }=\int
_{0}^{T}\left\langle \dot{h}\left(  t\right)  , \dot{k}\left(  t\right)
\right\rangle _{V}dt.
\]

\end{nota}

\begin{thm}
[Path space quasi-invariance]\label{t.5.2}Suppose $T>0,$ $k\left(
\cdot\right)  =\left(  A\left(  \cdot\right)  ,a\left(  \cdot\right)  \right)
\in\mathcal{H}_{T}\left(  \mathfrak{g}_{CM}\right)  $ (thought of as a finite
energy path in $G_{CM}),$ and $g\left(  \cdot\right)  $ is the $G$--valued
Brownian motion in Equation \eqref{e.4.10}. Then over the finite time
interval, $\left[  0,T\right]  ,$ the laws of $k\cdot g$ and $g$ are
equivalent, i.e. they are mutually absolutely continuous relative to one
another. More precisely, if $F:W_{T}\left(  G\right)  \rightarrow\left[
0,\infty\right]  $ is a measurable function, then%
\begin{equation}
\mathbb{E}\left[  F\left(  k\cdot g\right)  \right]  =\mathbb{E}\left[
\tilde{Z}_{k}\left(  B,B_{0}\right)  F\left(  g\right)  \right]  ,
\label{e.5.1}%
\end{equation}
where
\begin{equation}
\tilde{Z}_{k}\left(  B,B_{0}\right)  :=\exp\left(
\begin{array}
[c]{c}%
\int_{0}^{T}\left\langle \dot{A}\left(  t\right)  ,dB\left(  t\right)
\right\rangle _{H}-\frac{1}{2}\int_{0}^{T}\left\Vert \dot{A}\left(  t\right)
\right\Vert _{H}^{2}dt\\
+\int_{0}^{T}\left\langle \dot{a}\left(  t\right)  +\frac{1}{2}\omega\left(
A\left(  t\right)  -2B\left(  t\right)  ,\dot{A}\left(  t\right)  \right)
,dB_{0}\left(  t\right)  \right\rangle _{\mathbf{C}}\\
-\frac{1}{2}\int_{0}^{T}\left\Vert \dot{a}\left(  t\right)  +\frac{1}{2}%
\omega\left(  A\left(  t\right)  -2B\left(  t\right)  ,\dot{A}\left(
t\right)  \right)  \right\Vert _{\mathbf{C}}^{2}dt
\end{array}
\right)  . \label{e.5.2}%
\end{equation}
Moreover, Equation \eqref{e.5.2} is valid for all measurable functions,
$F:W_{T}\left(  G\right)  \rightarrow\mathbb{C}$ such that
\[
\mathbb{E}\left[  \left\vert F\left(  k\cdot g\right)  \right\vert \right]
=\mathbb{E}\left[  \tilde{Z}_{k}\left(  B,B_{0}\right)  \left\vert F\left(
g\right)  \right\vert \right]  <\infty.
\]

\end{thm}

\begin{proof}
The Cameron--Martin theorem states (see for example, \cite[Theorem 1.2 on p.
113]{Kuo75}) that%
\begin{equation}
\mathbb{E}\left[  F\left(  B,B_{0}\right)  \right]  =\mathbb{E}\left[
Z_{k}\left(  B,B_{0}\right)  F\left(  \left(  B,B_{0}\right)  -k\right)
\right]  , \label{e.5.3}%
\end{equation}
where%
\begin{equation}
Z_{k}\left(  B,B_{0}\right)  :=\exp\left(
\begin{array}
[c]{c}%
\int_{0}^{T}\left[  \left\langle \dot{A}\left(  t\right)  , dB\left(
t\right)  \right\rangle _{H}+\left\langle \dot{a}\left(  t\right)
,dB_{0}\left(  t\right)  \right\rangle _{\mathbf{C}}\right] \\
-\frac{1}{2}\int_{0}^{T}\left[  \left\Vert \dot{A}\left(  t\right)
\right\Vert _{H}^{2}+\left\Vert \dot{a}\left(  t\right)  \right\Vert
_{\mathbf{C}}^{2}\right]  dt
\end{array}
\right)  . \label{e.5.4}%
\end{equation}
Since%
\begin{align}
&  \left(  k\cdot g\right)  \left(  t\right)  =\nonumber\\
&  \left(  B\left(  t\right)  +A\left(  t\right)  , B_{0}\left(  t\right)
+a\left(  t\right)  +\frac{1}{2}\int_{0}^{t}\omega\left(  B\left(
\tau\right)  , dB\left(  \tau\right)  \right)  +\frac{1}{2}\omega\left(
A\left(  t\right)  ,B\left(  t\right)  \right)  \right)  \label{e.5.5}%
\end{align}
is mapped to
\[
\left(  B\left(  t\right)  , B_{0}\left(  t\right)  +\frac{1}{2}\int_{0}%
^{t}\omega\left(  \left(  B-A\right)  \left(  \tau\right)  , d\left(
B-A\right)  \left(  \tau\right)  \right)  +\frac{1}{2}\omega\left(  A\left(
t\right)  , \left(  B-A\right)  \left(  t\right)  \right)  \right)
\]
under the transformation $B\rightarrow B-A$ and $B_{0}\rightarrow B_{0}-a,$ we
may conclude from Equation \eqref{e.5.3} that%
\begin{equation}
\mathbb{E}\left[  F\left(  k\cdot g\right)  \right]  =\mathbb{E}\left[
Z_{k}\left(  B,B_{0}\right)  F\left(  B,B_{0}+c\right)  \right]  ,
\label{e.5.6}%
\end{equation}
where%
\[
c\left(  t\right)  =\frac{1}{2}\int_{0}^{t}\omega\left(  \left(  B-A\right)
\left(  \tau\right)  , d\left(  B-A\right)  \left(  \tau\right)  \right)
+\frac{1}{2}\omega\left(  A\left(  t\right)  ,\left(  B-A\right)  \left(
t\right)  \right)  .
\]
By taking the differential of $c,$ one easily shows that
\[
c\left(  t\right)  =\frac{1}{2}\int_{0}^{t}\omega\left(  B\left(  \tau\right)
,dB\left(  \tau\right)  \right)  +u_{B}\left(  t\right)  ,
\]
where
\begin{equation}
u_{B}\left(  t\right)  :=\frac{1}{2}\int_{0}^{t}\omega\left(  A\left(
\tau\right)  -2B\left(  \tau\right)  , \dot{A}\left(  \tau\right)  \right)
d\tau. \label{e.5.7}%
\end{equation}
Hence Equation \eqref{e.5.6} may be rewritten as
\begin{equation}
\mathbb{E}\left[  F\left(  k\cdot g\right)  \right]  =\mathbb{E}\left[
Z_{k}\left(  B,B_{0}\right)  F\left(  B,B_{0}+u_{B}+\frac{1}{2}\int_{0}%
^{\cdot}\omega\left(  B\left(  t\right)  , dB\left(  t\right)  \right)
\right)  \right]  . \label{e.5.8}%
\end{equation}
Freezing the integration over $B$ (i.e. using Fubini's theorem) we may use the
Cameron-Martin theorem one more time to make the transformation,
$B_{0}\rightarrow B_{0}-u_{B}.$ Doing so gives
\begin{align}
\mathbb{E}\left[  F\left(  k\cdot g\right)  \right]   &  =\mathbb{E}\left[
\tilde{Z}_{k}\left(  B,B_{0}\right)  F\left(  \left(  B,B_{0}+\frac{1}{2}%
\int_{0}^{\cdot}\omega\left(  B\left(  t\right)  , dB\left(  t\right)
\right)  \right)  \right)  \right] \nonumber\\
&  =\mathbb{E}\left[  \tilde{Z}_{k}\left(  B,B_{0}\right)  F\left(  g\right)
\right]  , \label{e.5.9}%
\end{align}
where
\begin{equation}
\tilde{Z}_{k}\left(  B,B_{0}\right)  :=Z_{k}\left(  B,B_{0}-u_{B}\right)
\exp\left(  \int_{0}^{T}\left\langle \dot{u}_{B}\left(  t\right)
,dB_{0}\left(  t\right)  \right\rangle _{\mathbf{C}}-\frac{1}{2}\int_{0}%
^{T}\left\Vert \dot{u}_{B}\left(  t\right)  \right\Vert _{\mathbf{C}}%
^{2}dt\right)  . \label{e.5.10}%
\end{equation}
A little algebra shows that $\tilde{Z}_{k}\left(  B,B_{0}\right)  $ defined in
Equation \eqref{e.5.10} may be expressed as in Equation \eqref{e.5.2}.
\end{proof}

\begin{rem}
\label{r.5.3}The above proof fails if we try to use it to prove the right
quasi-invariance on the path space measure, i.e. that $g\cdot k$ has a law
which is absolutely continuous to that of $g.$ In this case
\[
\left(  g\cdot k\right)  \left(  t\right)  =\left(  B\left(  t\right)
+A\left(  t\right)  ,B_{0}\left(  t\right)  +a\left(  t\right)  +\frac{1}%
{2}\int_{0}^{t}\omega\left(  B\left(  \tau\right)  ,dB\left(  \tau\right)
\right)  -\frac{1}{2}\omega\left(  A\left(  t\right)  ,B\left(  t\right)
\right)  \right)
\]
and then making the transformation, $B\rightarrow B-A$ and $B_{0}\rightarrow
B_{0}-a$ gives%
\[
\mathbb{E}\left[  F\left(  g\cdot k\right)  \right]  =\mathbb{E}\left[
Z_{k}\left(  B,B_{0}\right)  F\left(  B,B_{0}+c\right)  \right]
\]
where%
\[
c\left(  t\right)  =\frac{1}{2}\int_{0}^{t}\omega\left(  B\left(  \tau\right)
,dB\left(  \tau\right)  \right)  +u_{B}\left(  t\right)
\]
and%
\[
u_{B}\left(  t\right)  =\frac{1}{2}\int_{0}^{t}\left[  \omega\left(
A,dA\right)  -2\omega\left(  A,dB\right)  \right]  .
\]
The argument breaks down at this point since $u_{B}$ is no longer absolutely
continuous in $t.$ Hence we can no longer use the Cameron -- Martin theorem to
translate away the $u_{B}$ term.
\end{rem}

\begin{prop}
\label{p.5.4}There exists a $\delta>0$ and a function $C\left(  p,u\right)
\in(0,\infty],$ for $1<p<\infty$ and $0\leqslant u<\infty,$ which is
non-decreasing in each of its variables, $C\left(  p,u\right)  <\infty$
whenever%
\begin{equation}
p\leqslant\frac{1}{2}\left(  1+\sqrt{1+\delta/u}\right)  , \label{e.5.11}%
\end{equation}
and,%
\begin{equation}
\mathbb{E}\left[  \tilde{Z}_{k}\left(  B,B_{0}\right)  ^{p}\right]  \leqslant
C\left(  p,\left\Vert k\right\Vert _{\mathcal{H}_{T}\left(  \mathfrak{g}%
_{CM}\right)  }\right)  \text{ for all }k\in\mathcal{H}_{T}\left(
\mathfrak{g}_{CM}\right)  . \label{e.5.12}%
\end{equation}

\end{prop}

\begin{proof}
For the purposes of this proof, let $\mathbb{E}_{B_{0}}$ and $\mathbb{E}_{B}$
denote the expectation relative to $B_{0}$ and $B$ respectively, so that by
Fubini's theorem $\mathbb{E=E}_{B_{0}}\mathbb{E}_{B}=\mathbb{E}_{B}%
\mathbb{E}_{B_{0}},$ We may write $\tilde{Z}_{k}\left(  B,B_{0}\right)  $ as
\[
\tilde{Z}_{k}\left(  B,B_{0}\right)  :=\zeta\left(  B\right)  \exp\left(
\int_{0}^{T}\left\langle \dot{a}\left(  t\right)  +\dot{u}_{B}\left(
t\right)  , dB_{0}\left(  t\right)  \right\rangle _{\mathbf{C}}\right)
\]
where
\[
\zeta\left(  B\right)  :=\exp\left(  \int_{0}^{T}\left\langle \dot{A}\left(
t\right)  , dB\left(  t\right)  \right\rangle _{H}-\frac{1}{2}\int_{0}%
^{T}\left\Vert \dot{A}\left(  t\right)  \right\Vert _{H}^{2}dt-\frac{1}{2}%
\int_{0}^{T}\left\Vert \dot{a}\left(  t\right)  +\dot{u}_{B}\left(  t\right)
\right\Vert _{\mathbf{C}}^{2}dt\right)
\]
and $u_{B}\left(  t\right)  $ is as in Equation \eqref{e.5.7}. Hence it
follows that,%
\begin{align}
\mathbb{E}_{B_{0}}\left[  \tilde{Z}_{k}\left(  B,B_{0}\right)  ^{p}\right]
&  =\zeta^{p}\left(  B\right)  \mathbb{E}_{B_{0}}\left[  \exp\left(  p\int
_{0}^{T}\left\langle \dot{a}\left(  t\right)  +\dot{u}_{B}\left(  t\right)
,dB_{0}\left(  t\right)  \right\rangle _{\mathbf{C}}\right)  \right]
\nonumber\\
&  =\zeta^{p}\left(  B\right)  \exp\left(  \frac{p^{2}}{2}\int_{0}%
^{T}\left\Vert \dot{a}\left(  t\right)  +\dot{u}_{B}\left(  t\right)
\right\Vert _{C}^{2}dt\right)  =U V \label{e.5.13}%
\end{align}
where%
\[
U:=\exp\left(  p\left(  \int_{0}^{T}\left\langle \dot{A}\left(  t\right)
,dB\left(  t\right)  \right\rangle _{H}-\frac{1}{2}\int_{0}^{T}\left\Vert
\dot{A}\left(  t\right)  \right\Vert _{H}^{2}dt\right)  \right)
\]
and
\[
V=\exp\left(  \frac{p^{2}-p}{2}\int_{0}^{T}\left\Vert \dot{a}\left(  t\right)
+\dot{u}_{B}\left(  t\right)  \right\Vert _{\mathbf{C}}^{2}dt\right)  .
\]
Note that when $p=1,$ Equation \eqref{e.5.13} becomes%
\[
\mathbb{E}_{B_{0}}\left[  \tilde{Z}_{k}\left(  B,B_{0}\right)  \right]
=\exp\left(  \int_{0}^{T}\left\langle \dot{A}\left(  t\right)  , dB\left(
t\right)  \right\rangle _{H}-\frac{1}{2}\int_{0}^{T}\left\Vert \dot{A}\left(
t\right)  \right\Vert _{H}^{2}dt\right)  ,
\]
from which it easily follows that
\[
\mathbb{E}\left[  \tilde{Z}_{k}\left(  B,B_{0}\right)  \right]  =\mathbb{E}%
_{B}\mathbb{E}_{B_{0}}\left[  \tilde{Z}_{k}\left(  B,B_{0}\right)  \right]
=1.
\]

Now suppose that $p>1.$ By the Cauchy -- Schwarz inequality,%
\[
\mathbb{E}\left[  \tilde{Z}_{k}\left(  B,B_{0}\right)  ^{p}\right]
=\mathbb{E}_{B}\left[  UV\right]  \leqslant\left(  \mathbb{E}_{B}\left[
V^{2}\right]  \right)  ^{1/2}\left(  \mathbb{E}_{B}\left[  U^{2}\right]
\right)  ^{1/2}.
\]
Because%
\begin{align*}
\mathbb{E}U^{2}  &  =\exp\left(  -p\int_{0}^{T}\left\Vert \dot{A}\left(
t\right)  \right\Vert _{H}^{2}dt\right)  \mathbb{E}\left[  \exp\left(
2p\int_{0}^{T}\left\langle \dot{A}\left(  t\right)  ,dB\left(  t\right)
\right\rangle _{H}\right)  \right] \\
&  =\exp\left(  p\left\Vert A\right\Vert _{\mathcal{H}_{T}\left(  H\right)
}^{2}\right)  \leqslant\exp\left(  p\left\Vert k\right\Vert _{\mathcal{H}%
_{T}\left(  \mathfrak{g}_{CM}\right)  }^{2}\right)  <\infty,
\end{align*}
we have reduced the problem to estimating $\mathbb{E}V^{2}.$ By elementary
estimates we have%
\begin{align*}
\left\Vert \dot{u}_{B}\left(  t\right)  \right\Vert _{\mathbf{C}}^{2}  &
=\frac{1}{4}\left\Vert \omega\left(  A\left(  t\right)  -2B\left(  t\right)
,\dot{A}\left(  t\right)  \right)  \right\Vert _{\mathbf{C}}^{2}\\
&  \leqslant\frac{1}{4}\left\Vert \omega\right\Vert _{0}^{2}\left\Vert
A\left(  t\right)  -2B\left(  t\right)  \right\Vert _{W}^{2}\left\Vert \dot
{A}\left(  t\right)  \right\Vert _{W}^{2}\\
&  \leqslant\frac{1}{2}\left\Vert \omega\right\Vert _{0}^{2}\left\Vert \dot
{A}\left(  t\right)  \right\Vert _{W}^{2}\left(  \left\Vert A\left(  t\right)
\right\Vert _{W}^{2}+4\left\Vert B\left(  t\right)  \right\Vert _{W}%
^{2}\right)
\end{align*}
and hence
\begin{align}
\left\Vert \dot{a}\left(  t\right)  +\dot{u}_{B}\left(  t\right)  \right\Vert
^{2}  &  \leqslant2\left\Vert \dot{a}\left(  t\right)  \right\Vert
_{\mathbf{C}}^{2}+2\left\Vert \dot{u}_{B}\left(  t\right)  \right\Vert
_{\mathbf{C}}^{2}\nonumber\\
&  =2\left\Vert \dot{a}\left(  t\right)  \right\Vert _{\mathbf{C}}%
^{2}+\left\Vert \omega\right\Vert _{0}^{2}\left\Vert \dot{A}\left(  t\right)
\right\Vert _{W}^{2}\left(  \left\Vert A\left(  t\right)  \right\Vert _{W}%
^{2}+4\sup_{0\leqslant t\leqslant T}\left\Vert B\left(  t\right)  \right\Vert
_{W}^{2}\right)  . \label{e.5.14}%
\end{align}
By Equation \eqref{e.2.10} there exits $c<\infty$ such that $\left\Vert
\cdot\right\Vert _{W}\leqslant c\left\Vert \cdot\right\Vert _{H}.$ Since%
\[
\left\Vert A\left(  t\right)  \right\Vert _{H}\leqslant\int_{0}^{T}\left\Vert
\dot{A}\left(  \tau\right)  \right\Vert _{H}d\tau\leqslant\sqrt{T}\left\Vert
A\right\Vert _{\mathcal{H}_{T}\left(  H\right)  },
\]
we find%
\[
V^{2}\leqslant C\exp\left(  4\left(  p^{2}-p\right)  c^{2}\left\Vert
\omega\right\Vert _{0}^{2}\left\Vert A\right\Vert _{\mathcal{H}_{T}\left(
H\right)  }^{2}\sup_{0\leqslant t\leqslant T}\left\Vert B\left(  t\right)
\right\Vert _{W}^{2}\right)
\]
where%
\begin{align*}
C  &  =\exp\left(  \left(  p^{2}-p\right)  \left(  2\left\Vert a\right\Vert
_{\mathcal{H}_{T}\left(  \mathbf{C}\right)  }^{2}+c^{4}T\left\Vert
\omega\right\Vert _{0}^{2}\left\Vert A\right\Vert _{\mathcal{H}_{T}\left(
H\right)  }^{4}\right)  \right) \\
&  \leqslant C^{\prime}\left(  p,\left\Vert k\right\Vert _{\mathcal{H}%
_{T}\left(  \mathfrak{g}_{CM}\right)  }\right)  <\infty.
\end{align*}
Now by Fernique's theorem as in Equation \ref{e.2.4} there exists
$\delta^{\prime}>0$ such that
\[
M:=\mathbb{E}\left[  \exp\left(  \delta^{\prime}\sup_{0\leqslant t\leqslant
T}\left\Vert B\left(  t\right)  \right\Vert _{W}^{2}\right)  \right]  <\infty
\]
and hence it follows that%
\[
\mathbb{E}V^{2}\leqslant C^{\prime}\left(  p,\left\Vert k\right\Vert
_{\mathcal{H}_{T}\left(  \mathfrak{g}_{CM}\right)  }\right)  \cdot M<\infty
\]
provided%
\[
4\left(  p^{2}-p\right)  c^{2}\left\Vert \omega\right\Vert _{0}^{2}\left\Vert
A\right\Vert _{\mathcal{H}_{T}\left(  H\right)  }^{2}\leqslant4\left(
p^{2}-p\right)  c^{2}\left\Vert \omega\right\Vert _{0}^{2}\left\Vert
k\right\Vert _{\mathcal{H}_{T}\left(  \mathfrak{g}_{CM}\right)  }^{2}%
\leqslant\delta^{\prime}.
\]
The latter condition holds provided%
\[
p\leqslant\frac{1+\sqrt{1+\delta/\left\Vert k\right\Vert _{\mathcal{H}%
_{T}\left(  \mathfrak{g}_{CM}\right)  }^{2}}}{2}%
\]
where $\delta:=\left(  c^{2}\left\Vert \omega\right\Vert _{0}^{2}\right)
^{-1}\delta^{\prime}>0.$
\end{proof}

\begin{df}
\label{d.5.5}We will say that a function, $F:W_{T}\left(  G\right)
\rightarrow\mathbb{R}$ ($W_{T}\left(  G\right)  $ as in Notation \ref{n.5.1})
is \textbf{polynomially bounded} if there exist constants $K,M<\infty$ such
that
\begin{equation}
\left\vert F\left(  g\right)  \right\vert \leqslant K\left(  1+\sup
_{t\in\left[  0,T\right]  }\left\Vert g\left(  t\right)  \right\Vert
_{G}\right)  ^{M}\text{ for all }g\in W_{T}\left(  G\right)  . \label{e.5.15}%
\end{equation}
Given a finite energy path, $k\left(  t\right)  =\left(  A\left(  t\right)
,a\left(  t\right)  \right)  \in\mathfrak{g}_{CM},$ we say that $F$ is right
$k$ differentiable if
\[
\frac{d}{ds}\Big|_{0}F\left(  \left(  sk\right)  \cdot g\right)  =:\left(
\hat{k}F\right)  \left(  g\right)
\]
exists for all $g\in W_{T}\left(  G\right)  .$
\end{df}

\begin{cor}
[Path space integration by parts]\label{c.5.6}Let $k\left(  \cdot\right)
=\left(  A\left(  \cdot\right)  ,a\left(  \cdot\right)  \right)
\in\mathcal{H}_{T}\left(  \mathfrak{g}_{CM}\right)  $ and $F:W_{T}\left(
G\right)  \rightarrow\mathbb{R}$ be a $k$ -- differentiable function such that
$F$ and $\hat{k}F$ are polynomial bounded functions on $W_{T}\left(  G\right)
.$ Then%
\begin{equation}
\mathbb{E}\left[  \left(  \hat{k}F\right)  \left(  g\right)  \right]
=\mathbb{E}\left[  F\left(  g\right)  z_{k}\right]  \label{e.5.16}%
\end{equation}
where%
\begin{equation}
z_{k}:=\int_{0}^{T}\left[  \left\langle \dot{A}\left(  t\right)  ,dB\left(
t\right)  \right\rangle _{H}+\left\langle \dot{a}\left(  t\right)
-\omega\left(  B\left(  t\right)  ,\dot{A}\left(  t\right)  \right)
,dB_{0}\left(  t\right)  \right\rangle _{\mathbf{C}}\right]  . \label{e.5.17}%
\end{equation}
Moreover, $\mathbb{E}\left\vert z_{k}\right\vert ^{p}<\infty$ for all
$p\in\lbrack1,\infty).$
\end{cor}

\begin{proof}
From Theorem \ref{t.5.2}, we have that for any $s\in\mathbb{R}$
\begin{equation}
\mathbb{E}\left[  F\left(  \left(  sk\right)  \cdot g\right)  \right]
=\mathbb{E}\left[  \tilde{Z}_{sk}\left(  B,B_{0}\right)  F\left(  g\right)
\right]  . \label{e.5.18}%
\end{equation}
Formally differentiating this identity at $s=0$ and interchanging the
derivatives with the expectations immediately leads to Equation
\eqref{e.5.16}. To make this rigorous we need only to verify that derivative
interchanges are permissible. From Equations \eqref{e.3.8} and \eqref{e.5.15},
there exists $C\left(  k\right)  <\infty$ such that%
\begin{align*}
\sup_{\left\vert s\right\vert \leqslant1}\left\vert \frac{d}{ds}F\left(
\left(  sk\right)  \cdot g\right)  \right\vert  &  =\sup_{\left\vert
s\right\vert \leqslant1}\left\vert \left(  \hat{k}F\right)  \left(  \left(
sk\right)  \cdot g\right)  \right\vert \\
&  \leqslant K\sup_{\left\vert s\right\vert \leqslant1}\left(  1+\sup
_{t\in\left[  0,T\right]  }\left\Vert \left[  sk\left(  t\right)  \right]
\cdot g\left(  t\right)  \right\Vert _{G}\right)  ^{M}\\
&  \leqslant C\left(  k\right)  \left(  1+\sup_{t\in\left[  0,T\right]
}\left\Vert g\left(  t\right)  \right\Vert _{G}\right)  ^{M}%
\end{align*}
wherein the last expression is integrable by Fernique's theorem and the moment
estimate in Proposition \ref{p.4.1}. Therefore,%
\[
\frac{d}{ds}|_{0}\mathbb{E}\left[  F\left(  \left(  sk\right)  \cdot g\right)
\right]  =\mathbb{E}\left[  \frac{d}{ds}|_{0}F\left(  \left(  sk\right)  \cdot
g\right)  \right]  =\mathbb{E}\left[  \left(  \hat{k}F\right)  \left(
g\right)  \right]  .
\]

To see that we may also differentiate the right side of Equation
\eqref{e.5.18}, observe that%
\[
\tilde{Z}_{sk}\left(  B,B_{0}\right)  =\exp\left(  sz_{k}+s^{2}\beta
+s^{3}\gamma+s^{4}\kappa\right)
\]
where%
\begin{align*}
\beta=  &  -\frac{1}{2}\int_{0}^{T}\left\Vert \dot{A}\left(  t\right)
\right\Vert _{H}^{2}dt+\frac{1}{2}\int_{0}^{T}\left\langle \omega\left(
A\left(  t\right)  ,\dot{A}\left(  t\right)  \right)  ,dB_{0}\left(  t\right)
\right\rangle _{\mathbf{C}}\\
&  \qquad-\frac{1}{2}\int_{0}^{T}\left\Vert \dot{a}\left(  t\right)
-\omega\left(  B\left(  t\right)  ,\dot{A}\left(  t\right)  \right)
\right\Vert _{\mathbf{C}}^{2}dt,\\
\gamma=  &  -\frac{1}{2}\int_{0}^{T}\operatorname{Re}\left\langle \dot
{a}\left(  t\right)  -\omega\left(  B\left(  t\right)  ,\dot{A}\left(
t\right)  \right)  ,\omega\left(  A\left(  t\right)  ,\dot{A}\left(  t\right)
\right)  \right\rangle _{\mathbf{C}}dt,
\end{align*}
and%
\[
\kappa=-\frac{1}{8}\int_{0}^{T}\left\Vert \omega\left(  A\left(  t\right)
,\dot{A}\left(  t\right)  \right)  \right\Vert _{\mathbf{C}}^{2}dt.
\]
Using Fernique's theorem again and estimates similar to those used in the
proof of Proposition \ref{p.5.4}, one shows for any $p\in\lbrack1,\infty)$
that there exists $s_{0}\left(  p\right)  >0$ such that
\[
\mathbb{E}\left[  \sup_{\left\vert s\right\vert \leqslant s_{0}\left(
p\right)  }\left\vert \frac{d}{ds}\tilde{Z}_{sk}\left(  B,B_{0}\right)
\right\vert ^{p}\right]  <\infty.
\]
Therefore we may differentiate past the expectation to find%
\[
\frac{d}{ds}|_{0}\mathbb{E}\left[  F\left(  g\right)  \tilde{Z}_{sk}\left(
B,B_{0}\right)  \right]  =\mathbb{E}\left[  F\left(  g\right)  \frac{d}%
{ds}|_{0}\tilde{Z}_{sk}\left(  B,B_{0}\right)  \right]  =\mathbb{E}\left[
F\left(  g\right)  z_{k}\right]  .
\]
The fact that $z_{k}$ has finite moments of all orders follows by the
martingale arguments along with Nelson's theorem as described in the proof of
Lemma \ref{l.4.7}. Alternatively, observe that $\int_{0}^{T}\left\langle
\dot{A}\left(  t\right)  ,dB\left(  t\right)  \right\rangle _{H}$ is Gaussian
and hence has finite moments of all orders. If we let $M_{t}:=\int_{0}%
^{t}\left\langle \dot{a}-\omega\left(  B,\dot{A}\right)  ,dB_{0}\right\rangle
_{\mathbf{C}},$ then $M$ is a martingale such that%
\[
\left\langle M\right\rangle _{T}=\int_{0}^{T}\left\Vert \dot{a}\left(
t\right)  -\omega\left(  B\left(  t\right)  ,\dot{A}\left(  t\right)  \right)
\right\Vert _{C}^{2}dt\leqslant C\left(  1+\max_{0\leqslant t\leqslant
T}\left\Vert B\left(  t\right)  \right\Vert _{W}^{2}\right)  .
\]
So by Fernique's theorem, $\mathbb{E}\left[  \left\langle M\right\rangle
_{T}^{p}\right]  <\infty$ for all $p<\infty$ and hence by the
Burkholder-Davis-Gundy inequalities, $\mathbb{E}\left\vert M_{T}\right\vert
^{p}<\infty$ for all $1\leqslant p<\infty.$
\end{proof}

\section{Heat Kernel Quasi-Invariance\label{s.6}}

In this section we will use the results of Section \ref{s.5} to prove both
quasi-invariance of the heat kernel measures, $\left\{  \nu_{T}\right\}
_{T>0},$ relative to left and right translations by elements of $G_{CM}.$

\begin{thm}
[Left quasi-invariance of the heat kernel measure]\label{t.6.1} Let $T>0$ and
$\left(  A,a\right)  \in G_{CM}.$ Then $\left(  A,a\right)  \cdot g\left(
T\right)  $ and $g\left(  T\right)  $ have equivalent laws. More precisely, if
$f:G\rightarrow\left[  0,\infty\right]  $ is a measurable function, then%
\begin{equation}
\mathbb{E}\left[  f\left(  \left(  A,a\right)  \cdot g\left(  T\right)
\right)  \right]  =\mathbb{E}\left[  f\left(  g\left(  T\right)  \right)
\bar{Z}_{k}\left(  g(T\right)  \right]  , \label{e.6.1}%
\end{equation}
where%
\begin{equation}
\bar{Z}_{k}\left(  g\left(  T\right)  \right)  =\mathbb{E}\left[  \left.
\zeta_{\left(  A,a\right)  }\left(  B,B_{0}\right)  \right\vert \sigma\left(
g\left(  T\right)  \right)  \right]  \label{e.6.2}%
\end{equation}
and%
\begin{align}
\ln\zeta_{\left(  A,a\right)  }\left(  B,B_{0}\right)  :=  &  \frac{1}%
{T}\left\langle A,B\left(  T\right)  \right\rangle _{H}-\frac{\left\Vert
A\right\Vert _{H}^{2}}{2T^{2}}+\frac{1}{T}\int_{0}^{T}\left\langle
a-\omega\left(  B\left(  t\right)  ,A\right)  ,dB_{0}\left(  t\right)
\right\rangle _{\mathbf{C}}\nonumber\\
&  \qquad-\frac{1}{2T^{2}}\int_{0}^{T}\left\Vert a-\omega\left(  B\left(
t\right)  ,A\right)  \right\Vert _{\mathbf{C}}^{2}dt. \label{e.6.3}%
\end{align}

\end{thm}

\begin{proof}
An application of Theorem \ref{t.5.2} with $F\left(  g\right)  :=f\left(
g\left(  T\right)  \right)  $ and $k\left(  t\right)  :=\frac{t}{T}\left(
A,a\right)  $ implies%
\begin{align}
\mathbb{E}\left[  f\left(  \left(  A,a\right)  \cdot g\left(  T\right)
\right)  \right]   &  =\mathbb{E}\left[  F\left(  \mathbf{k}\cdot g\right)
\right]  =\mathbb{E}\left[  \tilde{Z}_{k}\left(  B,B_{0}\right)  \cdot
F\left(  g\right)  \right] \nonumber\\
&  =\mathbb{E}\left[  \tilde{Z}_{k}\left(  B,B_{0}\right)  f\left(  g\left(
T\right)  \right)  \right]  , \label{e.6.4}%
\end{align}
where after a little manipulation one shows, $\tilde{Z}_{k}\left(
B,B_{0}\right)  =\zeta_{\left(  A,a\right)  }\left(  B,B_{0}\right)  .$ By
conditioning on $\sigma\left(  g\left(  T\right)  \right)  $ we can also write
Equation \eqref{e.6.4} as in Equation \eqref{e.6.1}.
\end{proof}

\begin{cor}
[Right quasi-invariance of the heat kernel measure]\label{c.6.2}The heat
kernel measure, $\nu_{T},$ is also quasi-invariant under right translations,
and
\begin{equation}
\frac{d\nu_{T}\circ r_{k}^{-1}}{d\nu_{T}}\left(  g\right)  =\bar{Z}_{k^{-1}%
}\left(  g^{-1}\right)  . \label{e.6.5}%
\end{equation}
where%
\[
\bar{Z}_{k}=d\nu_{T}\circ l_{k}^{-1}/d\nu_{T}%
\]
is as in Theorem \ref{t.6.1}.
\end{cor}

\begin{proof}
Recall from Corollary \ref{c.4.9} that $\nu_{T}$ is invariant under the
inversion map, $g\rightarrow g^{-1}.$ From this observation and Theorem
\ref{t.6.1} it follows that $\nu_{T}$ is also quasi-invariant under right
translations of elements of $G_{CM}.$ In more detail, if $k\in G_{CM}$ and
$f:G\rightarrow\mathbb{R}$ is a bounded measurable function, then
\begin{align*}
\int_{G}f\left(  g\cdot k\right)  d\nu_{T}\left(  g\right)   &  =\int
_{G}f\left(  g^{-1}\cdot k\right)  d\nu_{T}\left(  g\right)  =\int_{G}f\left(
\left(  k^{-1}g\right)  ^{-1}\right)  d\nu_{T}\left(  g\right) \\
&  =\int_{G}f\left(  g^{-1}\right)  \bar{Z}_{k^{-1}}\left(  g\right)  d\nu
_{T}\left(  g\right)  =\int_{G}f\left(  g\right)  \bar{Z}_{k^{-1}}\left(
g^{-1}\right)  d\nu_{T}\left(  g\right)  .
\end{align*}
Equation \eqref{e.6.5} is a consequence of this identity.
\end{proof}

Just like in the case abstract Wiener spaces we have the following strong
converses of Theorem \ref{t.6.1} and Corollary \ref{c.6.2}.

\begin{prop}
\label{p.6.3}Suppose that $k\in G\setminus G_{CM}$ and $T>0,$ then
$\nu_{T}\circ l_{k}^{-1}$ and $\nu_{T}$ are singular and $\nu_{T}\circ
r_{k}^{-1}$ and $\nu_{T}$ are singular.
\end{prop}

\begin{proof}
Let $k=\left(  A,a\right)  \in G\setminus G_{CM}$ with $a\in\mathbf{C}$ and
$A\in W\setminus H.$ Given a measurable subset, $V\subset W,$ we have
\[
\nu_{T}\left(  V\times\mathbf{C}\right)  =P\left(  B\left(  T\right)  \in
V\right)  =:\mu_{T}\left(  V\right)  ,
\]
where $\mu_{T}$ is Wiener measure on $W$ with variance $T.$ It is well known
(see e.g. Corollary 2.5.3 in \cite{Bog98}) that if $A\in W\setminus H$ that
$\mu_{T}\left(  \cdot-A\right)  $ is singular relative to $\mu_{T}\left(
\cdot\right)  ,$ i.e. we may partition $W$ into two disjoint measurable sets,
$W_{0}$ and $W_{1}$ such that $\mu_{T}\left(  W_{0}\right)  =1=\mu_{T}\left(
W_{1}-A\right)  .$ A simple computation shows for any $V\subset W$ that%
\[
l_{k}^{-1}\left(  V\times\mathbf{C}\right)  =r_{k}^{-1}\left(  V\times
\mathbf{C}\right)  =\left(  V-A\right)  \times\mathbf{C}.
\]
Thus if we define $G_{i}:=W_{i}\times\mathbf{C}$ for $i=0,1,$ we have that $G$
is the disjoint union of $G_{0}$ and $G_{1}$ and $\nu_{T}\left(  G_{0}\right)
=\mu_{T}\left(  W_{0}\right)  =1$ while
\[
\nu_{T}\left(  r_{k}^{-1}\left(  G_{1}\right)  \right)  =\nu_{T}\left(
l_{k}^{-1}\left(  G_{1}\right)  \right)  =\nu_{T}\left(  \left(
W_{1}-A\right)  \times\mathbf{C}\right)  =\mu_{T}\left(  W_{1}-A\right)  =1.
\]

\end{proof}

\begin{cor}
[Right heat kernel integration by parts]\label{c.6.4}Let $k:=\left(
A,a\right)  \in\mathfrak{g}_{CM}$ and suppose that $f:G\rightarrow\mathbb{C}$
is a smooth function such that $f$ and $\hat{k}f$ are polynomially bounded.
Then%
\[
\mathbb{E}\left[  \left(  \hat{k}f\right)  \left(  g\left(  T\right)  \right)
\right]  =\mathbb{E}\left[  f\left(  g\left(  T\right)  \right)  z_{k}\right]
\]
where $\hat{k}f\left(  g\right)  :=\frac{d}{ds}\Big|_{0}f\left(  \left(
sk\right)  g\right)  $ and
\[
z_{k}:=\frac{1}{T}\left[  \left\langle A,B\left(  T\right)  \right\rangle
_{H}+\left\langle a,B_{0}\left(  T\right)  \right\rangle _{\mathbf{C}}%
-\int_{0}^{T}\left\langle \omega\left(  B\left(  t\right)  ,A\right)
,dB_{0}\left(  t\right)  \right\rangle _{\mathbf{C}}\right]  .
\]
Moreover, with $\nu_{T}:=\operatorname{Law}\left(  g\left(  T\right)  \right)
,$ the above formula gives,%
\[
\int_{G}\left(  \hat{k}f\right)  d\nu_{T}\left(  g\right)  =\int_{G}f\left(
g\right)  \bar{z}_{k}\left(  g\right)  d\nu_{T}\left(  g\right)  ,
\]
where
\begin{equation}
\bar{z}_{k}\left(  g\left(  T\right)  \right)  :=\mathbb{E}\left(
z_{k}|\sigma\left(  g\left(  T\right)  \right)  \right)  . \label{e.6.6}%
\end{equation}

\end{cor}

\begin{proof}
This is a special case of Corollary \ref{c.5.6}, with $k\left(  t\right)
:=\frac{t}{T}\left(  A,a\right)  $ and $F\left(  g\right)  :=f\left(  g\left(
T\right)  \right)  .$
\end{proof}

\begin{cor}
[Left heat kernel integration by parts]\label{c.6.5}Let $k:=\left(
A,a\right)  \in\mathfrak{g}_{CM}$ and suppose that $f:G\rightarrow\mathbb{C}$
is a smooth function such that $f$ and $\tilde{k}f$ are polynomially bounded.
Then%
\[
\int_{G}\left(  \tilde{k}f\right)  d\nu_{T}\left(  g\right)  =\int_{G}f\left(
g\right)  \bar{z}_{k}^{l}\left(  g\right)  d\nu_{T}\left(  g\right)  ,
\]
where $\tilde{k}f\left(  g\right)  :=\frac{d}{ds}\Big|_{0}f\left(  g\left(
sk\right)  \right)  $ and%
\begin{equation}
\bar{z}_{k}^{l}\left(  g\right)  =-\bar{z}_{k}\left(  g^{-1}\right)
\label{e.6.7}%
\end{equation}
where $\bar{z}_{k}$ is defined in Equation \eqref{e.6.6}.
\end{cor}

\begin{proof}
Let $u\left(  g\right)  :=f\left(  g^{-1}\right)  $ so that $f\left(
g\right)  =u\left(  g^{-1}\right)  .$ Then
\[
\left(  \tilde{k}f\right)  \left(  g\right)  =\frac{d}{ds}|_{0}f\left(
g\cdot\left(  sk\right)  \right)  =\frac{d}{ds}|_{0}u\left(  \left(
-sk\right)  \cdot g^{-1}\right)  =-\left(  \hat{k}u\right)  \left(
g^{-1}\right)  .
\]
Therefore by Corollary \ref{c.6.4} and two uses of Corollary \ref{c.4.9} we
find%
\begin{align*}
\int_{G}\left(  \tilde{k}f\right)  d\nu_{T}\left(  g\right)   &  =-\int
_{G}\left(  \hat{k}u\right)  \left(  g^{-1}\right)  d\nu_{T}\left(  g\right)
=-\int_{G}\left(  \hat{k}u\right)  \left(  g\right)  d\nu_{T}\left(  g\right)
\\
&  =-\int_{G}u\left(  g\right)  \,\bar{z}_{k}\left(  g\right)  \,d\nu
_{T}\left(  g\right)  =-\int_{G}f\left(  g^{-1}\right)  \,\bar{z}_{k}\left(
g\right)  \,d\nu_{T}\left(  g\right) \\
&  =-\int_{G}f\left(  g\right)  \,\bar{z}_{k}\left(  g^{-1}\right)  \,d\nu
_{T}\left(  g\right)  .
\end{align*}

\end{proof}

\begin{df}
\label{d.6.6}A \textbf{cylinder polynomial} is a cylinder function,
$f=F\circ\pi_{P}:G\rightarrow\mathbb{C},$ where $P\in\operatorname*{Proj}%
\left(  W\right)  $ and $F$ is a real or complex polynomial function on
$PH\times\mathbf{C.}$
\end{df}

\begin{cor}
[Closability of the Dirichlet Form]\label{c.6.7} Given real--valued
cylindrical polynomials, $u,v$ on $G,$ let
\[
\mathcal{E}_{T}^{0}\left(  u,v\right)  :=\int_{G}\left\langle
\operatorname*{grad}u,\operatorname*{grad}v\right\rangle _{H}d\nu_{T},
\]
where $\operatorname*{grad}u:G\rightarrow\mathfrak{g}_{CM}$ is the gradient of
$u$ defined by%
\[
\left\langle \operatorname*{grad}u,k\right\rangle _{\mathfrak{g}_{CM}}%
=\tilde{k}u\text{ for all }k\in\mathfrak{g}_{CM}.
\]
Then $\mathcal{E}_{T}^{0}$ is closable and its closure, $\mathcal{E}_{T},$ is
a Dirichlet form on $\operatorname{Re}L^{2}\left(  G,\nu_{T}\right)  .$
\end{cor}

\begin{proof}
The closability of $\mathcal{E}_{T}^{0}$ is equivalent to the closability of
the gradient operator, $\operatorname{grad}:L^{2}\left(  \nu_{T}\right)
\rightarrow L^{2}\left(  \nu_{T}\right)  \otimes\mathfrak{g}_{CM},$ with the
domain, $\mathcal{D}\left(  \operatorname*{grad}\right)  ,$ being the space of
cylinder polynomials on $G.$ To check the latter statement it suffices to show
that $\operatorname*{grad}$ has a densely defined adjoint which is easily
accomplished. Indeed, if $k\in\mathfrak{g}_{CM}$ and $u$ and $v$ are cylinder
polynomials, then
\begin{align*}
\left\langle \operatorname*{grad}u,v\cdot k\right\rangle _{L^{2}\left(
\nu_{T}\right)  \otimes\mathfrak{g}_{CM}}  &  =\int_{G}\tilde{k}u\cdot
v~d\nu_{T}\\
&  =\int_{G}\left[  \tilde{k}\left(  u\cdot v\right)  -u\cdot\tilde
{k}v\right]  ~d\nu_{T}\\
&  =\left\langle u,-\tilde{k}v+\bar{z}_{k}^{l}v\right\rangle _{L^{2}\left(
\nu_{T}\right)  },
\end{align*}
wherein we have used the product rule in the second equality and Corollary
\ref{c.6.5} for the third. This shows that $v\cdot k$ is contained in the
domain of $\operatorname*{grad}^{\ast}$ and $\operatorname*{grad}^{\ast
}\left(  v\cdot k\right)  =-\tilde{k}v+\bar{z}_{k}^{l}v,$ where $z_{k}^{l}$ is
as in Eq. (\ref{e.6.7}). This completes the proof since linear combination of
functions of the form $v\cdot k$ with $k\in\mathfrak{g}_{CM}$ and $v$ being a
cylinder polynomial is dense in $L^{2}\left(  \nu_{T}\right)  \otimes
\mathfrak{g}_{CM}.$
\end{proof}

\section{The Ricci Curvature on Heisenberg type groups\label{s.7}}

In this section we compute the Ricci curvature for $G\left(  \omega\right)  $
and its finite-dimensional approximations. This information will be used in
Section \ref{s.8} to prove a logarithmic Sobolev inequality for $\nu_{T}$ and
to get detailed $L^{p}$ --bounds on the Radon-Nikodym derivatives of $\nu_{T}$
under translations by elements from $G_{CM}.$

\begin{nota}
\label{n.7.1}Let $\left(  W,H,\omega\right)  $ be as in Notation \ref{n.3.1},
$P\in\operatorname*{Proj}\left(  W\right)  ,$ and $G_{P}=PW\times
\mathbf{C}\subset G_{CM}$ as in Notation \ref{n.3.24}. We equip $G_{P}$ with
the left invariant Riemannian metric induced from restriction of the (real
part of the) inner product on $\mathfrak{g}_{CM}=H\times\mathbf{C}$ to
$\operatorname*{Lie}\left(  G_{P}\right)  =PH\times\mathbf{C}.$ Further, let
$\operatorname{Ric}^{P}$ denote the associated Ricci tensor at the identity in
$G_{P}.$
\end{nota}

\begin{prop}
\label{p.7.2}If $\left(  W,H,\omega,P\right)  $ as in Notation \ref{n.7.1},
$P\in\operatorname*{Proj}\left(  W\right)  $ is as in Eq. (\ref{e.3.42}), and
$\left(  A,a\right)  \in PH\times\mathbf{C},$\textrm{ }then%
\begin{align}
\left\langle \operatorname{Ric}^{P}\left(  A,a\right)  ,\left(  A,a\right)
\right\rangle _{H\times\mathbf{C}}  &  =\frac{1}{4}\sum_{j,k=1}^{n}\left\vert
\left\langle \omega\left(  e_{k},e_{j}\right)  ,a\right\rangle _{\mathbf{C}%
}\right\vert ^{2}-\frac{1}{2}\sum_{k=1}^{n}\left\Vert \omega\left(
A,e_{k}\right)  \right\Vert _{\mathbf{C}}^{2}\label{e.7.1}\\
&  =\frac{1}{4}\left\Vert \left\langle \omega\left(  \cdot,\cdot\right)
,a\right\rangle _{\mathbf{C}}\right\Vert _{\left(  PH\right)  ^{\ast}%
\otimes\left(  PH\right)  ^{\ast}}^{2}-\frac{1}{2}\left\Vert \omega\left(
A,\cdot\right)  \right\Vert _{\left(  PH\right)  ^{\ast}\otimes\mathbf{C}}%
^{2}. \label{e.7.2}%
\end{align}

\end{prop}

\begin{proof}
We are going to compute $\operatorname{Ric}^{P}$ using the formula in Equation
\eqref{e.11.3} of Appendix \ref{s.11}. If $\left\{  f_{\ell}\right\}  _{\ell
=1}^{\dim\mathbf{C}}$ is an orthonormal basis for $\mathbf{C},$ then%
\begin{equation}
\sum_{k=1}^{n}\left\Vert ad_{\left(  e_{k},0\right)  }\left(  A,a\right)
\right\Vert _{H\times\mathbf{C}}^{2}+\sum_{\ell=1}^{\dim\mathbf{C}}\left\Vert
ad_{\left(  0,f_{\ell}\right)  }\left(  A,a\right)  \right\Vert _{H\times
\mathbf{C}}^{2}=\sum_{k=1}^{n}\left\Vert \omega\left(  e_{k},A\right)
\right\Vert _{\mathbf{C}}^{2}. \label{e.7.3}%
\end{equation}
If $\left(  B,b\right)  \in PH\times\mathbf{C},$ then%
\begin{align*}
ad_{\left(  B,b\right)  }^{\ast}\left(  A,a\right)  =  &  \sum_{j=1}%
^{n}\left\langle ad_{\left(  B,b\right)  }^{\ast}\left(  A,a\right)  ,\left(
e_{j},0\right)  \right\rangle _{\mathfrak{g}_{CM}}\left(  e_{j},0\right) \\
&  +\sum_{\ell=1}^{\dim\mathbf{C}}\left\langle ad_{\left(  B,b\right)  }%
^{\ast}\left(  A,a\right)  ,\left(  0,f_{\ell}\right)  \right\rangle
_{\mathfrak{g}_{CM}}\left(  0,f_{\ell}\right) \\
=  &  \sum_{j=1}^{n}\left\langle \left(  A,a\right)  ,\left[  \left(
B,b\right)  ,\left(  e_{j},0\right)  \right]  \right\rangle _{\mathfrak{g}%
_{CM}}\left(  e_{j},0\right) \\
&  +\sum_{\ell=1}^{\dim\mathbf{C}}\left\langle \left(  A,a\right)  ,\left[
\left(  B,b\right)  ,\left(  0,f_{\ell}\right)  \right]  \right\rangle
_{\mathfrak{g}_{CM}}\left(  0,f_{\ell}\right) \\
=  &  \sum_{j=1}^{n}\left\langle \left(  A,a\right)  ,\left(  0,\omega\left(
B,e_{j}\right)  \right)  \right\rangle _{\mathfrak{g}_{CM}}\left(
e_{j},0\right)  =\sum_{j=1}^{n}\left(  a,\omega\left(  B,e_{j}\right)
\right)  _{\mathbf{C}}\left(  e_{j},0\right)  .
\end{align*}
This then immediately implies
\begin{equation}
\sum_{k=1}^{n}\left\Vert ad_{\left(  e_{k},0\right)  }^{\ast}\left(
A,a\right)  \right\Vert _{\mathfrak{g}_{CM}}^{2}+\sum_{\ell=1}^{\dim
\mathbf{C}}\left\Vert ad_{\left(  0,f_{\ell}\right)  }^{\ast}\left(
A,a\right)  \right\Vert _{\mathfrak{g}_{CM}}^{2}=\sum_{k=1}^{n}\sum_{j=1}%
^{n}\left\langle a,\omega\left(  e_{k},e_{j}\right)  \right\rangle
_{\mathbf{C}}^{2}. \label{e.7.4}%
\end{equation}
Using Equations \eqref{e.7.3} and \eqref{e.7.4} with the formula for the Ricci
tensor in Equation \eqref{e.11.3} of Appendix \ref{s.11} implies Equation \eqref{e.7.1}.
\end{proof}

\begin{cor}
\label{c.7.3}For $P\in\operatorname*{Proj}\left(  W\right)  $ as in
\eqref{e.3.42}, let%
\begin{equation}
k_{P}\left(  \omega\right)  :=-\frac{1}{2}\sup\left\{  \left\Vert
\omega\left(  \cdot,A\right)  \right\Vert _{\left(  PH\right)  ^{\ast}%
\otimes\mathbf{C}}^{2}:A\in PH,\ \left\Vert A\right\Vert _{PH}=1\right\}  .
\label{e.7.5}%
\end{equation}
Also let
\begin{equation}
k\left(  \omega\right)  :=-\frac{1}{2}\sup_{\left\Vert A\right\Vert _{H}%
=1}\left\Vert \omega\left(  \cdot,A\right)  \right\Vert _{H^{\ast}%
\otimes\mathbf{C}}^{2}\geqslant-\frac{1}{2}\left\Vert \omega\right\Vert
_{2}^{2}>-\infty. \label{e.7.6}%
\end{equation}
Then $k_{P}\left(  \omega\right)  $ is the largest constant $k\in\mathbb{R}$
such that
\begin{equation}
\left\langle \operatorname{Ric}^{P}\left(  A,a\right)  ,\left(  A,a\right)
\right\rangle _{PH\times\mathbf{C}}\geqslant k\left\Vert \left(  A,a\right)
\right\Vert _{PH\times\mathbf{C}}^{2}\text{~}\text{ for all }\text{ }\left(
A,a\right)  \in PH\times\mathbf{C} \label{e.7.7}%
\end{equation}
and $k\left(  \omega\right)  $ is the largest constant $k\in\mathbb{R}$ such
that Equation \eqref{e.7.7} holds uniformly for all $P\in\operatorname*{Proj}%
\left(  W\right)  .$
\end{cor}

\begin{proof}
Let us observe that by Equation \eqref{e.7.1}
\[
\frac{\left\langle \operatorname{Ric}^{P}\left(  A,a\right)  ,\left(
A,a\right)  \right\rangle _{PH\times\mathbf{C}}}{\left\Vert \left(
A,a\right)  \right\Vert _{PH\times\mathbf{C}}^{2}}\geqslant\frac{\left\langle
\operatorname{Ric}^{P}\left(  A,0\right)  ,\left(  A,0\right)  \right\rangle
_{PH\times\mathbf{C}}}{\left\Vert \left(  A,0\right)  \right\Vert
_{PH\times\mathbf{C}}^{2}}%
\]
the optimal lower bound, $k_{P}\left(  \omega\right)  ,$ for
$\operatorname{Ric}^{p}$ is determined by
\begin{align*}
k_{P}\left(  \omega\right)   &  =\inf_{A\in PH\setminus\left\{  0\right\}
}\frac{\left\langle \operatorname{Ric}^{P}\left(  A,0\right)  ,\left(
A,0\right)  \right\rangle _{PH\times\mathbf{C}}}{\left\Vert \left(
A,0\right)  \right\Vert _{PH\times\mathbf{C}}^{2}}\\
&  =\inf_{A\in PH\setminus\left\{  0\right\}  }\left(  -\frac{1}{2}%
\frac{\left\Vert \omega\left(  \cdot,A\right)  \right\Vert _{\left(
PH\right)  ^{\ast}\otimes\mathbf{C}}^{2}}{\left\Vert A\right\Vert _{PH}^{2}%
}\right)
\end{align*}
which is equivalent to Equation \eqref{e.7.5}. It is now a simple matter to
check that $k\left(  \omega\right)  =\inf_{P\in\operatorname*{Proj}\left(
W\right)  }k_{P}\left(  \omega\right)  $ which is the content of the last
assertion of the theorem.
\end{proof}

In revisiting the examples from Section \ref{s.3.3} we will have a number of
cases where $H$ and $\mathbf{C}$ are complex Hilbert spaces and $\omega
:H\times H\rightarrow\mathbf{C}$ will be a complex bilinear form. In these
cases it will be convenient to express the Ricci curvature in terms of these
complex structures.

\begin{prop}
\label{p.7.4}Suppose that $H$ and $\mathbf{C}$ are complex Hilbert spaces,
$\omega:H\times H\rightarrow\mathbf{C}$ is complex bi-linear, and
$P:H\rightarrow H$ is a finite rank (complex linear) orthogonal projection. We
make $G_{P}=PH\times\mathbf{C}$ into a Lie group using the group law in
Equation \eqref{e.3.6}. Letting and endow $G_{P}$ with the left invariant
Riemannian metric which agrees with $\left\langle \cdot,\cdot\right\rangle
_{\left[  \mathfrak{g}_{P}\right]  _{\operatorname{Re}}}:=\operatorname{Re}%
\left\langle \cdot,\cdot\right\rangle _{\mathfrak{g}_{P}}$ on $\mathfrak{g}%
_{P}=PH\times\mathbf{C}$ at the identity in $G_{P}.$ Then for all $\left(
A,a\right)  \in\mathfrak{g}_{P},$%
\begin{align}
\left\langle \operatorname{Ric}^{P}\left(  A,a\right)  , \left(  A,a\right)
\right\rangle _{\left[  \mathfrak{g}_{P}\right]  _{\operatorname{Re}}}  &
=\frac{1}{2}\left\Vert \left\langle \omega\left(  \cdot,\cdot\right)
,a\right\rangle _{\mathbf{C}}\right\Vert _{\left(  PH\right)  ^{\ast}%
\otimes\left(  PH\right)  ^{\ast}}^{2}-\left\Vert \omega\left(  A,\cdot
\right)  \right\Vert _{\left(  PH\right)  ^{\ast}\otimes\mathbf{C}}%
^{2}\label{e.7.8}\\
&  =\frac{1}{2}\sum_{j,k=1}^{n}\left\vert \left\langle \omega\left(
e_{k},e_{j}\right)  , a\right\rangle _{\mathbf{C}}\right\vert ^{2}-\sum
_{k=1}^{n}\left\Vert \omega\left(  A,e_{k}\right)  \right\Vert _{\mathbf{C}%
}^{2}, \label{e.7.9}%
\end{align}
where $\left\{  e_{j}\right\}  _{j=1}^{n}$ is any orthonormal basis for $PH.$
\end{prop}

\begin{proof}
Applying Equation \eqref{e.7.2} with $PH,$ $\mathbf{C},$ and $\mathfrak{g}%
_{P}$ being replaced by $\left(  PH\right)  _{\operatorname{Re}},$
$\mathbf{C}_{\operatorname{Re}},$ and $\left[  \mathfrak{g}_{P}\right]
_{\operatorname{Re}}$ implies%
\begin{multline*}
\left\langle \operatorname{Ric}^{P}\left(  A,a\right)  ,\left(  A,a\right)
\right\rangle _{\left[  \mathfrak{g}_{P}\right]  _{\operatorname{Re}}}%
=\frac{1}{4}\left\Vert \left\langle \omega\left(  \cdot,\cdot\right)
,a\right\rangle _{\mathbf{C}_{\operatorname{Re}}}\right\Vert _{\left(
PH\right)  _{\operatorname{Re}}^{\ast}\otimes\left(  PH\right)
_{\operatorname{Re}}^{\ast}}^{2}\\
-\frac{1}{2}\left\Vert \omega\left(  A,\cdot\right)  \right\Vert _{\left(
PH\right)  _{\operatorname{Re}}^{\ast}\otimes\mathbf{C}_{\operatorname{Re}}%
}^{2}.
\end{multline*}
However, by Lemma \ref{l.3.17} we also know that%
\[
\left\Vert \left\langle \omega\left(  \cdot,\cdot\right)  ,a\right\rangle
_{\mathbf{C}_{\operatorname{Re}}}\right\Vert _{\left(  PH\right)
_{\operatorname{Re}}^{\ast}\otimes\left(  PH\right)  _{\operatorname{Re}%
}^{\ast}}^{2}=2\left\Vert \left\langle \omega\left(  \cdot,\cdot\right)
,a\right\rangle _{\mathbf{C}}\right\Vert _{\left(  PH\right)  ^{\ast}%
\otimes\left(  PH\right)  ^{\ast}}^{2}%
\]
and%
\[
\left\Vert \omega\left(  A,\cdot\right)  \right\Vert _{\left(  PH\right)
_{\operatorname{Re}}^{\ast}\otimes\mathbf{C}_{\operatorname{Re}}}%
^{2}=2\left\Vert \omega\left(  A,\cdot\right)  \right\Vert _{\left(
PH\right)  ^{\ast}\otimes\mathbf{C}}^{2}%
\]
which completes the proof of Equation \eqref{e.7.8}.
\end{proof}

\begin{rem}
\label{r.7.5}By letting $n\rightarrow\infty$ in Propositions \ref{p.7.2} and
\ref{p.7.4}, it is reasonable to interpret the Ricci tensor on $G_{CM}$ to be
determined by
\begin{align}
\left\langle \operatorname{Ric}\mathrm{\,}\left(  A,a\right)  ,\left(
A,a\right)  \right\rangle _{\left[  \mathfrak{g}_{CM}\right]
_{\operatorname{Re}}}  &  =\alpha_{\mathbb{F}}\left(  \frac{1}{4}\left\Vert
\left\langle a,\omega\left(  \cdot,\cdot\right)  \right\rangle _{\mathbf{C}%
}\right\Vert _{H^{\ast}\otimes H^{\ast}}^{2}-\frac{1}{2}\left\Vert
\omega\left(  \cdot,A\right)  \right\Vert _{H^{\ast}\otimes\mathbf{C}}%
^{2}\right) \label{e.7.10}\\
&  =\alpha_{\mathbb{F}}\left(  \frac{1}{4}\sum_{j,k=1}^{\infty}\left\vert
\left\langle a,\omega\left(  e_{k},e_{j}\right)  \right\rangle _{\mathbf{C}%
}\right\vert ^{2}-\frac{1}{2}\sum_{k=1}^{\infty}\left\Vert \omega\left(
e_{k},A\right)  \right\Vert _{\mathbf{C}}^{2}\right)  , \label{e.7.11}%
\end{align}
where $\left\{  e_{j}\right\}  _{j=1}^{\infty}$ is an orthonormal basis for
$H,$ $\mathbb{F}$ is either $\mathbb{R}$ or $\mathbb{C}$ and $\alpha
_{\mathbb{F}}$ is one or two respectively. Moreover if $\mathbf{C}%
=\mathbb{F},$ then Equation \eqref{e.7.10} may be written as%
\begin{equation}
\left\langle \operatorname{Ric}\mathrm{\,}\left(  A,a\right)  ,\left(
A,a\right)  \right\rangle _{\left[  \mathfrak{g}_{CM}\right]
_{\operatorname{Re}}}=\alpha_{\mathbb{F}}\left(  \frac{1}{4}\left\Vert
\omega\left(  \cdot,\cdot\right)  \right\Vert _{H^{\ast}\otimes H^{\ast}}%
^{2}\cdot\left\vert a\right\vert ^{2}-\frac{1}{2}\left\Vert \omega\left(
\cdot,A\right)  \right\Vert _{H^{\ast}}^{2}\right)  . \label{e.7.12}%
\end{equation}

\end{rem}

\subsection{Examples revisited\label{s.7.1}}

Using Equation \eqref{e.7.10}, it is straight forward to compute the Ricci
tensor on $G$ for each of the Examples \ref{ex.3.18} -- \ref{ex.3.23}.

\begin{lem}
\label{l.7.6}The Ricci tensor for $G_{CM}$ associated to each of the
structures introduced in Examples \ref{ex.3.18} and \ref{ex.3.19} are given
(respectively) by%
\begin{equation}
\left\langle \operatorname{Ric}\left(  z,c\right)  , \left(  z,c\right)
\right\rangle _{\mathfrak{h}_{\mathbb{R}}^{n}}=\frac{n c^{2}}{2}-\frac{1}%
{2}\left\Vert z\right\Vert _{\mathbb{C}^{n}}^{2}~\text{ for all }\left(
z,c\right)  \in\mathbb{C}^{n}\times\mathbb{R}, \label{e.7.13}%
\end{equation}
and
\begin{equation}
\left\langle \operatorname{Ric}\left(  z,c\right)  , \left(  z,c\right)
\right\rangle _{\left[  \mathfrak{h}_{\mathbb{C}}^{n}\right]
_{\operatorname{Re}}}=n\left\vert c\right\vert ^{2}-\left\Vert z\right\Vert
_{\mathbb{C}^{2n}}^{2}~\text{ for all }\left(  z,c\right)  \in\mathbb{C}%
^{2n}\times\mathbb{C}. \label{e.7.14}%
\end{equation}

\end{lem}

\begin{proof}
We omit the proof of this lemma as it can be deduced from the next proposition
by taking $Q=I.$
\end{proof}

\begin{prop}
\label{p.7.7}The Ricci tensor for $G_{CM}$ associated to each of the
structures introduced in Examples \ref{ex.3.20} and \ref{ex.3.21} are given
(respectively) by%
\begin{equation}
\left\langle \operatorname{Ric}\left(  h,c\right)  , \left(  h,c\right)
\right\rangle _{\mathfrak{g}_{CM}}=\frac{1}{2}\left[  c^{2}\operatorname{tr}%
Q^{2} -\left\Vert Qh\right\Vert _{H}^{2}\right]  ~\text{ for all }\left(
h,c\right)  \in H\times\mathbb{R}, \label{e.7.15}%
\end{equation}
and
\begin{equation}
\left\langle \operatorname{Ric}\left(  k_{1},k_{2},c\right)  , \left(
k_{1},k_{2},c\right)  \right\rangle _{\left[  \mathfrak{g}_{CM}\right]
_{\operatorname{Re}}}=\left\vert c\right\vert ^{2}\operatorname{tr}%
Q^{2}-\left\Vert Qk_{1}\right\Vert _{K}^{2}-\left\Vert Qk_{2}\right\Vert
_{K}^{2} \label{e.7.16}%
\end{equation}
for all $\left(  k_{1},k_{2},c\right)  \in K\times K\times\mathbb{R}.$
\end{prop}

\begin{proof}
We start with the proof of Equation \eqref{e.7.15}. In this case,%
\begin{align*}
\left\Vert \omega\left(  \cdot,h\right)  \right\Vert _{H_{\operatorname{Re}%
}^{\ast}}^{2}  &  =\sum_{j=1}^{\infty}\left[  \omega\left(  e_{j},A\right)
^{2}+\omega\left(  ie_{j},A\right)  ^{2}\right] \\
&  =\sum_{j=1}^{\infty}\left[  \left(  \operatorname{Im}\left\langle
h,e_{j}\right\rangle _{Q}\right)  ^{2}+\left(  \operatorname{Im}\left\langle
h,ie_{j}\right\rangle _{Q}\right)  ^{2}\right] \\
&  =\sum_{j=1}^{\infty}\left[  \left(  \operatorname{Im}\left\langle
h,e_{j}\right\rangle _{Q}\right)  ^{2}+\left(  \operatorname{Re}\left\langle
h,e_{j}\right\rangle _{Q}\right)  ^{2}\right] \\
&  =\sum_{j=1}^{\infty}\left\vert \left\langle h,e_{j}\right\rangle
_{Q}\right\vert ^{2}=\left\Vert Qh\right\Vert _{H}^{2}%
\end{align*}
and from Equation \eqref{e.3.36} $\left\Vert \omega\right\Vert _{2}%
^{2}=2\operatorname{tr}\left(  Q^{2}\right)  .$ Using these results in
Equation \eqref{e.7.12} with $\mathbb{F}=\mathbb{R}$ gives Equation
\eqref{e.7.15} with $\mathbb{F}=\mathbb{C}$ and $H=K\times K,$ Equation
\eqref{e.7.16} follows from Equation \eqref{e.7.12} with $\mathbb{F}%
=\mathbb{C},$ Equation \eqref{e.3.36}, and the following identity;%
\begin{align}
\left\Vert \omega\left(  \left(  k_{1},k_{2}\right)  ,\cdot\right)
\right\Vert _{H^{\ast}}^{2}  &  =\sum_{j=1}^{\infty}\left(  \left\vert
\omega\left(  \left(  k_{1},k_{2}\right)  ,\left(  e_{j},0\right)  \right)
\right\vert ^{2}+\left\vert \omega\left(  \left(  k_{1},k_{2}\right)  ,\left(
0,e_{j}\right)  \right)  \right\vert ^{2}\right) \nonumber\\
&  =\sum_{j=1}^{\infty}\left\vert \left\langle k_{2},Q\bar{e}_{j}\right\rangle
\right\vert ^{2}+\sum_{j=1}^{\infty}\left\vert \left\langle k_{1},Q\bar{e}%
_{j}\right\rangle \right\vert ^{2}\nonumber\\
&  =\left\Vert Qk_{1}\right\Vert _{K}^{2}+\left\Vert Qk_{2}\right\Vert
_{K}^{2}. \label{e.7.17}%
\end{align}

\end{proof}

\begin{prop}
\label{p.7.8}The Ricci tensor for $G_{CM}$ associated to the structure
introduced in Example \ref{ex.3.22} is given by%
\begin{equation}
\left\langle \operatorname{Ric}\left(  v,c\right)  ,\left(  v,c\right)
\right\rangle _{\left[  \mathfrak{g}_{CM}\right]  _{\operatorname{Re}}}%
=\sum_{j=1}^{\infty}q_{j}^{2}\left\langle \operatorname{Ric}^{\alpha}\left(
v_{j},c\right)  ,\left(  v_{j},c\right)  \right\rangle _{V_{\operatorname{Re}%
}\times\mathbf{C}_{\operatorname{Re}}}~\forall~\left(  v,c\right)  \in
H\times\mathbb{F},\label{e.7.18}%
\end{equation}
where $\operatorname{Ric}^{\alpha}$ denotes the Ricci tensor on $G\left(
\alpha\right)  :=V\times\mathbf{C}$ as is defined by Equation \eqref{e.7.19} below.
\end{prop}

\begin{proof}
Using Equation \eqref{e.3.38} along with the identity,%
\begin{align*}
\left\Vert \omega\left(  \cdot,v\right)  \right\Vert _{H^{\ast}\otimes
\mathbf{C}}^{2} &  =\sum_{j=1}^{\infty}\sum_{a=1}^{d}\left\Vert \omega\left(
u_{a}\left(  j\right)  ,v\right)  \right\Vert _{\mathbf{C}}^{2}=\sum
_{j=1}^{\infty}\sum_{a=1}^{d}q_{j}^{2}\left\Vert \alpha\left(  u_{a}%
,v_{j}\right)  \right\Vert _{\mathbf{C}}^{2}\\
&  =\sum_{j=1}^{\infty}q_{j}^{2}\left\Vert \alpha\left(  \cdot,v_{j}\right)
\right\Vert _{V^{\ast}\otimes\mathbf{C}}^{2},
\end{align*}
in Equation \eqref{e.7.10} shows
\[
\left\langle \operatorname{Ric}\left(  v,c\right)  ,\left(  v,c\right)
\right\rangle _{\left[  \mathfrak{g}_{CM}\right]  _{\operatorname{Re}}}%
=\alpha_{\mathbb{F}}\sum_{j=1}^{\infty}q_{j}^{2}\left(  \frac{1}{4}\left\Vert
\left\langle \alpha\left(  \cdot,\cdot\right)  ,c\right\rangle _{\mathbf{C}%
}\right\Vert _{2}^{2}-\frac{1}{2}\left\Vert \alpha\left(  \cdot,v_{j}\right)
\right\Vert _{V^{\ast}\otimes\mathbf{C}}^{2}\right)  .
\]
Moreover, by a completely analogous finite-dimensional application of Equation
\eqref{e.7.10}, we have%
\begin{equation}
\left\langle \operatorname{Ric}^{\alpha}\left(  v_{j},c\right)  ,\left(
v_{j},c\right)  \right\rangle _{V_{\operatorname{Re}}\times\mathbf{C}%
_{\operatorname{Re}}}=\alpha_{\mathbb{F}}\left(  \frac{1}{4}\left\Vert
\left\langle \alpha\left(  \cdot,\cdot\right)  ,c\right\rangle _{\mathbf{C}%
}\right\Vert _{2}^{2}-\frac{1}{2}\left\Vert \alpha\left(  \cdot,v_{j}\right)
\right\Vert _{V^{\ast}\otimes\mathbf{C}}^{2}\right)  .\label{e.7.19}%
\end{equation}
Combining these two identities completes the proof.
\end{proof}

\begin{prop}
\label{p.7.9}Let $\alpha:V\times V\rightarrow\mathbf{C}$ be as in Example
\ref{ex.3.23}. For $v\in V,$ let $\alpha_{v}:V\rightarrow\mathbf{C}$ be
defined by $\alpha_{v}w=\alpha\left(  v,w\right)  $ and $\alpha_{v}^{\ast
}:\mathbf{C}\rightarrow V$ be its adjoint. The Ricci tensor for $G_{CM}$
associated to the structure introduced in Example \ref{ex.3.23} is then given
by%
\begin{align}
\left\langle \operatorname{Ric}\left(  h,c\right)  ,\left(  h,c\right)
\right\rangle _{\left[  \mathfrak{g}_{CM}\right]  _{\operatorname{Re}}}=  &
\frac{1}{2}\left[  \int_{\left[  0,1\right]  ^{2}}\left(  s\wedge t\right)
^{2}d\eta\left(  s\right)  d\bar{\eta}\left(  t\right)  \right]
\cdot\left\Vert \left\langle c,\alpha\left(  \cdot,\cdot\right)  \right\rangle
_{\mathbf{C}}\right\Vert _{V^{\ast}\otimes V^{\ast}}^{2}\label{e.7.20}\\
&  \qquad-\int_{\left[  0,1\right]  ^{2}}\left(  s\wedge t\right)
\operatorname{tr}\left(  \alpha_{h\left(  t\right)  }^{\ast}\alpha_{h\left(
s\right)  }\right)  d\eta\left(  s\right)  d\bar{\eta}\left(  t\right)
.\nonumber
\end{align}

\end{prop}

\begin{proof}
In this example we have%
\begin{align*}
\left\Vert \omega\left(  h,\cdot\right)  \right\Vert _{H^{\ast}\otimes
\mathbf{C}}^{2}  &  =\sum_{j=1}^{\infty}\sum_{a=1}^{d}\left\Vert \omega\left(
h,l_{j}u_{a}\right)  \right\Vert _{\mathbf{C}}^{2}\\
&  =\sum_{j=1}^{\infty}\sum_{a=1}^{d}\left\Vert \int_{0}^{1}\alpha\left(
h\left(  s\right)  ,l_{j}\left(  s\right)  u_{a}\right)  d\eta\left(
s\right)  \right\Vert _{\mathbf{C}}^{2}\\
&  =\sum_{j=1}^{\infty}\sum_{a=1}^{d}\left\langle \int_{0}^{1}\alpha\left(
h\left(  s\right)  ,l_{j}\left(  s\right)  u_{a}\right)  d\eta\left(
s\right)  ,\int_{0}^{1}\alpha\left(  h\left(  t\right)  ,l_{j}\left(
t\right)  u_{a}\right)  d\eta\left(  t\right)  \right\rangle _{\mathbf{C}}\\
&  =\int_{0}^{1}d\eta\left(  s\right)  \int_{0}^{1}d\bar{\eta}\left(
t\right)  \left(  s\wedge t\right)  \sum_{a=1}^{d}\left\langle \alpha
_{h\left(  s\right)  }u_{a},\alpha_{h\left(  t\right)  }u_{a}\right\rangle
_{\mathbf{C}}\\
&  =\int_{\left[  0,1\right]  ^{2}}s\wedge t\left[  \operatorname{tr}\left(
\alpha_{h\left(  t\right)  }^{\ast}\alpha_{h\left(  s\right)  }\right)
\right]  d\eta\left(  s\right)  d\bar{\eta}\left(  t\right)  .
\end{align*}
Using this identity along with Equation \eqref{e.3.40} in Equation
\eqref{e.7.10} with $\alpha_{\mathbb{F}}=\alpha_{\mathbb{C}}=2$ implies
Equation \eqref{e.7.20}.
\end{proof}

\section{Heat Inequalities\label{s.8}}

\subsection{Infinite-dimensional Radon-Nikodym derivative
estimates\label{s.8.1}}

Recall from Theorem \ref{t.6.1} and Corollary \ref{c.6.2}, we have already
shown that $\nu_{T}\circ l_{h}^{-1}$ and $\nu_{T}\circ r_{h}^{-1}$ are
absolutely continuous to $\nu_{T}$ for all $h\in G_{CM}$ and $T>0.$ These
results were based on the path space quasi-invariance formula given Theorem
\ref{t.5.2}. However, in light of the results in Malliavin
\cite{Malliavin1990} it is surprising that Theorem \ref{t.5.2} holds at all
and we do not expect it to extend to many other situations. Therefore, it is
instructive to give an independent proof of Theorem \ref{t.6.1} and Corollary
\ref{c.6.2} which will work for a much larger class of examples. The
alternative proof have the added advantage of giving detailed size estimates
on the resulting Randon-Nikodym derivatives.

\begin{thm}
\label{t.8.1}For all $h\in G_{CM}$ and $T>0,$ $\nu_{T}\circ l_{h}^{-1}$ and
$\nu_{T}\circ r_{h}^{-1}$ are absolutely continuous to $\nu_{T}.$ Let
$Z_{h}^{l}:=\frac{d\left(  \nu_{T}\circ l_{h}^{-1}\right)  }{d\nu_{T}}$ and
$Z_{h}^{r}:=\frac{d\left(  \nu_{T}\circ r_{h}^{-1}\right)  }{d\nu_{T}}$ be the
respective Randon-Nikodym derivatives, $k\left(  \omega\right)  $ is given in
Equation \eqref{e.7.6}, and%
\[
c\left(  t\right)  :=\frac{t}{e^{t}-1}~\text{ for all }t\in\mathbb{R}%
\]
with the convention that $c\left(  0\right)  =1.$ Then for all $1\leqslant
p<\infty,$ $Z_{h}^{l}$ and $Z_{h}^{r}\ $are both in $L^{p}\left(  \nu
_{T}\right)  $ and satisfy the estimate
\begin{equation}
\left\Vert Z_{h}^{\ast}\right\Vert _{L^{p}\left(  \nu_{T}\right)  }%
\leqslant\exp\left(  \frac{c\left(  k\left(  \omega\right)  T\right)  \left(
p-1\right)  }{2T}d_{G_{CM}}^{2}\left(  \mathbf{e},h\right)  \right)  ,
\label{e.8.1}%
\end{equation}
where $\ast=l$ or $\ast=r.$
\end{thm}

\begin{proof}
The proof of this theorem is an application Theorem 7.3 and Corollary 7.4 in
\cite{DG07a} on quasi-invariance of the heat kernel measures for inductive
limits of finite-dimensional Lie groups. In applying these results the reader
should take: $G_{0}=G_{CM},$ $A=\operatorname*{Proj}\left(  W\right)  ,$
$s_{P}:=\pi_{P},$ $\nu_{P}=\operatorname{Law}\left(  g_{P}\left(  T\right)
\right)  ,$ and $\nu=\nu_{T}=\operatorname{Law}\left(  g\left(  T\right)
\right)  .$ We now verify that the hypotheses \cite[Theorem 7.3]{DG07a} are
satisfied. These assumptions include a densness condition on the inductive
limit group, consistency of the heat kernel measures on finite-dimensional Lie
groups, uniform bound on the Ricci curvature, and finally that the length of a
path in the inductive limit group can be approximated by the lengths of paths
in finite-dimensional groups.

\begin{enumerate}
\item By Proposition \ref{p.3.10}, $\cup_{P\in\operatorname*{Proj}\left(
W\right)  }G_{P}$ is a dense subgroup of $G_{CM}.$

\item From Lemma \ref{l.4.7}, for any $\left\{  P_{n}\right\}  _{n=1}^{\infty
}\subset\operatorname*{Proj}\left(  W\right)  $ with $P_{n}|_{H}\uparrow
I_{H}$ and $f\in BC\left(  G,\mathbb{R}\right)  $ (the bounded continuous maps
from $G$ to $\mathbb{R)},$ we have
\[
\int_{G}fd\nu=\lim_{n\rightarrow\infty}\int_{G_{P_{n}}}\left(  f\circ
i_{P}\right)  d\nu_{P_{n}}.
\]

\item Corollary \ref{c.7.3} shows that $\operatorname{Ric}_{P}\geqslant
k\left(  \omega\right)  g_{P}$ for all $P\in\operatorname*{Proj}\left(
W\right)  .$

\item Lastly we have to verify that for any $P_{0}\in\operatorname*{Proj}%
\left(  W\right)  ,$ and $k\in C^{1}\left(  \left[  0,1\right]  ,G_{CM}%
\right)  $ with $k\left(  0\right)  =e,$ there exists an increasing sequence,
$\left\{  P_{n}\right\}  _{n=1}^{\infty}\subset\operatorname*{Proj}\left(
W\right)  $ such that $P_{0}\subset P_{n},$ $P_{n}\uparrow I$ on $H,$ and%
\begin{equation}
\ell_{G_{CM}}\left(  k\right)  =\lim_{n\rightarrow\infty}\ell_{G_{P_{n}}%
}\left(  \pi_{n}\circ k\right)  , \label{e.8.2}%
\end{equation}
where $\pi_{n}:=\pi_{P_{n}}$ and $\ell_{G_{CM}}\left(  k\right)  $ is the
length of $k$ (see Notation \ref{n.3.9} with $T=1$). However, with $k\left(
t\right)  =\left(  A\left(  t\right)  ,c\left(  t\right)  \right)  ,$ using
the dominated convergence theorem applied to the identity (see Equation
\eqref{e.3.21});
\begin{align*}
\ell_{G_{P_{n}}}\left(  \pi_{n}\circ k\right)   &  =\int_{0}^{1}\left\Vert
\pi_{n}\dot{k}\left(  t\right)  -\frac{1}{2}\left[  \pi_{n}k\left(  t\right)
,\pi_{n}\dot{k}\left(  t\right)  \right]  \right\Vert _{\mathfrak{g}_{CM}}dt\\
&  =\int_{0}^{1}\sqrt{\left\Vert P_{n}\dot{A}\left(  t\right)  \right\Vert
_{H}^{2}+\left\Vert \dot{c}\left(  t\right)  -\frac{1}{2}\omega\left(
P_{n}A\left(  t\right)  ,P_{n}\dot{A}\left(  t\right)  \right)  \right\Vert
_{\mathbf{C}}^{2}}dt
\end{align*}
shows Equation \eqref{e.8.2} holds for any such choice of $P_{n}|_{H}\uparrow
I_{H}$ with $P_{0}\subset P_{n}\in\operatorname*{Proj}\left(  W\right)  .$
\end{enumerate}
\end{proof}

\begin{rem}
\label{r.8.2} In the case of infinite-dimensional matrix groups three out of
four assumptions hold as has been shown in \cite{Gordina2000b}. The condition
that fails is the uniform bounds on the Ricci curvature which is one of the
main results in \cite{Gordina2005a}.
\end{rem}

\subsection{Logarithmic Sobolev Inequality\label{s.8.2}}

\begin{thm}
\label{t.8.3}Let $\left(  \mathcal{E}_{T},\mathcal{D}\left(  \mathcal{E}%
_{T}\right)  \right)  $ be the closed Dirichlet form in Corollary \ref{c.6.7}
and $k\left(  \omega\right)  $ be as in Equation \eqref{e.7.6}. Then for all
real-valued $f\in\mathcal{D}\left(  \mathcal{E}_{T}\right)  ,$ the following
logarithmic Sobolev inequality holds
\begin{equation}
\int_{G}\left(  f^{2}\ln f^{2}\right)  d\nu_{T}\leqslant2\frac{1-e^{-k\left(
\omega\right)  T}}{k\left(  \omega\right)  }\mathcal{E}_{T}\left(  f,f\right)
+\int_{G}f^{2}d\nu_{T}\cdot\ln\int_{G}f^{2}d\nu_{T}, \label{e.8.3}%
\end{equation}
where $\nu_{T}=\operatorname{Law}\left(  g\left(  T\right)  \right)  $ is the
heat kernel measure on $G$ as in Definition \ref{d.4.2}.
\end{thm}

\begin{proof}
Let $f:G\rightarrow\mathbb{R}$ be a cylinder polynomial as in Definition
\ref{d.6.6}. Following the method of Bakry and Ledoux applied to $G_{P}$ (see
\cite[Theorem 2.9]{Driver1996b} for the case needed here) shows
\begin{align}
\mathbb{E}\left[  \left(  f^{2}\log f^{2}\right)  \left(  g_{P}\left(
T\right)  \right)  \right]  \leqslant &  2\frac{1-e^{-k_{P}\left(
\omega\right)  T}}{k_{P}\left(  \omega\right)  }\mathbb{E}\left\Vert \left(
\operatorname*{grad}{}^{P}f\right)  \left(  g_{P}\left(  T\right)  \right)
\right\Vert _{G_{P}}^{2}\nonumber\\
&  \qquad+\mathbb{E}\left[  f^{2}\left(  g_{P}\left(  T\right)  \right)
\right]  \log\mathbb{E}\left[  f^{2}\left(  g_{P}\left(  T\right)  \right)
\right]  \label{e.8.4}%
\end{align}
where $k_{P}\left(  \omega\right)  $ is as in Equation \eqref{e.7.6}. Since
the function, $x\rightarrow x^{-1}\left(  1-e^{-x}\right)  ,$ is decreasing
and $k\left(  \omega\right)  \leqslant k_{P}\left(  \omega\right)  $ for all
$P\in\operatorname*{Proj}\left(  W\right)  ,$ Equation \eqref{e.8.4} also
holds with $k_{P}\left(  \omega\right)  $ replaced by $k\left(  \omega\right)
.$ With this observation along with Lemma \ref{l.4.7}, we may pass to the
limit at $P\uparrow I$ in Equation \eqref{e.8.4} to find%
\begin{align*}
\mathbb{E}\left[  \left(  f^{2}\log f^{2}\right)  \left(  g\left(  T\right)
\right)  \right]  \leqslant &  2\frac{1-e^{-k\left(  \omega\right)  T}%
}{k\left(  \omega\right)  }\mathbb{E}|\operatorname*{grad}f\left(  g\left(
T\right)  \right)  |^{2}\\
&  \qquad+\mathbb{E}\left[  f^{2}\left(  g\left(  T\right)  \right)  \right]
\log\mathbb{E}\left[  f^{2}\left(  g\left(  T\right)  \right)  \right]  .
\end{align*}
This is equivalent to Equation \eqref{e.8.3} when $f$ is a cylinder
polynomial. The result for general $f\in\mathcal{D}\left(  \mathcal{E}%
_{T}\right)  $ then holds by a standard (and elementary) limiting argument --
see the end of Example 2.7 in \cite{Gross93}.
\end{proof}

\section{Future directions\label{s.9}}

In this last section we wish to speculate on a number of ways that the results
in this paper might be extended.

\begin{enumerate}
\item It would be interesting to see what happens if we set $B_{0}$ to be
identically zero so that $g\left(  t\right)  $ in Equation \eqref{e.4.2}
becomes
\begin{equation}
g\left(  t\right)  =\left(  B\left(  t\right)  ,\frac{1}{2}\int_{0}^{t}%
\omega\left(  B\left(  \tau\right)  ,\dot{B}\left(  \tau\right)  \right)
d\tau\right)  . \label{e.9.1}%
\end{equation}
The generator now is $L=\frac{1}{2}\sum_{k=1}^{\infty}\widetilde{\left(
e_{k},0\right)  }^{2}$ and if $\omega\left(  \mathfrak{g}_{CM}\times
\mathfrak{g}_{CM}\right)  $ is a total subset of $\mathbf{C},$ $L$ would
satisfy H\"{o}rmander's condition for hypoellipticity. If $\dim H$ were
finite, it would follow that the heat kernel measure, $\nu_{T},$ is a smooth
positive measure and hence quasi-invariant. When $\dim H$ is infinite we do
not know if $\nu_{T}$ is still quasi-invariant. Certainly both proofs which
were given above when $B_{0}$ was not zero now break down.

\item It should be possible to remove the restriction on $\mathbf{C}$ being
finite-dimensional, i.e. we expect much of what we have done to go through
when $\mathbf{C}$ is a separable Hilbert space. In doing so one would have to
modify the finite-dimensional approximations used in the theory to truncate
$\mathbf{C}$ as well.

\item It should be possible to widen the class of admissible $\omega$s
substantially. The idea is to assume that $\omega$ is only defined from
$H\times H\rightarrow\mathbf{C}$ such that $\left\Vert \omega\right\Vert
_{2}<\infty.$ Under this relaxed assumption, we will no longer have a group
structure on $G:=W\times\mathbf{C}.$ Nevertheless, with a little work one can
still make sense of Brownian motion process defined in Definition \ref{d.4.2}
by letting
\begin{equation}
\int_{0}^{t}\omega\left(  B\left(  \tau\right)  ,dB\left(  \tau\right)
\right)  :=L^{2}\text{ -- }\lim_{n\rightarrow\infty}\int_{0}^{t}\omega\left(
P_{n}B\left(  \tau\right)  ,dP_{n}B\left(  \tau\right)  \right)  .
\label{e.9.2}%
\end{equation}
In fact, using Nelson's hypercontractivity and the fact that
\[
\int_{0}^{t}\omega\left(  P_{n}B\left(  \tau\right)  ,dP_{n}B\left(
\tau\right)  \right)
\]
is in the second homogeneous chaos subspace, the convergence in Equation
\eqref{e.9.2} is in $L^{p}$ for all $p\in\lbrack1,\infty).$ In this setting we
expect the path space quasi-invariance results of Section \ref{s.5} to remain
valid. Similarly, as the lower bound on the Ricci curvature only depends on
$\omega|_{H\times H},$ we expect the results of Section \ref{s.8.1} to go
through as well. As a consequence, $G$ should carry a measurable left and
right actions by element of $G_{CM}$ and these actions should leave the heat
kernel measures (end point distributions of the Brownian motion on $G)$
quasi-invariant. One might call the resulting structure a \textbf{quasigroup}.
Unfortunately, this term has already been used in abstract algebra.

\item We also expect that level of non-commutativity of $G$ may be increased.
To be more precise, under suitable hypotheses, it should be possible to handle
more general graded nilpotent Lie groups. However, when the level of
nilpotency of $G$ is increased, there will likely be trouble with the path
space quasi-invariance in section \ref{s.5}. Nevertheless, the methods of
Section \ref{s.8.1} should survive and therefore we still expect the heat
kernel measure to be quasi-invariant.
\end{enumerate}

\begin{acknowledgement}
The first author would like to thank the Berkeley mathematics department and
the Miller Institute for Basic Research in Science for their support of this
project in its latter stages.
\end{acknowledgement}

\appendix

\section{Wiener Space Results\label{s.10}}

The well known material presented in this Appendix may be (mostly) found in
the books \cite{Kuo75} and \cite{Bog98}. In particular, the following theorem
is based in part on Lemma 2.4.1 on p. 59 of \cite{Bog98}, and Theorem 3.9.6 on
p. 138 \cite{Bog98}.

\begin{thm}
\label{t.10.1}Let $(X,\mathcal{B}_{X},\mu)$ be a Gaussian measure space as in
Definition \ref{d.2.1}. Then

\begin{enumerate}
\item $\left(  H,\left\Vert \cdot\right\Vert _{H}\right)  $ is a normed space
such that%
\begin{equation}
\left\Vert h\right\Vert _{X}\leqslant\sqrt{C_{2}}\left\Vert h\right\Vert
_{H}~\text{ for all }h\in H, \label{e.10.1}%
\end{equation}
where $C_{2}$ is as in \eqref{e.2.2}.

\item Let $K$ be the closure of $X^{\ast}$ in $\operatorname{Re}L^{2}(\mu)$
and for $f\in K$ let
\[
\iota f:=h_{f}:=\int\limits_{X}xf(x)d\mu(x)\in X,
\]
where the integral is to be interpreted as a Bochner integral. Then
$\iota(K)=H$ and $\iota:K\rightarrow H$ is an isometric isomorphism of real
Banach spaces. Since $K$ is a real Hilbert space it follows that $\left\Vert
\cdot\right\Vert _{H}$ is a Hilbertian norm on $H.$

\item $H$ is a separable Hilbert space and%
\begin{equation}
(\iota u,h)_{H}=u(h)\text{ for all } u\in X^{\ast}\text{ and }h\in H.
\label{e.10.2}%
\end{equation}

\item The Cameron-Martin space, $H,$ is dense in $X.$

\item The quadratic form $q$ may be computed as
\begin{equation}
q(u,v)=\sum\limits_{i=1}^{\infty}u(e_{i})v(e_{i}) \label{e.10.3}%
\end{equation}
where $\{e_{i}\}_{i=1}^{\infty}$ is any orthonormal basis for $H$.
\end{enumerate}

Notice that by Item 1. $H\overset{i}{\hookrightarrow}X$ is continuous and
hence so is $X^{\ast}\overset{i^{tr}}{\hookrightarrow}H^{\ast}\cong
H=(\cdot,\cdot)_{H^{\ast}}.$ Equation \eqref{e.10.3} asserts that
\[
q=(\cdot,\cdot)_{H^{\ast}}\big|_{X^{\ast}\times X^{\ast}}.
\]

\end{thm}

\begin{proof}
1. Using Equation \eqref{e.2.6} we find
\[
\left\Vert h\right\Vert _{X}=\sup_{u\in X^{\ast}\setminus\left\{  0\right\}
}\frac{\left\vert u\left(  h\right)  \right\vert }{\left\Vert u\right\Vert
_{X^{\ast}}}\leqslant\sup_{u\in X^{\ast}\setminus\left\{  0\right\}  }%
\frac{\left\vert u\left(  h\right)  \right\vert }{\sqrt{q\left(  u,u\right)
/C_{2}}}\leqslant\sqrt{C_{2}}\left\Vert h\right\Vert _{H},
\]
and hence if $\left\Vert h\right\Vert _{H}=0$ then $\left\Vert h\right\Vert
_{X}=0$ and so $h=0.$ If $h,k\in H,$ then for all $u\in X^{\ast},$ $\left\vert
u(h)\right\vert \leqslant\left\Vert h\right\Vert _{H}\sqrt{q(u)}$ and
$\left\vert u(k)\right\vert \leqslant\left\Vert k\right\Vert _{H}\sqrt{q(u)}$
so that
\[
\left\vert u(h+k)\right\vert \leqslant\left\vert u(h)\right\vert +\left\vert
u(k)\right\vert \leqslant\left(  \left\Vert h\right\Vert _{H}+\left\Vert
k\right\Vert _{H}\right)  \sqrt{q(u)}.
\]
This shows $h+k\in H$ and $\left\Vert h+k\right\Vert _{H}\leqslant\left\Vert
h\right\Vert _{H}+\left\Vert k\right\Vert _{H}.$ Similarly, if $\lambda
\in\mathbb{R}$ and $h\in H,$ then $\lambda h\in H$ and $\left\Vert \lambda
h\right\Vert _{H}=\left\vert \lambda\right\vert \left\Vert h\right\Vert _{H}.$
Therefore $H$ is a subspace of $W$ and $\left(  H,\left\Vert \cdot\right\Vert
_{H}\right)  $ is a normed space.

2. For $f\in K$ and $u\in X^{\ast}$%
\begin{equation}
u\left(  \iota f\right)  =u\left(  \int\limits_{X}xf(x)d\mu(x)\right)
=\int\limits_{X}u(x)f(x)d\mu(x) \label{e.10.4}%
\end{equation}
and hence
\[
\left\vert u\left(  \iota f\right)  \right\vert \leqslant\left\Vert
u\right\Vert _{L^{2}(\mu)}\left\Vert f\right\Vert _{L^{2}(\mu)}=\sqrt
{q(u)}\left\Vert f\right\Vert _{K}%
\]
which shows that $\iota f\in H$ and $\left\Vert \iota f\right\Vert
_{H}\leqslant\left\Vert f\right\Vert _{K}.$ Moreover, by choosing $u_{n}\in
X^{\ast}$ such that $L^{2}(\mu)-\lim_{n\rightarrow\infty}u_{n}=f,$ we find
\[
\lim_{n\rightarrow\infty}\frac{\left\vert u_{n}\left(  \iota f\right)
\right\vert }{\sqrt{q(u_{n})}}=\lim_{n\rightarrow\infty}\frac{\left\vert
\int\limits_{X}u_{n}(x)f(x)d\mu(x)\right\vert }{\left\Vert u_{n}\right\Vert
_{L^{2}(\mu)}}=\frac{\left\Vert f\right\Vert _{L^{2}(\mu)}^{2}}{\left\Vert
f\right\Vert _{L^{2}(\mu)}}=\left\Vert f\right\Vert _{L^{2}(\mu)}%
\]
from which it follows $\left\Vert \iota f\right\Vert _{H}=\left\Vert
f\right\Vert _{K}.$ So we have shown that $\iota:K\rightarrow H$ is an
isometry. Let us now show that $\iota(K)=H.$ Given $h\in H$ and $u\in X^{\ast
}$ let $\hat{h}(u)=u(h).$ Since
\[
\left\vert \hat{h}(u)\right\vert \leqslant\left\vert u(h)\right\vert
\leqslant\sqrt{q(u)}\left\Vert h\right\Vert _{H}=\left\Vert u\right\Vert
_{L^{2}(\mu)}\left\Vert h\right\Vert _{H}=\left\Vert u\right\Vert
_{K}\left\Vert h\right\Vert _{H}%
\]
the functional $\hat{h}$ extends continuously to $K.$ We will continue to
denote this extension by $\hat{h}\in K^{\ast}.$ Since $K$ is a Hilbert space,
there exists $f\in K$ such that
\[
\hat{h}(u)=\int\limits_{X}f(x)u(x)d\mu(x)
\]
for all $u\in X^{\ast}$ (and in fact all $u\in K).$ Thus we have, for all
$u\in X^{\ast},$ that%
\[
u(h)=\int_{X}u(x)f(x)d\mu(x)=u\left(  \int_{X}xf(x)d\mu(x)\right)  =u(\iota
f).
\]
From this equation we conclude that $h=\iota f$ and hence $\iota(K)=H.$

3. $H$ is a separable since it is unitarily equivalent to $K\subset
L^{2}(X,\mathcal{B},\mu)$ and $L^{2}(X,\mathcal{B},\mu)$ is separable. Suppose
that $u\in X^{\ast},$ $f\in K$ and $h=\iota f\in H.$ Then
\begin{align*}
(\iota u,h)_{H}  &  =(\iota u,\iota f)_{H}=(u,f)_{K}\\
&  =\int_{X}u(x)f(x)d\mu(x)=u\left(  \int_{X}xf(x)d\mu(x)\right) \\
&  =u(\iota f)=u(h).
\end{align*}
4. For sake of contradiction, if $H\subset X$ were not dense, then, by the
Hahn--Banach theorem, there would exist $u\in X^{\ast}\setminus\{0\}$ such
that $u(H)=0.$ However from Equation \eqref{e.10.2}, we would then have%
\[
q(u,u)=(\iota u,\iota u)_{H}=u(\iota u)=0.
\]
Because we have assumed that $q$ to be an inner product on $X^{\ast},$ $u$
must be zero contrary to $u$ being in $X^{\ast}\setminus\left\{  0\right\}  .
$

5. Let $\left\{  e_{i}\right\}  _{i=1}^{\infty}$ be an orthonormal basis for
$H,$ then for $u,v\in X^{\ast},$
\begin{align*}
q(u,v)  &  =(u,v)_{K}=(\iota u,\iota v)_{H}=\sum_{i=1}^{\infty}(\iota
u,e_{i})_{H}(e_{i},\iota v)_{H}\\
&  =\sum\limits_{i=1}^{\infty}u(e_{i})v(e_{i})
\end{align*}
wherein the last equality we have again used Equation \eqref{e.10.2}.
\end{proof}

\section{The Ricci tensor on a Lie group\label{s.11}}

In this appendix we recall a formula for the Ricci tensor relative to a left
invariant Riemannian metric, $\left\langle \cdot,\cdot\right\rangle ,$ on any
finite-dimensional Lie Group, $G.$ Let $\nabla$ be the Levi-Civita covariant
derivative on $TG,$ for any $X\in\mathfrak{g}$ let $\tilde{X}\left(  g\right)
=l_{g\ast}X$ be the left invariant vector field on $G\ $such that $\tilde
{X}\left(  e\right)  =X.$, and for $X,Y\in\mathfrak{g},$ let $D_{X}%
Y:=\nabla_{X}\tilde{Y}\in\mathfrak{g}.$ Since $\nabla_{\tilde{X}}\tilde{Y}$ is
a left invariant vector field and $\left(  \nabla_{\tilde{X}}\tilde{Y}\right)
\left(  e\right)  =\nabla_{X}\tilde{Y}=D_{X}Y,$ we have the identity;
$\nabla_{\tilde{X}}\tilde{Y}=\widetilde{D_{X}Y}.$ Similarly the Ricci
curvature tensor, $\operatorname{Ric},$ (and more generally the full curvature
tensor) is invariant under left translations, i.e. $\operatorname{Ric}%
_{g}=l_{g^{-1}\ast}\operatorname{Ric}_{e}l_{g\ast}$ for all $g\in G.$ Hence it
suffices to compute the Ricci tensor at $e\in G.$ We will abuse notation and
simply write $\operatorname{Ric}$ for $\operatorname{Ric}_{e}.$

\begin{prop}
[The Ricci tensor on $G$]\label{p.11.1}Continuing the notation above, for all
$X,Y\in\mathfrak{g}$ we have
\begin{equation}
D_{X}Y:={{\frac{1}{2}}}\left(  [X,Y]-ad_{X}^{\ast}Y-ad_{Y}^{\ast}X\right)
\in\mathfrak{g}, \label{e.11.1}%
\end{equation}
where $ad_{X}^{\ast}$ denotes the adjoint of $ad_{X}$ relative to
$\left\langle \cdot,\cdot\right\rangle _{e}.$ We also have,%
\begin{equation}
\left\langle \operatorname{Ric}\mathrm{\,}X,X\right\rangle =\operatorname{tr}%
\left(  ad_{ad_{X}^{\ast}X}\right)  -\frac{1}{2}\operatorname*{tr}\left(
ad_{X}^{2}\right)  +\frac{1}{4}\sum_{Y\in\Gamma}\left\vert ad_{Y}^{\ast
}X\right\vert ^{2}-\frac{1}{2}\sum_{Y\in\Gamma}\left\vert ad_{Y}X\right\vert
^{2}, \label{e.11.2}%
\end{equation}
where $\Gamma\subset\mathfrak{g}$ is any orthonormal basis for $\mathfrak{g}.$
In particular if $\mathfrak{g}$ is nilpotent then$\operatorname{tr}\left(
ad_{ad_{X}^{\ast}X}\right)  =0$ and $\operatorname*{tr}\left(  ad_{X}%
^{2}\right)  =0$ and therefore Equation \eqref{e.11.2} reduces to
\begin{equation}
\left\langle \operatorname{Ric}\mathrm{\,}X,X\right\rangle =\frac{1}{4}%
\sum_{Y\in\Gamma}\left\vert ad_{Y}^{\ast}X\right\vert ^{2}-\frac{1}{2}%
\sum_{Y\in\Gamma}\left\vert ad_{Y}X\right\vert ^{2}\geqslant-\frac{1}{2}%
\sum_{Y\in\Gamma}\left\vert ad_{Y}X\right\vert ^{2}. \label{e.11.3}%
\end{equation}

\end{prop}

These results may be found in \cite{Besse1987}, see Lemma 7.27, Theorem 7.30,
and Corollary 7.33 for the computations of the Levi -Civita covariant
derivative, the curvature tensor, and the Ricci curvature tensor respectively.

\section{Proof of Theorem \ref{t.3.12}\label{s.12}}

Before giving the proof of Theorem \ref{t.3.12} it will be necessary to
develop the notion Carnot-Carath\'{e}odory distance function, $\delta,$ in
this infinite dimensional context.

\begin{nota}
\label{n.12.1}Let $T>0$ and $C_{HCM}^{1}$ denote the horizontal elements in
$C_{CM}^{1},$ where $g\in C_{CM}^{1}$ is \textbf{horizontal} iff $l_{g\left(
s\right)  ^{-1}\ast}g^{\prime}\left(  s\right)  \in H\times\left\{  0\right\}
$ for all $s.$ We then define,
\[
\delta(x,y)=\inf\left\{  \ell_{G_{CM}}\left(  g\right)  :g\in C_{HCM}%
^{1}\text{ such that }g\left(  0\right)  =x\text{ and }g\left(  T\right)
=y\right\}
\]
with the infimum of the empty set is taken to be infinite.
\end{nota}

Observe that $\delta\left(  x,y\right)  \geqslant d_{CM}\left(  x,y\right)  $
for all $x,y\in G_{CM}.$ The following theorem describes the behavior of
$\delta.$

\begin{thm}
\label{t.12.2}If $\left\{  \omega\left(  A,B\right)  :A,B\in H\right\}  $ is a
total subset of $\mathbf{C},$ then there exists $c\in\left(  0,1\right)  $
such that
\begin{equation}
c\left(  \left\Vert A\right\Vert _{H}+\sqrt{\left\Vert a\right\Vert
_{\mathbf{C}}}\right)  \leqslant\delta\left(  \mathbf{e},\left(  A,a\right)
\right)  \leqslant c^{-1}\left(  \left\Vert A\right\Vert _{H}+\sqrt{\left\Vert
a\right\Vert _{\mathbf{C}}}\right)  \text{ for all }\left(  A,a\right)
\in\mathfrak{g}_{CM}. \label{e.12.1}%
\end{equation}

\end{thm}

\begin{proof}
Our proof will be modeled on the standard proof of this result in the finite
dimensional context, see for example \cite{VSCC192,Montgomery2002}. The only
thing we must be careful of is to avoid using any compactness arguments.

For any left invariant metric, $d,$ (e.g. $d=\delta$ or $d=d_{CM})$ on
$G_{CM}$ we have
\begin{equation}
d\left(  \mathbf{e},xy\right)  \leqslant d\left(  \mathbf{e},x\right)
+d\left(  x,xy\right)  =d\left(  \mathbf{e},x\right)  +d\left(  \mathbf{e}%
,y\right)  \text{~}\forall~x,y\in G_{CM}. \label{e.12.2}%
\end{equation}

Given any path $g=\left(  w, c\right)  \in C_{CM}^{1}$ joining $\mathbf{e}$ to
$\left(  A, a\right) ,$ we have from Eq. (\ref{e.3.20}) that
\begin{align*}
\ell_{G_{CM}}\left(  g\right)   &  =\int_{0}^{1}\sqrt{\left\Vert w^{\prime
}(s)\right\Vert _{H}^{2}+\left\Vert c^{\prime}(s)-\omega(w(s),w^{\prime
}(s))/2\right\Vert _{\mathbf{C}}^{2}}ds\\
&  \geqslant\int_{0}^{1}\left\Vert w^{\prime}(s)\right\Vert _{H}%
ds\geqslant\left\Vert A\right\Vert _{H}%
\end{align*}
from which it follows that
\begin{equation}
\delta\left(  \mathbf{e},\left(  A,a\right)  \right)  \geqslant d_{CM}\left(
\mathbf{e},\left(  A,0\right)  \right)  \geqslant\left\Vert A\right\Vert _{H}.
\label{e.12.3}%
\end{equation}
Since the path $g\left(  t\right)  =\left(  tA,0\right)  $ is horizontal and
\[
\left\Vert A\right\Vert _{H}=\ell_{G_{CM}}\left(  g\right)  \geqslant
\delta\left(  \mathbf{e},\left(  A,0\right)  \right)  \geqslant d_{CM}\left(
\mathbf{e},\left(  A,0\right)  \right)  \geqslant\left\Vert A\right\Vert _{H}%
\]
it follows that
\begin{equation}
\delta\left(  \mathbf{e},\left(  A,0\right)  \right)  =d\left(  \mathbf{e}%
,\left(  A,0\right)  \right)  =\left\Vert A\right\Vert _{H} \text{ for all }
A\in H. \label{e.12.4}%
\end{equation}

Given $A,B\in H,$ let $\xi\left(  t\right)  =A\cos t+B\sin t$ for $0\leqslant
t\leqslant2\pi$ and
\[
g\left(  t\right)  =\left(  \xi\left(  t\right)  -A,\frac{1}{2}\int_{0}%
^{t}\omega\left(  \xi\left(  \tau\right)  -A,\dot{\xi}\left(  \tau\right)
\right)  d\tau\right)
\]
so that $l_{g\left(  t\right)  _{\ast}^{-1}}\dot{g}\left(  t\right)  =\left(
\xi\left(  t\right)  ,0\right)  ,$ $g\left(  0\right)  =\mathbf{e},$ and
\begin{align*}
g\left(  2\pi\right)   &  =\left(  0,\frac{1}{2}\int_{0}^{2\pi}\omega\left(
\xi\left(  \tau\right)  ,\dot{\xi}\left(  \tau\right)  \right)  d\tau\right)
\\
&  =\left(  0,\frac{1}{2}\int_{0}^{2\pi}\omega\left(  A,B\right)
d\tau\right)  =\left(  0,\pi\omega\left(  A,B\right)  \right)  .
\end{align*}
From this one horizontal curve we may conclude that
\begin{align}
\delta\left(  \mathbf{e},\left(  0,\pi\omega\left(  A,B\right)  \right)
\right)   &  \leqslant\ell_{G_{CM}}\left(  g\right)  =\int_{0}^{2\pi
}\left\Vert -A\sin t+B\cos t\right\Vert _{H}dt\nonumber\\
&  \leqslant2\pi\left(  \left\Vert A\right\Vert _{H}+\left\Vert B\right\Vert
_{H}\right)  . \label{e.12.5}%
\end{align}

Choose $\left\{  A_{\ell},B_{\ell}\right\}  _{\ell=1}^{d}\subset H$ such that
$\left\{  \pi\omega\left(  A_{\ell},B_{\ell}\right)  \right\}  _{\ell=1}^{d}$
is a basis for $\mathbf{C.}$ Let $\left\{  \varepsilon^{\ell}\right\}
_{\ell=1}^{d}$ be the corresponding dual basis. Hence for any $a\in\mathbf{C}$
we have%
\begin{align*}
\delta\left(  \mathbf{e},\left(  0,a\right)  \right)   &  =\delta\left(
\mathbf{e},\prod_{\ell=1}^{d}\left(  0,\varepsilon^{\ell}\left(  a\right)
\pi\omega\left(  A_{\ell},B_{\ell}\right)  \right)  \right) \\
&  \leqslant\sum_{\ell=1}^{d}\delta\left(  \mathbf{e},\left(  0,\varepsilon
^{\ell}\left(  a\right)  \pi\omega\left(  A_{\ell},B_{\ell}\right)  \right)
\right) \\
&  =\sum_{\ell=1}^{d}\delta\left(  \mathbf{e},\left(  0,\pi\omega\left(
\mathrm{sgn}(\varepsilon^{\ell}\left(  a\right)  )\sqrt{\left\vert
\varepsilon^{\ell}\left(  a\right)  \right\vert }A_{\ell},\sqrt{\left\vert
\varepsilon^{\ell}\left(  a\right)  \right\vert }B_{\ell}\right)  \right)
\right) \\
&  \leqslant2\pi\sum_{\ell=1}^{d}\left(  \left\Vert \sqrt{\left\vert
\varepsilon^{\ell}\left(  a\right)  \right\vert }A_{\ell}\right\Vert
_{H}+\left\Vert \sqrt{\left\vert \varepsilon^{\ell}\left(  a\right)
\right\vert }B_{\ell}\right\Vert _{H}\right)  ,
\end{align*}
wherein we have used Eq. (\ref{e.12.2}) for the first inequality and Eq.
(\ref{e.12.5}) for the second inequality. It now follows by simple estimates
that
\begin{equation}
\delta\left(  \mathbf{e},\left(  0,a\right)  \right)  \leqslant C_{1}%
\sum_{\ell=1}^{d}\sqrt{\left\vert \varepsilon^{\ell}\left(  a\right)
\right\vert }\leqslant C_{2}\sqrt{\sum_{\ell=1}^{d}\left\vert \varepsilon
^{\ell}\left(  a\right)  \right\vert }\leqslant C\left(  \omega\right)
\sqrt{\left\Vert a\right\Vert _{\mathbf{C}}}. \label{e.12.6}%
\end{equation}
for some constants $C_{1}\leqslant C_{2}\leqslant C\left(  \omega\right)
<\infty.$ Combining Eqs. (\ref{e.12.2}), (\ref{e.12.4}), and (\ref{e.12.6})
gives,%
\begin{align}
\delta\left(  \mathbf{e},\left(  A,a\right)  \right)   &  =\delta\left(
\mathbf{e},\left(  A,0\right)  \left(  0,a\right)  \right) \nonumber\\
&  \leqslant\delta\left(  \mathbf{e},\left(  A,0\right)  \right)
+\delta\left(  \mathbf{e},\left(  0,a\right)  \right)  \leqslant\left\Vert
A\right\Vert _{H}+C\left(  \omega\right)  \sqrt{\left\Vert a\right\Vert
_{\mathbf{C}}}. \label{e.12.7}%
\end{align}

To prove the analogous lower bound we will make use of the dilation
homomorphisms defined for each $\lambda>0$ by $\varphi_{\lambda}\left(
w,c\right)  =\left(  \lambda w,\lambda^{2}c\right)  $ for all $\left(
w,c\right)  \in\mathfrak{g}_{CM}=G_{CM}.$ One easily verifies that
$\varphi_{\lambda}$ is both a Lie algebra homomorphism on $\mathfrak{g}_{CM}$
and a group homomorphism on $G_{CM}.$ Using the homomorphism property it
follows that
\[
l_{\varphi_{\lambda}\left(  g\left(  t\right)  \right)  _{\ast}^{-1}}\frac
{d}{dt}\varphi_{\lambda}\left(  g\left(  t\right)  \right)  =\varphi_{\lambda
}\left(  l_{g\left(  t\right)  _{\ast}^{-1}}\dot{g}\left(  t\right)  \right)
\]
and consequently; if $g$ is any horizontal curve, then $\varphi_{\lambda}\circ
g$ is again horizontal and $\ell_{G_{CM}}\left(  \varphi_{\lambda}\circ
g\right)  =\lambda\ell_{G_{CM}}\left(  g\right)  .$ From these observations we
may conclude that
\begin{equation}
\delta\left(  \varphi_{\lambda}\left(  x\right)  ,\varphi_{\lambda}\left(
y\right)  \right)  =\lambda\delta\left(  x,y\right)  \text{ for all } x,y\in
G_{CM}. \label{e.12.8}%
\end{equation}

By Proposition \ref{p.3.10}, we know there exists $\varepsilon>0$ and
$K<\infty$ such that%
\begin{equation}
K\delta\left(  \mathbf{e},x\right)  \geq Kd_{G_{CM}}\left(  \mathbf{e}%
,x\right)  \geq\left\Vert x\right\Vert _{\mathfrak{g}_{CM}}\text{ whenever
}\left\Vert x\right\Vert _{\mathfrak{g}_{CM}}\leq\varepsilon. \label{e.12.9}%
\end{equation}
For arbitrary $x=\left(  A,a\right)  \in G_{CM},$ choose $\lambda>0$ such
that
\[
\varepsilon^{2}=\left\Vert \varphi_{\lambda}\left(  x\right)  \right\Vert
^{2}=\lambda^{2}\left\Vert A\right\Vert _{H}^{2}+\lambda^{4}\left\Vert
a\right\Vert _{\mathbf{C}}^{2},
\]
i.e.
\[
\lambda^{2}=\frac{\sqrt{\left\Vert A\right\Vert _{H}^{4}+4\left\Vert
a\right\Vert _{\mathbf{C}}^{2}\varepsilon^{2}}-\left\Vert A\right\Vert
_{H}^{2}}{2\left\Vert a\right\Vert _{\mathbf{C}}^{2}}.
\]
It then follows from Eqs. (\ref{e.12.8}) and (\ref{e.12.9}) that $\lambda
K\delta\left(  \mathbf{e},x\right)  =K\delta\left(  \mathbf{e},\varphi
_{\lambda}\left(  x\right)  \right)  \geq\varepsilon,$ i.e.%
\begin{align}
\delta^{2}\left(  \mathbf{e},x\right)   &  \geq\frac{\varepsilon^{2}}%
{K^{2}\lambda^{2}}=2\frac{\varepsilon^{2}}{K^{2}}\frac{\left\Vert a\right\Vert
_{\mathbf{C}}^{2}}{\sqrt{\left\Vert A\right\Vert _{H}^{4}+4\left\Vert
a\right\Vert _{\mathbf{C}}^{2}\varepsilon^{2}}-\left\Vert A\right\Vert
_{H}^{2}}\nonumber\\
&  =2\frac{\varepsilon^{2}\left\Vert a\right\Vert _{\mathbf{C}}^{2}}%
{K^{2}\left\Vert A\right\Vert _{H}^{2}}\frac{1}{\sqrt{1+\frac{4\left\Vert
a\right\Vert _{\mathbf{C}}^{2}\varepsilon^{2}}{\left\Vert A\right\Vert
_{H}^{4}}}-1}. \label{e.12.10}%
\end{align}
Since $\sqrt{1+x}-1\leq\min\left(  x/2,\sqrt{x}\right)  $ we have
\[
\frac{1}{\sqrt{1+x}-1}\geq\max\left(  \frac{2}{x},\frac{1}{\sqrt{x}}\right)
\geq\frac{1}{x}+\frac{1}{2\sqrt{x}}.
\]
Using this estimate with $x=4\left\Vert a\right\Vert _{\mathbf{C}}%
^{2}\left\Vert A\right\Vert _{H}^{-4}\varepsilon^{2}$ in Eq. (\ref{e.12.10})
shows%
\[
\delta^{2}\left(  \mathbf{e},x\right)  \geq2\frac{\varepsilon^{2}\left\Vert
a\right\Vert _{\mathbf{C}}^{2}}{K^{2}\left\Vert A\right\Vert _{H}^{2}}\left(
\frac{\left\Vert A\right\Vert _{H}^{4}}{4\left\Vert a\right\Vert _{\mathbf{C}%
}^{2}\varepsilon^{2}}+\frac{\left\Vert A\right\Vert _{H}^{2}}{4\varepsilon
\left\Vert a\right\Vert _{\mathbf{C}}}\right)  =\frac{1}{2K^{2}}\left(
\left\Vert A\right\Vert _{H}^{2}+\varepsilon\left\Vert a\right\Vert
_{\mathbf{C}}\right)  ,
\]
which implies the lower bound in Eq. (\ref{e.12.1}).
\end{proof}

We are now ready to give the proof of Theorem \ref{t.3.12}

\subsection{Proof of Theorem \ref{t.3.12}.}

\begin{proof}
The first assertion in Eq. (\ref{e.3.24}) of Theorem \ref{t.3.12} follows from
Theorem \ref{t.12.2} and the previously observed fact that $d_{CM}%
\leqslant\delta.$ To prove Eq. (\ref{e.3.25}), let $\varepsilon_{0}%
<\varepsilon/2$ where $\varepsilon>0$ is as in Proposition \ref{p.3.10}. Then
according to that proposition, if $d_{CM}\left(  \mathbf{e},x\right)
\leqslant\varepsilon_{0}$ then $\left\Vert x\right\Vert _{\mathfrak{g}_{CM}%
}\leqslant Kd_{CM}\left(  \mathbf{e},x\right)  \leqslant K\varepsilon_{0}.$ So
if $x=\left(  A,a\right)  ,$ we have $\left\Vert A\right\Vert _{H}\leqslant
K\varepsilon_{0}$ and $\left\Vert a\right\Vert _{\mathbf{C}}\leqslant
K\varepsilon_{0}$ and hence by Theorem \ref{t.12.2}, $\delta\left(
\mathbf{e},x\right)  \leqslant c^{-1}\left(  K\varepsilon_{0}+\sqrt
{K\varepsilon_{0}}\right)  .$ This implies that%
\begin{equation}
M\left(  \varepsilon_{0}\right)  :=\sup\left\{  \delta\left(  \mathbf{e}%
,x\right)  :x~\ni~d_{CM}\left(  \mathbf{e},x\right)  \leq\varepsilon
_{0}\right\}  \leqslant c^{-1}\left(  K\varepsilon_{0}+\sqrt{K\varepsilon_{0}%
}\right)  <\infty. \label{e.12.11}%
\end{equation}

Now suppose that $x\in G_{CM}$ with $d_{CM}\left(  \mathbf{e},x\right)
\geqslant\varepsilon_{0}.$ Choose a curve, $g\in C_{CM}^{1}\ $such that
$g\left(  0\right)  =\mathbf{e},$ $g\left(  1\right)  =x,$ and $\ell_{G_{CM}%
}\left(  g\right)  \leqslant d_{CM}\left(  \mathbf{e},x\right)  +\varepsilon
_{0}/4.$ Also choose $\varepsilon_{1}\in(\varepsilon_{0}/2,\varepsilon_{0}]$
such that $\ell_{G_{CM}}\left(  g\right)  =n\varepsilon_{1}$ with
$n\in\mathbb{N}$ and let $0=t_{0}<t_{1}<t_{2}<\dots<t_{n}=1$ be a partition of
$\left[  0,1\right]  $ such that $\ell_{G_{CM}}\left(  g|_{\left[
t_{i-1},t_{i}\right]  }\right)  =\varepsilon_{1}$ for $i=1,2,\dots,n.$ If
$x_{i}:=g\left(  t_{i}\right)  $ for $i=0,\dots,n,$ then $\varepsilon
_{0}\geqslant\varepsilon_{1}=\ell_{G_{CM}}\left(  g|_{\left[  t_{i-1}%
,t_{i}\right]  }\right)  \geqslant d_{CM}\left(  x_{i-1},x_{i}\right)  $ and
therefore from Eq. (\ref{e.12.11}) and the left invariance of $d_{CM}$ and
$\delta$ we have $1\geqslant M\left(  \varepsilon_{0}\right)  ^{-1}%
\delta\left(  x_{i-1},x_{i}\right)  $ for $i=1,2,\dots,n.$ Hence we may
conclude that
\begin{align*}
2d_{CM}\left(  \mathbf{e},x\right)   &  \geqslant d_{CM}\left(  \mathbf{e}%
,x\right)  +\varepsilon_{0}/4\geqslant\ell_{G_{CM}}\left(  g\right)
=\varepsilon_{1}n\\
&  \geqslant\varepsilon_{1}\sum_{i=1}^{n}M\left(  \varepsilon_{0}\right)
^{-1}\delta\left(  x_{i-1},x_{i}\right)  \geqslant\frac{\varepsilon_{0}}%
{2}M\left(  \varepsilon_{0}\right)  ^{-1}\delta\left(  \mathbf{e},x\right)  .
\end{align*}
Combining this estimate with the lower bound in Eq. (\ref{e.12.1}) shows Eq.
(\ref{e.3.25}) holds for all $\varepsilon_{0}$ sufficiently small which is
enough to complete the proof.
\end{proof}

\def\cprime{$'$}
\providecommand{\bysame}{\leavevmode\hbox to3em{\hrulefill}\thinspace}
\providecommand{\MR}{\relax\ifhmode\unskip\space\fi MR }
% \MRhref is called by the amsart/book/proc definition of \MR.
\providecommand{\MRhref}[2]{%
  \href{http://www.ams.org/mathscinet-getitem?mr=#1}{#2}
}
\providecommand{\href}[2]{#2}

%\bibliographystyle{amsplain}
%\bibliography{heis}

\end{document}